# Ricci flow on Poincare' and Thurston's Geometrization Conjecture

A Dissertation, Presented

By

Hassan Jolany

To

The Graduate School

in Partial Fulfillment of the

Requirements

for the Degree of

Master of Philosophy
In

Mathematics

**Supervisor:M.R.Darafsheh**

University of Tehran

2010

# Contents

# INTRODUCTION

## Chapter 1

# Chapter 1

## Preliminaries

*We will sweep through the basics of Riemannian geometry in this chapter, with a focus on the concepts that will be important for the Ricci flow later. I will quickly review the basic notions of infinitesimal Riemannian geometry, and in particular defining the Riemann, Ricci, and scalar curvature of a Riemannian manifold. This is a review only, in particular omitting any leisurely discussion of examples or motivation for Riemannian geometry. For more details see* [1-8]

### Connection

*A definition frequently used by differential geometries goes like this.*

*Let $\chi(M) = \{vector\ field\ on\ M\},\ \Gamma(E) = \{Smooth\ sections\ of\ E\}$*

*Definition 1.1. A connection on a vector bundle E is a map*

$$\nabla : \chi(M) \times \Gamma(E) \to \Gamma(E)$$

*Which satisfies the following axioms (where $\nabla_X \sigma = \nabla(X, \sigma)$)*

$$\nabla_X(f\sigma + \tau) = (Xf)\sigma + f\nabla_X\sigma + \nabla_X\tau$$

$$\nabla_{fX+Y}\sigma = f\nabla_X\sigma + \nabla_Y\sigma$$

*Here f is a real-valued function if E is a real vector bundle, A complex valued function if E is a complex vector bundle. Given a connection $\nabla$ in the sense of definition (1.1), we can define a map*

$$d_A : \Gamma(E) \to \Gamma(T^*M \otimes E) = \Gamma\big(Hom(TM, E)\big)$$

*by*

$$d_A(\sigma)(X) = \nabla_X\sigma$$

*Then $d_A$ satisfies .*

*Definition1.2. A connection on a vector bundle E is a map*

$$d_A : \Gamma(E) \to \Gamma(T^*M \otimes E)$$

*Which satisfies the following axiom*

$$d_A(f\sigma + \tau) = (df) \otimes \sigma + f d_A \sigma + d_A \tau$$

*Definition1.3 In more frequently used in gauge theory. The simplest example of a connection occurs on the bundle $E = M \times R^m$, the trivial vector bundle of rank m over M. A section of this bundle*

$$\sigma = \begin{pmatrix} \sigma^1 \\ \sigma^2 \\ \vdots \\ \sigma^m \end{pmatrix} : M \to R^m$$

*We can use the exterior derivative to define the "trivial" flat connection $d_A$ on :*

$$d_A \begin{pmatrix} \sigma^1 \\ \sigma^2 \\ \vdots \\ \sigma^m \end{pmatrix} = \begin{pmatrix} d\sigma^1 \\ d\sigma^2 \\ \vdots \\ d\sigma^m \end{pmatrix}$$

*More generally, given an $m \times m$ matrix*

$$\omega = \begin{pmatrix} \omega_1^1 & \cdots & \omega_m^1 \\ \vdots & \cdots \vdots \cdots & \vdots \\ \omega_1^m & \cdots & \omega_m^m \end{pmatrix}$$

*We can define a connection $d_A$ by*

$$d_A \begin{pmatrix} \sigma^1 \\ \sigma^2 \\ \vdots \\ \sigma^m \end{pmatrix} = \begin{pmatrix} d\sigma^1 \\ d\sigma^2 \\ \vdots \\ d\sigma^m \end{pmatrix} + \begin{pmatrix} \omega_1^1 & \cdots & \omega_m^1 \\ \vdots & \cdots \vdots \cdots & \vdots \\ \omega_1^m & \cdots & \omega_m^m \end{pmatrix} \begin{pmatrix} \sigma^1 \\ \sigma^2 \\ \vdots \\ \sigma^m \end{pmatrix}$$

*We can write this last equation in a more abbreviated fashion: $d_A \sigma = d\sigma + \omega\sigma$, Matrix multiplication being understand in the last term.*

*If $f \in \Gamma(V)$ (f is a section of bundle V) is such that $\nabla_X f = 0$ for all vector field , we say that say that f is parallel to the connection $\nabla$ .*

*A connection on the tangent bundle TM is known as an affine connection.*

*Remark1.1. Suppose $\phi_v$ be a parallel transport map i.e. $\phi_v : V_x \to V_{x+v}$ to each infinitesimal tangent vector $v \in V$ and for fixed x . We have*

$$\nabla_X f(x) = \lim_{t \to 0} \frac{\phi_{tX(x)}^{-1}\big(f(x+tX(x))\big) - f(x)}{t}$$

*Let $\nabla$ be a connection on $M$. We say that this connection is torsion-free if we have the pleasent identity*

$$\nabla_\alpha \nabla_\beta f = \nabla_\beta \nabla_\alpha f$$

*For all scalar field $f \in C^\infty(M)$.*

**Riemannian metrics and Levi-Civita connection**

*Definition 1.4. Let $M$ be an n-dimensional manifold. A Riemannian metric $g$ on $M$ is a smooth section of $T^*M \otimes T^*M$ defining a positive definite symmetric bilinear form on $T_p M$ for each $p \in M$. In local coordinates $(x^1, \ldots, x^n)$, one has a natural local basis $\{\partial_1, \ldots, \partial_n\}$ for $TM$, where $\partial_i = \frac{\partial}{\partial x_i}$. The metric tensor $g = g_{ij} dx^i \otimes dx^j$ is represented by a smooth matrix –valued function $g_{ij} = g(\partial_i, \partial_j)$. The pair $(M, g)$ is Riemannian manifold. We denote by $(g^{ij})$ the inverse of the matrix $(g_{ij})$*

*Definition 1.5. The Levi-Civita connection $\nabla: TM \times C^\infty(TM) \to C^\infty(TM)$ is the unique connection on $TM$ that is compatible with the metric and torsion free.*

$$X(g(Y,Z)) = g(\nabla_X Y, Z) + g(Y, \nabla_X Z)$$

$$\nabla_X Y - \nabla_Y X = [X, Y]$$

*Where $[X,Y]f \coloneqq X(Y(f)) - Y(X(f))$ defines the Lie bracket acting on functions.*

*From this one can show by taking a linear combination of the above equations with permutations of the vector field $X, Y$ and $Z$ that*

$$2g(\nabla_X Y, Z) = X\big(g(Y,Z)\big) + Y\big(g(X,Z)\big) - Z\big(g(X,Y)\big) + g([X,Y],Z) - g([X,Z],Y) - g([Y,Z],X) \qquad (1.1)$$

*Let $\{x^i\}_{i=1}^n$ be a local coordinate system defined in an open set $U$ in $M^n$. The christoffel symbols are defined in $U$ by $\nabla_{\frac{\partial}{\partial x^i}} \frac{\partial}{\partial x^j} := \Gamma_{ij}^k \frac{\partial}{\partial x^k}$. Here and through the thesis; we follow the Einstein summation convention of summing over repeated indices. By (1.1) and $\left[\frac{\partial}{\partial x^i}, \frac{\partial}{\partial x^j}\right] = 0$ we see that they are given by*

$$\Gamma_{ij}^k = \frac{1}{2} g^{kl} \left( \frac{\partial}{\partial x^i} g_{jl} + \frac{\partial}{\partial x^j} g_{il} - \frac{\partial}{\partial x^l} g_{ij} \right)$$

*Also for $= a^i \partial_i$, $X = b^j \partial_j$ we have a formula for $\nabla_X Y$ in local coordinates*

$$\nabla_X Y = b^i \partial_i(a^j)\partial_j + b^i a^j \Gamma_{ij}^k \partial_k$$

*Let $f$ be a smooth real-valued function on $M$. We define the Hessian of $f$, denoted $Hess(f)$, as follows:*

$$Hess(f)(X,Y) = X(Y(f)) - \nabla_X Y(f)$$

*Note that the Hessian is a contra variant, symmetric two-tensor i.e., for vector fields $X, Y$, we have*

$$Hess(f)(X,Y) = Hess(f)(Y,X)$$

*and*

$$Hess(f)(\phi X, \psi Y) = \phi\psi Hess(f)(X,Y)$$

*For all smooth function $\phi, \psi$. Other formulas for the Hessian are*

$$Hess(f)(X,Y) = \langle \nabla_X(\nabla f), Y \rangle = \nabla_X(\nabla_Y(f)) = \nabla^2 f(X,Y)$$

*Also, in local coordinate we have*

$$Hess(f)_{ij} = \partial_i \partial_j f - (\partial_k f)\Gamma_{ij}^k$$

*The Laplacian $\Delta f$ is defined as the Hessian. That is to say*

$$\Delta f(p) = \sum_{ij} g^{ij} Hess(f)(\partial_i, \partial_j)$$

*Also, if $\{X_i\}$ is an orthogonal basis for $T_P M$ then*

$$\Delta f(p) = \sum_i Hess(f)(X_i, X_i)$$

**Curvature tensor**

*Let $(M, g)$ be a Riemannian manifold and $\nabla$ the Riemannian connection, the curvature tensor is a $(1,3)$-tensor defined by*

$$R(X,Y)Z = \nabla^2_{X,Y} Z - \nabla^2_{Y,X} Z$$

$$= \nabla_X \nabla_Y Z - \nabla_Y \nabla_X - \nabla_{[X,Y]} Z = [\nabla_X, \nabla_Y]Z - \nabla_{[X,Y]} Z$$

*On vector fields $, Y, Z$.*

*Using metric g we can change this to a* $(0,4)-$ *tensor as follows*

$$R(X,Y,Z,W) = g(R(X,Y)Z,W)$$

*The Riemannian curvature tensor* $R(X,Y,Z,W)$ *satisfies the following properties*

1) *R is skew-symmetric in the first two and last two entries*

$$R(X,Y,Z,W) = -R(Y,X,Z,W) = R(Y,X,W,Z)$$

2) $R(X,Y,Z,W) = R(Z,W,X,Y)$
3) *R satisfies a cyclic permutation property called Bianchi's first :*
$$R(X,Y)Z + R(Z,X)Y + R(Y,Z)X = 0$$
4) $\nabla R$ *satisfies a cyclic permutation property called Bianchi's second identity*
$$(\nabla_Z R)(X,Y)W + (\nabla_X R)(Y,Z)W + (\nabla_Y R)(Z,X)W = 0$$

*In local coordinates, the Riemann curvature tensor can be represented as*

$$R(\partial_i, \partial_j, \partial_k, \partial_l) = R_{ijkl}$$

*Definition1.6. The sectional curvature of a 2-plance* $P \subset T_p M$ *is defined as* $K(P) = R(X,Y,X,Y)$ *where* $\{X,Y\}$ *is an orthogonal basis of P . We say that* $(M,g)$ *has positive sectional curvature for every 2-plane .*

*A Riemannian manifold is said to have constant sectional curvature if* $K(P)$ *is the same for all* $p \in M$ *and all 2-planes* $P \subseteq T_p M$.

*Definition1.7. The Riemannian manifold* $(M,g)$ *is said to be an Einstein manifold with Einstein constant* $\lambda$ *if* $icg = \lambda g$ .

*Definition1.8. Using the metric, one can replace the Riemannian curvature tensor R by a symmetric bilinear form Rm on* $\Omega^2 TM$ . *in local sections of* $\Omega^2 TM$ .*the formula for Rm is*

$$Rm(\varphi, \psi) = R_{ijkl}\varphi^{ij}\psi^{kl}$$

*We call Rm the curvature operator. We say* $(M,g)$ *has positive curvature operator if* $Rm(\varphi, \varphi) \geq 0$ *for any* $\varphi \in \Omega^2 TM$ .

*In local coordinate the Ricci curvature tensor is defined by*

$$Ric(X,Y) = g^{kl}R(X,\partial_k, Y, \partial_l)$$

*Cleary Ric is a symmetric bilinear form on TM,given in local coordinates by*

$$Ric = Ric_{ij}dx^i \otimes dx^j$$

*Where $Ric_{ij} = Ric(\partial_i, \partial_j)$. the scalar curvature is defined by*

$$R = R_g = tr_g Ric = g^{ij} Ric_{ij}$$

*Definition1.9 Let $\Delta$ denote the laplacian (or laplace-Beltrami operator) acting on functions, which is globally defined as the divergence of the gradient and given in local coordinates by:*

$$\Delta := div \nabla = g^{ij} \nabla_i \nabla_j = g^{ij} \left( \frac{\partial^2}{\partial x^i \partial x^j} - \Gamma_{ij}^k \frac{\partial}{\partial x^k} \right)$$

*There are other equivalent ways to define $\Delta$ such as*

$$\Delta f = \sum_{a=1}^{n} e_a(e_a f) - (\nabla_{e_a} e_a) f$$

*Where $\{e_a\}_{a=1}^n$ is an orthogonal frame.*

**Integration by parts**

*A basic tool is integration by parts recall that stokes theorem says that*

*Theorem1.1. If $\alpha$ is an $(n-1)$-form on a compact differentiable manifold $M^n$ with (possibly empty) boundary $\partial M$, then*

$$\int_M d\alpha = \int_{\partial M} \alpha$$

*The divergence theorem says*

*Theorem1.2. Let $(M, g)$ be a compact Riemannian manifold. If $X$ is a 1-form, then*

$$\int_{M^n} div(X) d\mu = \int_{\partial M^n} \langle X, v \rangle d\sigma$$

*Here $v$ is the unit outward normal, $d\mu$ denotes the volume form of $g$ and $d\sigma$ is the volume form of the boundary $\partial M^n$. Also we have following properties*

1) *On a closed manifold, $\int_{M^n} \Delta u\, d\mu = 0$*
2) *(Green) on a compact manifold*

$$\int_{M^n} (u \Delta v) d\mu = \int_{M^n} (v \Delta u) d\sigma$$

3) *If $f$ is a function and $X$ is a 1-form, then*

$$\int_{M^n} f \, div(X) d\mu = - \int_{M^n} \langle \nabla f, X \rangle d\mu + \int_{\partial M^n} f \langle X, v \rangle d\sigma$$

**Geodesics**

*We can now define the acceleration of a curve $\gamma: I \to M$ by the formula $\ddot{r} = \frac{d^2\gamma}{dt^2}$, in local coordinate this becomes*

$$\ddot{r} = \frac{d^2\gamma^k}{dt^2} \partial_k + \frac{d\gamma^i}{dt} \frac{d\gamma^j}{dt} \Gamma_{ij}^k \partial_k$$

*A $C^\infty$ curve $\gamma: I \to M$ is called a geodesic if $\ddot{r} = 0$. this means a smooth curve $\gamma: I \to M$ is called a geodesic if $\nabla_{\dot\gamma} \dot\gamma = 0$.*

**Projective linear groups**

*The group $PSL_2(R)$ is the group of automorphisms of the hyperbolic plane H and the elements of this group are transformation of the form*

$$z \to \frac{az + b}{cz + d}$$

*$a, b, c, d \in R$ and $ad - bc = 1$*

*This group also be made into a topological space by identifying the transformation with the point $(a, b, c, d) \in R^4$ in the subspace $\{(a, b, c, d) \in R^4 : ad - bc = 1\}$ and in fact one show that $PSL_2(R)$ is homeomorphic to $R^2 \times S^1$ and thus being a 3-dimensional manifold, moreover the group multiplication and taking of inverses are continuous with the topology on $PSL_2(R)$ and so this is a topological Group.*

## Lens space

*In the 3-manifold case, a picturesque description of a lens space is that of a space resulting from gluing two solid tori together by a homeomorphism of their boundaries. Of course, to be consistent, we should exclude the 3-sphere and $S^2 \times S^1$, both of which can be obtained as just described.*

There is a complete classification of three-dimensional lens spaces.

Three-dimensional lens spaces arise as quotients of $S^3 \subset \mathbb{C}^2$ by the action of the group that is generated by elements of the form $\begin{pmatrix} w & 0 \\ 0 & w^q \end{pmatrix}$

## Seifert fiber space

*A Seifert manifold is a closed 3-manifold together with a decomposition into a disjoint union of circles (called fibers) such that each fiber has a tubular neighborhood that forms a standard fibered torus.*

*A standard fibered torus corresponding to a pair of coprime integers $(a, b)$ with $a > b$ is the surface bundle of the automorphism of a disk given by rotation by an angle of $2\pi b/a$ (with the natural fibering by circles). If $a = 1$ the middle fiber is called ordinary, while if $a > 1$ the middle fiber is called exceptional. A compact Seifert fiber space has only a finite number of exceptional fibers.*

RICCI FLOW

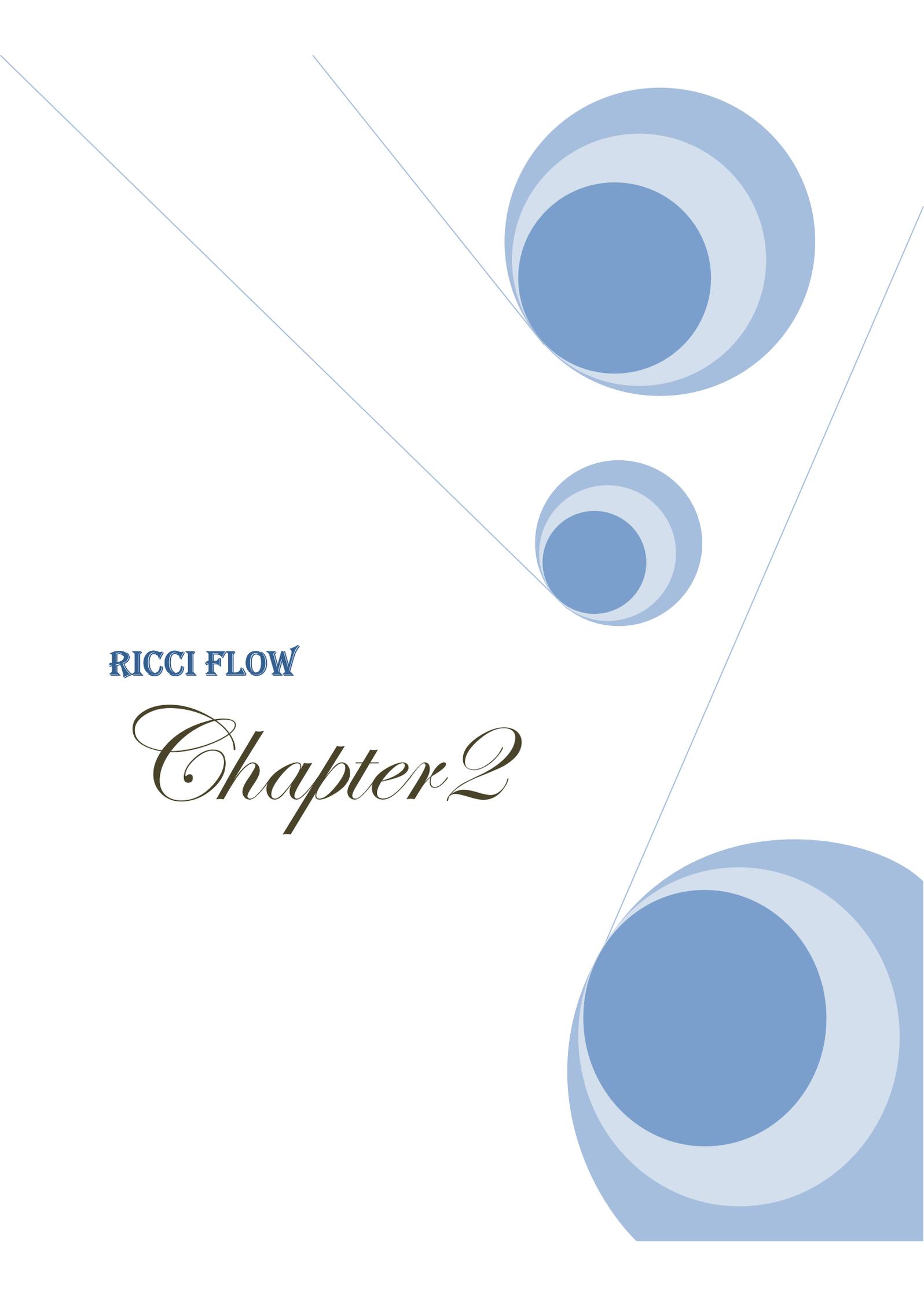

Chapter 2

# Chapter 2

## *Ricci flow*

### Flows on Riemannian manifold

*At first we introduce flows $t \to (M(t), g(t))$ on Riemannian manifolds $(M, g)$, which are recipes for describing smooth deformations of such manifolds over time, and derive the basic first variation formulae for how various structures on such manifolds change by such flows.*

*In this section we get a one-parameter family of such manifolds $t \to (M(t), g(t))$, parametrized by a "time" parameter $t$ in the usual manner the time derivatives $\dot{g}(t) = \frac{d}{d\tau} g(t)$ is*

$$\frac{d}{d\tau} g(t) = \lim_{dt \to 0} \frac{g(t + dt) - g(t)}{dt}$$

*and also the analogue of the time derivative $\frac{d}{d\tau} g(t)$ is then the Lie derivative $\mathcal{L}_{\partial_t} g$.*

*Definition 2.1.(Ricci flow) A one-parameter family of metrics $g(t)$ on a smooth manifold M for all time t in an interval I is said to obey Ricci flow if we have*

$$\frac{d}{d\tau} g(t) = -2 Ric(t)$$

*Note that this equation makes tensorial sense since g and Ric are both symmetric rank 2 tensor. The factor of 2 here is just a notational convenience and is not terribly important, but the minus sign − is crucial.*

*Let us give a quick introduction of what the Ricci flow equation means .In the harmonic coordinate $(x^1, x^2, \ldots, x^n)$ about p, that is to say coordinates where $\Delta x^i = 0$ for all i ,we have*

$$Ric_{ij} = Ric\left(\frac{\partial}{\partial x^i}, \frac{\partial}{\partial x^j}\right) = -\frac{1}{2}\Delta g_{ij} + Q_{ij}(g^{-1}, \partial g)$$

*Where Q is a quadratic form in $g^{-1}$ and $\partial g$ .See lemma (3.32) on page (92) of* [9]

*So in this coordinates, the Ricci flow equation is actually a heat equation for the Riemannian metric*

$$\frac{\partial g}{\partial t} = \Delta g + 2Q(g^{-1}, \partial g)$$

**Special solutions of flows**

*The Ricci flow $\frac{\partial g}{\partial t} = -Ric$ introduced by Hamilton is a degenerate parabolic evolution system on metrics.*

*Theorem2.1 (Hamilton* [9]*) Let $(M, g_{ij}(x))$ be a compact Riemannian manifold .Then there exists a constant $T > 0$ Such that the Ricci flow $\frac{\partial g}{\partial t} = -2Ric$, with $g_{ij}(x, 0) = g_{ij}(x)$ ,admits a unique smooth solution $g_{ij}(x, t)$ for all $x \in M$ and $t \in [0, T)$*

*Example 1(Einstein metric)*

*Recall that a Riemannian metric $g_{ij}$ is Einstein if $R_{ij} = \lambda g_{ij}$ for some for some constant $\lambda$ .Now since the initial metric is Einstein, we have*

$$R_{ij}(x, 0) = \lambda g_{ij}(x, 0)$$

and some $\lambda > 0$ .Let $g_{ij}(x, t) = \rho^2(t) g_{ij}(x, 0)$.

*Now from the definition of the Ricci tensor, one sees that*

$$R_{ij}(x, t) = R_{ij}(x, 0) = \lambda g_{ij}(x, 0)$$

*Thus the equation $\frac{\partial g_{ij}}{\partial t} = -2Ric$ corresponds to*

$$\frac{\partial\left(\rho^2(t) g_{ij}(x, 0)\right)}{\partial t} = -2\lambda g_{ij}(x, 0)$$

*This gives the ODE $\frac{\partial \rho}{\partial t} = -\frac{\lambda}{\rho}$.*

*Whose solution is given by $\rho^2(t) = 1 - 2\lambda t$ .Thus the evolving metric $g_{ij}(x,t)$ shrinks homothetically to a point as $\to T = \frac{1}{2\lambda}$ .*

*By contrast, if the initial metric is an Einstein metric of negative scalar curvature ,The metric will expand homothetically for all times .Indeed if $R_{ij}(x,0) = -\lambda g_{ij}(x,0)$ with $\lambda > 0$ and $g_{ij}(x,t) = \rho^2(t)g_{ij}(x,0)$ ,the $\rho(t)$ satisfies the ODE $\frac{d\rho}{dt} = \frac{\lambda}{\rho}$ . With the solution $\rho^2(t) = 1 + 2\lambda t$ , hence the evolving metric $g_{ij}(x,t) = \rho^2(t)g_{ij}(x,0)$ exists and expands homothetically formal times, and the curvature fall back to zero like $\frac{-1}{t}$ .*

*lemma2.1 (see [2]).Let $X \in T_p M$ be a unit vector. Suppose that X is contained in some orthonormal basis for $T_p M$. $Rc(X,X)$ is then the sum of the sectional curvatures of planes spanned by X and other elements of the basis.*

*Example 2 on an n-dimensional sphere of radius r (where $n > 1$), the metric is given by $g = r^2 \bar{g}$ where $\bar{g}$ is the metric on the unit sphere.The sectional curvatures are all $\frac{1}{r^2}$ . Thus for any unit vector v,the result of lemma (2.1) tells us that $Rc(v,v) = \frac{(n-1)}{r^2}$.Therefore $Rc = \frac{(n-1)}{r^2} g = (n-1)\bar{g}$*

*So the Ricci flow equation becomes an ODE*

$$\frac{\partial}{\partial t} g = -2Rc$$

$$\Rightarrow \frac{\partial}{\partial t}(r^2 \bar{g}) = -2(n-1)\bar{g}$$

$$\Rightarrow \frac{\partial}{\partial t}(r^2) = -2(n-1)$$

*We have the solution $(t) = \sqrt{R_0^2 - 2(n-1)t}$ , where $R_0$ is the initial radius of the sphere . The manifold shrinks to a point as $\to \frac{R_0^2}{2(n-1)}$ .*

*Similarly, for hyperbolic n-space $H^n$ (where $n > 1$),The Ricci flow reduces to the ODE*

$$\frac{\partial}{\partial t}(r^2) = 2(n-1)$$

*Which has the solution $r(t) = \sqrt{R_0^2 + 2(n-1)t}$ .So the solution expands out to infinity .*

*Example3 (Ricci Soliton) A Ricci soliton is a Ricci flow$(M, g(t))$ , $0 \leq t < T \leq \infty$ ,with the property that for each $t \in [0, T)$ there is a diffeomorphism $\varphi_t: M \to M$ and a constant $\sigma(t)$ such that $\sigma(t)\varphi_t^* g(0) = g(t)$. That is to say , in a Ricci soliton all the Riemannian manifold $(M, g(t))$ are isomorphic up to a scale factor that is allowed to vary with t .The soliton is said to be shrinking if $\varphi'(t) < 0$ for all t (note that$\varphi_0 = id$ and $\sigma(0) = 1$ )*

*Taking the derivative of the equation $g(t) = \sigma(t)\varphi_t^* g(0)$ and evaluating at $t = 0$ yields*

$$\frac{\partial}{\partial t} g(t) = \frac{\partial \sigma(t)}{\partial t} \varphi_t^* g(0) + \sigma(t) \frac{\partial}{\partial t} \varphi_t^* g(0)$$

$$-2Rc\big(g(0)\big) = \sigma'(0) g(0) + \mathcal{L}_v g(0)$$

*Where $v = \frac{d\varphi_t}{dt}$ .Let us set $\sigma'(0) = 2\lambda$. Now because on a Riemannian manifold $(M, g)$,we have $(\mathcal{L}_X g)_{ij} = \nabla_i X_j + \nabla_j X_i$*

*We get $-2R_{ij} = 2\lambda g_{ij} + \nabla_i X_j + \nabla_j X_i$*

*As a special case we can consider the case that $v$ is the gradient vector field of some scalar function f on $M^n$, i.e. $v_i = \nabla_i f$, the equation then become $R_{ij} + \lambda g_{ij} + \nabla_i \nabla_j f = 0$ such solutions are known as gradient Ricci solutions .*

*The gradient Ricci solitons play a role in motivating the definition Perelman's  $\mathcal{F}$ and $\mathcal{W}$ −functionals.*

*Proposition2.1(see[10]) suppose we have a complete Riemannian manifold $(M, g(0))$ ,a smooth function $f: M \to R$ ,and  a constant $\lambda > 0$ such that $-Ric\big(g(0)\big) = Hess(f) - \lambda g(0)$ then there is $T > 0$ and a gradient shrinking soliton $(M, g(t))$ defined for $0 \leq t < T$.*

*In 1989, Shi generalized the theorem of Hamilton about short-time existence and uniqueness theorem for the Ricci flow to complete non-compact manifolds with bounded curvature .*

*Theorem2.2(Shi[11]) Let $(M, g_{ij}(x))$ be a complete non compact Riemannian manifold of dimension n with bounded curvature. Then there exists a constant $T > 0$ such that the initial value problem*

$$\begin{cases} \frac{\partial}{\partial t} g_{ij}(x,t) = -2R_{ij}(x,t) \text{ on } M \\ g_{ij}(x,0) = g_{ij}(x) \text{ on } M \end{cases}$$

*Admits a smooth solution $g_{ij}(x,t)$, $t \in [0,T]$, with bounded curvature*

*Recently Chen and Zhu proved the following uniqueness Theorem.*

*Theorem2.3 (Chen-Zhu [11]) Let $(M, g_{ij}(x))$ be a complete non-compact Riemannian manifold of dimension n with bounded curvature. Let $g_{ij}(x,t)$ and $\overline{g_{ij}}(x,t)$ be two solutions to the Ricci flow on $M \times [0,T]$ with $g_{ij}(x)$ as the initial data and with bounded curvatures. Then Let $g_{ij}(x,t) = \overline{g_{ij}}(x,t)$ for $(x,t) \in M \times [0,T]$.*

**DELATIONS**
*Now we specialize to some specific flows $(M, g(t))$ of a Riemannian metric on a fixed backward manifold M. The simplest such flow (besides the trivial flow $g(t) = g(0)$, of course) is that of a dilation*

$$g(t) := A(t)g(0)$$

*Where $A(t) > 0$ is a positive scalar with $A(0) = 1$. The flow here is given by*

$$\dot{g}(t) = a(t)g(t) \qquad (2.1)$$

*where $a(t) := \frac{\dot{A}(t)}{A(t)} = \frac{d}{dt} \log A(t)$ is the logarithmic derivative of A (or equivalently, $A(t) = \exp\left(\int_0^t a(t')dt'\right)$)*

*In this case our variation formulas become very simple*

$$\frac{d}{dt} g^{\alpha\beta} = -a g^{\alpha\beta}$$

$$\dot{\Gamma}^{\gamma}_{\alpha\beta} = 0$$

$$\dot{Ric}_{\alpha\beta} = 0$$

$$\dot{R} = -aR$$

$$\frac{d}{dt}d(x,y) = \frac{1}{2}ad(x,y)$$

$$\frac{d}{dt}d\mu = \frac{d}{2}ad\mu$$

**Modifying Ricci flow**

*If $g(t)$ solves Ricci flow and we set $\tilde{g}(s) = A(s)g(t(s))$ for some reparameterized time $s = s(t)$ and some scalar $A = A(s) > 0$, Then the Ricci curvature here is $\widetilde{Ric}(s) = Ric(t(s))$. We then see from the chain rule that $\tilde{g}$ obeys the equation*

$$\frac{d}{ds}\tilde{g}(s) = -2A(s)\frac{dt}{ds}\widetilde{Ric}(s) + a(s)\tilde{g}(s)$$

*Where $a$ is the logarithmic derivative of $A$. If we normalize the time reparameterization by requiring $\frac{dt}{ds} = \frac{1}{A(s)}$, we thus see that $\tilde{g}$ obeys normalized Ricci flow*

$$\frac{d}{ds}\tilde{g} = -2\widetilde{Ric} + a\tilde{g}(s)$$

*Which can be viewed as a combination of (2.1) and $\frac{d}{dt}g(t) = -2Ric(t)$.*

*Definition 2.2 (Bounded curvature solution) we say that a solution $g(t)$, $t \in I$ of the Ricci flow has bounded curvature if on every compact subinterval $[a,b] \subset I$ the Riemann curvature tensor is bounded. In particular we do not assume the curvature bound is uniform in time on non-compact time intervals.*

### Time evolving metrics:
*A normal coordinates about the point p are defined by*

a) $\gamma_v(t) = (tv^1, tv^2, \ldots, tv^n)$ is a geodesic

b) $g_{ij}(p) = \delta_{ij}$
c) $\Gamma_{ij}^k(p) = 0, \partial_i g_{jk}(p) = 0$

*Suppose that $g_{ij}(t)$ is a time-dependent Riemannian metric and*

$$\frac{d}{dt}g_{ij}(t) = h_{ij}(t)$$

*Then the various geometric quantities evolve according to the following equations:*

1) Metric inverse

$$\frac{\partial}{\partial t}g^{ij}(t) = -h^{ij} = -g^{ik}g^{jl}h_{kl}$$

2) Christoffel symbols

$$\frac{\partial}{\partial t}\Gamma_{ij}^k = \frac{1}{2}g^{kl}(\nabla_i h_{jl} + \nabla_j h_{il} - \nabla_l h_{ij})$$

3) Riemann curvature tensor

$$\frac{\partial}{\partial t}R_{ijk}^l = \frac{1}{2}g^{lp}\begin{Bmatrix}\nabla_i\nabla_j h_{kp} + \nabla_i\nabla_k h_{jp} - \nabla_i\nabla_p h_{jk} \\ -\nabla_j\nabla_i h_{kp} + \nabla_j\nabla_k h_{ip} - \nabla_j\nabla_p h_{ik}\end{Bmatrix}$$

4) Ricci tensor

$$\frac{\partial}{\partial t}R_{ij} = \frac{1}{2}g^{pq}(\nabla_p\nabla_i h_{jp} + \nabla_q\nabla_j h_{ip} - \nabla_q\nabla_p h_{ij} - \nabla_i\nabla_j h_{qp})$$

5) Scalar curvature

$$\frac{\partial}{\partial t}R = -\Delta H + \nabla^p\nabla^q h_{pq} - h^{pq}R_{pq}$$

Where $H = g^{pq}h_{pq}$

6) Volume element

$$\frac{\partial}{\partial t}d\mu = \frac{H}{2}d\mu$$

7) Volume of manifold

$$\frac{\partial}{\partial t}\int_M d\mu = \int_M \frac{H}{2}d\mu$$

8) Total scalar curvature on a closed manifold M

$$\frac{\partial}{\partial t}\int_M R d\mu = \int_M \left(\frac{1}{2}RH - h^{ij}R_{ij}\right)d\mu$$

*We only proof 2), 5) and 6)*

*Proof2) We know*

$$\frac{d}{dt}\Gamma_{ij}^k = \frac{1}{2}g^{kl}(\partial_i g_{jl} + \partial_j g_{il} - \partial_l g_{ij})$$

*So we get*

$$\partial_t \Gamma_{ij}^k = \frac{1}{2}(\partial_t g^{kl})(\partial_i g_{jl} + \partial_j g_{il} - \partial_l g_{ij}) + \frac{1}{2}g^{kl}(\partial_i \partial_t g_{jl} + \partial_j \partial_t g_{il} - \partial_l \partial_t g_{ij})$$

*Now we work in normal coordinates about a point p. So according to properties b)and c) of normal coordinates we get $\partial_i g_{jk} = 0$ at, $\partial_i A = \nabla_i A$ at p for any tensor A. Hence*

$$\partial_t \Gamma_{ij}^k(p) = \frac{1}{2}g^{kl}(\nabla_i h_{jl} + \nabla_j h_{il} - \nabla_l h_{ij})(p)$$

*Now although the Christoffel symbols are not the coordinates of a tensor quantity, their derivative is. (This is true because the difference between the christoffel symbols of two connections is a tensor .Thus, by taking a fixed point connection with christoffel symbols $\widetilde{\Gamma_{ij}^k}$, we have $\partial_t \Gamma_{ij}^k = \partial_t(\Gamma_{ij}^k - \widetilde{\Gamma_{ij}^k})$ and the Right hand side is clearly a tensor) Hence both side of this equation are the coordinates of tonsorial quantities, so it does not matter what coordinates we evaluate them. in particular ,the equation is true for any coordinates , not just normal coordinates ,and about any point p.*

*Proof 5) from 4) and 1) we get*

$$\partial_t R = \partial_t(g^{ij}R_{ij}) = \partial_t(g^{ij})R_{ij} + g^{ij}(\partial_t R_{ij})$$

$$= -h^{ij}R_{ij} + g^{ij}\left(\frac{1}{2}g^{pq}(\nabla_q \nabla_i h_{jp} + \nabla_q \nabla_j h_{ip} - \nabla_q \nabla_p h_{ij} - \nabla_i \nabla_j h_{qp})\right)$$

$$= -\Delta H + \nabla^p \nabla^q h_{pq} - h^{pq}R_{pq}$$

*(Note that $\nabla g = 0$ and $\Delta = g^{ij}\nabla_i \nabla_j$).*

*Proof 6) we know*

$$d\mu = \sqrt{det g_{ij}}\, dx^1 \wedge dx^2 \wedge dx^n$$

*And*

$$\frac{d}{dt}detA = (A^{-1})^{ij}\left(\frac{dA_{ij}}{dt}\right)detA$$

*Now by chain rule formula we obtain*

$$\partial_t d\mu = \partial_t \sqrt{det g_{ij}} dx^1 \wedge dx^2 \wedge dx^n = \frac{1}{2\sqrt{det g_{ij}}} g^{ij} h_{ij} det g \, dx^1 \wedge dx^2 \wedge dx^n = \frac{H}{2} d\mu$$

*Where* $H = g^{ij} h_{ij}$.

## Evolution equations for derivatives of curvature

*The Ricci flow is an evolution equation on the metric. The evolution equation* $\frac{\partial g_{ij}}{\partial t} = -2R_{ij}$ *for the metric implies a heat equation for the Riemannian Curvature* $R_{ijkl}$ *which we now derive.*

*Theorem 2.4 (Hamilton[10]) under the Ricci flow, the curvature tensor satisfies the evaluation equation*

$$\frac{\partial}{\partial t} R_{ijkl} = \Delta R_{ijkl} + 2(B_{ijkl} - B_{ijlk} - B_{iljk} + B_{ikjl})$$

$$-g^{pq}(R_{pjkl}R_{qi} + R_{ipkl}R_{qj} + R_{ijpl}R_{qk} + R_{ijkp}R_{ql})$$

*Where* $B_{ijkl} = g^{pr} g^{qs} R_{piqj} R_{rksl}$ *and* $\Delta$ *is the Laplacian with respect to the evolving metric*

*Proof choose* $\{x^1, x^2, \dots, x^m\}$ *to be a normal coordinate system at a fixed point. At this point, we compute*

$$\frac{\partial}{\partial t} \Gamma^h_{jl} = \frac{1}{2} g^{hm} \left( \nabla_j(-2R_{lm}) + \nabla_l(-2R_{jm}) - \nabla_m(-2R_{jl}) \right)$$

$$\frac{\partial}{\partial t} R^h_{ijl} = \frac{\partial}{\partial x^i}\left(\frac{\partial}{\partial t}\Gamma^h_{jl}\right) - \frac{\partial}{\partial x^j}\left(\frac{\partial}{\partial t}\Gamma^h_{il}\right)$$

$$\frac{\partial}{\partial t} R_{ijkl} = g_{hk} \frac{\partial}{\partial t} R^h_{ijl} + \frac{\partial g_{hk}}{\partial t} R^h_{ijl}$$

*Combining these identities we get*

$$\frac{\partial}{\partial t} R_{ijkl} = g_{hk} \left[ \left( \frac{1}{2} \nabla_i \left[ g^{hm} \left( \nabla_j(-2R_{lm}) + \nabla_l(-2R_{jm}) - \nabla_m(-2R_{jl}) \right) \right] \right) \right]$$

$$-\left( \frac{1}{2} \nabla_j \left[ g^{hm} \left( \nabla_i(-2R_{lm}) + \nabla_l(-2R_{im}) - \nabla_m(-2R_{il}) \right) \right] \right) - 2R_{hk} R^h_{ijl}$$

$$= \nabla_i \nabla_k R_{jl} - \nabla_i \nabla_l R_{jk} - \nabla_j \nabla_k R_{il} + \nabla_j \nabla_l R_{ik} - R_{ijlp} g^{pq} R_{qk} - R_{ijkp} g^{pq} R_{ql}$$
$$- 2R_{ijpl} g^{pq} R_{qk}$$
$$= \nabla_i \nabla_k R_{jl} - \nabla_i \nabla_l R_{jk} - \nabla_j \nabla_k R_{il} + \nabla_j \nabla_l R_{ik} - g^{pq} (R_{ijkp} R_{ql} + R_{ijpl} R_{qk}).$$

*Note that we have used of this formula*

$$\nabla_i \nabla_j v_k - \nabla_j \nabla_i v_k = R_{ijkl} g^{lm} v_m,$$

*Now according to Simon's identity in extrinsic geometry,*

$$\Delta R_{ijkl} + 2(B_{ijkl} - B_{ijlk} - B_{iljk} + B_{ikjl}) = \nabla_i \nabla_k R_{jl} - \nabla_i \nabla_l R_{jk} - \nabla_j \nabla_k R_{il} + \nabla_j \nabla_l R_{ik} + g^{pq} (R_{pjkl} R_{qi} + R_{ipkl} R_{qj})$$

*So we obtain*

$$\Delta R_{ijkl} = g^{pq} \nabla_p \nabla_i R_{qjkl} - g^{pq} \nabla_p \nabla_j R_{qikl}$$

$$= \nabla_i \nabla_k R_{jl} - \nabla_i \nabla_l R_{jk} - (B_{ijkl} - B_{ijlk} - B_{iljk} + B_{ikjl}) + g^{pq} R_{pjkl} R_{qi} - \nabla_j \nabla_k R_{il}$$
$$+ \nabla_j \nabla_l R_{ik} + (B_{jikl} - B_{jilk} - B_{jlik} + B_{jkil}) - g^{pq} R_{pikl} R_{qj}$$
$$+ g^{pq} (R_{pjkl} R_{qi} + R_{ipkl} R_{qj}) - 2(B_{ijkl} - B_{ijlk} - B_{iljk} + B_{ikjl})$$

*As desired, where in the last step we used the symmetries*

$$B_{ijkl} = B_{klij} = B_{jilk}$$

*So proof is complete.*

*Corollary 2.1. (see [11]) The Ricci curvature satisfies the evolution equation*

$$\frac{\partial}{\partial t} R_{ik} = \Delta R_{ik} + 2g^{pr} g^{qs} R_{piqk} R_{rs} - 2g^{pq} R_{pi} R_{qk}$$

*Lemma 2.2 (Hamilton[11]) The scalar curvature satisfies the evolution equation*

$$\frac{\partial R}{\partial t} = \Delta R + 2|Ric|^2$$

*Proof*

$$\frac{\partial R}{\partial t} = g^{ik}\frac{\partial R_{ik}}{\partial t} + \left(-g^{ip}\frac{\partial g_{pq}}{\partial t}g^{qk}\right)R_{ik}$$

$$= g^{ik}\left(\Delta R_{ik} + 2g^{pr}g^{qs}R_{piqk}R_{rs} - 2g^{pq}R_{pi}R_{qk}\right) + R_{pq}R_{ik}g^{ip}g^{qk}$$

$$= \Delta R + 2|Ric|^2$$

*Lemma2.3. The volume form $dvol(x,t)$ satisfy the following evolution equation under Ricci flow*

$$\frac{d}{dt}dvol(x,t) = -R(x,t)\,dvol(x,t)$$

*and if $f: M \to R$ is a time dependent function, then under the Ricci flow $\frac{\partial}{\partial t}\Delta f = \Delta\frac{\partial f}{\partial t} + 2\langle Ric, Hess(f)\rangle$*

## Harmonic map

*Let $:(M^n, g) \to (N^m, g)$, the map Laplacian of $f$ is defined by*

$$\left(\Delta_{g,h}f\right)^\gamma = \Delta_g(f^\gamma) + g^{ij}\left(\Gamma(h)^\gamma_{\alpha\beta}\circ f\right)\frac{\partial f^\alpha}{\partial x^i}\frac{\partial f^\beta}{\partial x^j}$$

$$= g^{ij}\left(\frac{\partial^2 f^\gamma}{\partial x^i \partial x^j} - \Gamma^k_{ij}\frac{\partial f^\gamma}{\partial x^k} + \left(\Gamma(h)^\gamma_{\alpha\beta}\circ f\right)\frac{\partial f^\alpha}{\partial x^i}\frac{\partial f^\beta}{\partial x^j}\right) \qquad (2.2)$$

*Where $f^\gamma := y^\gamma \circ f$ and $\{x^i\}$ and $\{y^\alpha\}$ are coordinates on M and N, respectively. Note that $\Delta_{g,h}f \in C^\infty(f^*TN)$, where $f^*(TN) \to M$ Is the pullback vector bundle of TN by f. In (2.2) $(\Delta_g f)^\gamma$ denotes the Laplacian with respect to $g$ of the function $f^\gamma$. As a special case, if $M = N$ and $f$ is the identity map and we choose the x and y coordinates to be the same, then*

$$\left(\Delta_{g,h}id\right)^k = g^{ij}\left(-\Gamma(g)^k_{ij} + \Gamma(h)^k_{ij}\right)$$

*The derivative df of a map $f: M^n \to N^m$ is a section of the vector bundle $= T^*M \otimes f^*(TN)$. On E is a national metric and compatible connection $\nabla^{g,h}$ defined by the (dual of the) Riemannian metric $g$ and associated Levi-Civita connection on $T^*M$ and the pullback by $f$ on the metric h and its associated Levi-Civita connection on TN.*

So $\nabla^{g,h} df$ is a section of the bundle $T^*M \otimes_S T^*M \otimes f^*(TN)$. The map Laplacian is the trace with respect to g of $\nabla^{g,h} df$

$$\Delta_{g,h} f = tr_g(\nabla^{g,h} df)$$

A map $f: (M^n, g) \to (N^m, h)$ is called a harmonic map if $\Delta_{g,h} f = 0$.

In the case where $N = \mathbb{R}$ a harmonic map is the same as a harmonic function.

## Ricci flow on almost flat manifolds:

A compact Riemannian manifold $M^n$ is called $\varepsilon$ −flat if, its curvature is bounded in terms of the diameter as follows

$$|K| \leq \varepsilon . d(M)^{-2},$$

Where K denotes the sectional curvature and $d(M)$ the diameter of M. Here we show summary result of Ricci flow that act on almost flat manifolds.

Theorem2.5. (Ricci flow on almost flat manifolds[12])

In any dimension n there exists an $\varepsilon(n) > 0$ such that for any $\varepsilon \leq \varepsilon(n)$ an $\varepsilon$ −flat Riemannian manifold $(M^n, g)$ has the following properties :

i)      The solution $g(t)$ to the Ricci flow
$$\frac{\partial g}{\partial t} = -2ric_g, \quad g(0) = g \qquad (2.3)$$
exists for all $t \in [0, \infty]$

ii)      Along the flow () one has
$$\lim_{t \to \infty} |K|_{g_t} . d^2(M, g_t) = 0$$

iii)      $g(t)$ converges to a flat metric along Ricci flow (2.3), if and only if the fundamental group of M is almost abelian (abelian up to a subgroup of finite index).

In the case n=2 the normalized Ricci flow looks like

$$\frac{\partial g}{\partial t} = (sc(g) - K)g \qquad (2.4)$$

Where $sc(g)$ is average sectional curvature

Uniformization Theorem2.6 (see [12])

*On a closed Riemannian manifold $(M^2, g)$ with $\chi(M)$ (Euler characteristic of $M^2$), The normalized Ricci flow (2.4) with $g(0) = g_0$ has a unique solution for all time, moreover, as $t \to \infty$, the metrics $g(t)$ converge uniformly in any $C^k$ −norm to a flat metric $g_\infty$*

## The normalized Ricci flow

In case of 3-manifolds with positive Ricci curvature and higher dimensional manifolds with positive curvature operator, works by Hamilton, Huisken, Bohm and wilking showed that the normalized Ricci flow will evolve the metric to one with constant curvature. One way to avoid collapse is to add the condition that volume by preserved along the evaluation. To preserve the volume we need to change a little bit the equation for the Ricci flow.

We know

$$d\mu = \sqrt{\det g_{ij}} \, dx^1 . dx^2 \ldots dx^n$$

*Let us define the average scalar curvature*

$$r = \frac{\int_M R \, d\mu}{\int_M d\mu}$$

*Let $\tilde{g}(t) = \psi(t) g(t)$ with $\psi(0) = 1$. Let us choose $\psi(t)$ so that the volume of the manifold with respect to $\tilde{g}$ is constant. So*

$$vol(\tilde{g}(t)) = vol(\tilde{g}(0))$$

*But we know if $\tilde{g} = cg$ then $d\tilde{\mu} = c^{\frac{n}{2}} d\mu$ therefore*

$$\psi(t)^{\frac{n}{2}} vol(g(t)) = vol(g(0))$$

*So we obtain*

$$\psi(t) = \left( \frac{\int_{M^n} d\mu(t)}{\int_{M^n} d\mu(0)} \right)^{\frac{2}{-n}}$$

*But according to relation $\frac{\partial}{\partial t} \int_{M^n} d\mu = - \int_{M^n} R \, d\mu$ we get.*

$$\frac{d}{dt} \psi(t) = \frac{2}{-n} \left( \frac{\int_{M^n} d\mu(t)}{\int_{M^n} d\mu(0)} \right)^{\frac{2}{-n} - 1} \frac{\frac{d}{dt} \int_{M^n} d\mu(t)}{\int_{M^n} d\mu(0)}$$

$$= \frac{2}{n} \frac{\psi(t)}{\left(\frac{\int_{M^n} d\mu(t)}{\int_{M^n} d\mu(0)}\right)} \frac{\int_{M^n} R d\mu(t)}{\int_{M^n} d\mu(0)} = \frac{2r}{n}\psi(t)$$

*The normalized average scalar curvature, $\tilde{r}$, is defined by*

$$\tilde{r} = \frac{\int_{M^n} \tilde{R} d\mu}{\int_{M^n} d\mu} = \frac{r}{\psi(t)} \quad (Because\ if\ \tilde{g} = cg\ Then\ \tilde{R} = c^{-1}R)$$

*Where $\tilde{R} \coloneqq R(\tilde{g})$. So*

$$\frac{d}{dt}\psi(t) = \frac{2\tilde{r}}{n}(\psi(t))^2$$

*So we get*

$$\frac{\partial}{\partial t}\tilde{g} = \left(\frac{\partial}{\partial t}g\right)\psi(t) + g\frac{d}{dt}\psi(t)$$

$$= -2\psi(t)Rc(g(t)) + \frac{2\tilde{r}}{n}(\psi(t))^2 g$$

$$= \psi(t)\left(-2Rc(\tilde{g}(t)) + \frac{2\tilde{r}}{n}\tilde{g}\right)$$

*We define a rescaling of time to get rid of the $\psi(t)$ terms in this evolution equation:*

$$\tau = \int_0^t \psi(u)du$$

*So $\frac{\partial}{\partial t}\tilde{g} = -2\widetilde{Rc} + \frac{2\tilde{r}}{n}\tilde{g}$, that $\tilde{g}(t)$ called the normalized Ricci flow.*

*Therefore we checked that $\tilde{g}(t)$ is a solution to the normalized Ricci flow. Note that $\frac{\partial}{\partial t}g_{ij} = -2R_{ij} + \frac{2}{n}rg_{ij}$ that called unnormalized Ricci flow, may not have solution even for short time [Hamilton].*

## Laplacian spectrum under Ricci flow

*Let $(M^n, g)$ be a Riemannian manifold, we know*

$$\Delta f = g^{ij}\nabla_i \nabla_j f$$

$$= g^{ij}\left(\frac{\partial^2 f}{\partial x_i \partial x_j} - \Gamma_{ij}^k \frac{\partial f}{\partial x_k}\right)$$

*By an eigenvalue $\lambda$ of $\Delta$ we mean there exists a non-zero function $f$ such that*

$$\Delta f + \lambda f = 0$$

*The set of all eigenvalues is called the spectrum of the operator it is well-known that the spectrum of the laplacian –Beltrami operator is discrete and non-negative .So we can write for the Laplacian Eigen values is the following ways*

$$0 = \lambda_0 < \lambda_1 \leq \lambda_2 \leq \lambda_3 \leq \cdots$$

*Note that if $M^n$ is closed ,Then assuming $f \not\equiv 0$ we have*

$$\lambda = \frac{\int_{M^n} |\nabla f|^2 d\mu}{\int_{M^n} f^2 d\mu}$$

*Here we assume that if we have a smoothly varying one-parameter family of metrics $(t)$ , each Laplacian eigenvalue $\lambda_\alpha(g(t))$ will also vary smoothly.*

*Lemma2.4 On a closed manifold we have*

$$\int_{M^n} Rc(\nabla f, \nabla f) d\mu \leq \frac{n-1}{n} \int_{M^n} (\Delta f)^2 d\mu$$

*Theorem2.7 (Lichrerowicz) suppose f is an eigen function of the Laplacian with Eigen value$\lambda$:*

$$\Delta f + \lambda f = 0$$

*If$Rc \geq (n-1)kg$, where $k > 0$ is a constant, then*

$$\lambda \geq nk$$

*Proof: By applying lemma 2.4, we have*

$$(n-1)k \int_{M^n} |\nabla f|^2 d\mu \leq \int_{M^n} Rc(\nabla f, \nabla f) d\mu \leq$$

$$\leq \frac{n-1}{n} \int_{M^n} (\Delta f)^2 d\mu = \frac{n-1}{n} \lambda^2 \int_{M^n} f^2 d\mu$$

*and because $\lambda = \frac{\int_{M^n} |\nabla f|^2 d\mu}{\int_{M^n} f^2 d\mu}$, so proof is complete.*

*Theorem2.8. Let $M^n$ be a closed manifold and $g(t)$ be a smooth one-parameter family of metrics which vary along the direction of a symmetric 2-tensor $v_{ij}(t)$,i.e.*

$$\frac{\partial g_{ij}}{\partial t} = v_{ij}$$

*Then the $\alpha$ −th Laplacian Eigen value $\lambda_\alpha$ evolves according to the following equation*

$$\frac{d\lambda_\alpha}{dt} = -\int_{M^n} v_{ij}\nabla^i f_\alpha \nabla^j f_\alpha\, d\mu - \frac{\lambda_\alpha}{2}\int_{M^n} v f_\alpha^2\, d\mu + \frac{1}{2}\int_{M^n} v|\nabla f_\alpha|^2\, d\mu$$

*Where $v = g^{ij}v_{ij}$ and $f(x,t)$ is a time- dependent eigen function of $\lambda_\alpha(g(t))$ with $\int_{M^n} f_\alpha(x,t)^2\, d\mu_{g(t)} = 1$ for all t.*

*Proof since $\int_{M^n} f^2 d\mu = 1$, we have*

$$\lambda = \int_{M^n} |\nabla f|^2\, d\mu$$

*at first we compute $\frac{\partial}{\partial t}|\nabla f|^2$.*

$$\frac{\partial}{\partial t}|\nabla f|^2 = \frac{\partial}{\partial t}\left(g^{ij}\nabla_i f \nabla_j f\right) = \left(\frac{\partial}{\partial t}g^{ij}\right)\nabla_i f \nabla_j f + 2g^{ij}\nabla_i\left(\frac{\partial}{\partial t}f\right)\nabla_j f$$

$$= -g^{ik}g^{jl}\left(\frac{\partial}{\partial t}g_{kl}\right)\nabla_i f \nabla_j f + 2g^{ij}\nabla_i\left(\frac{\partial f}{\partial t}\right)\nabla_j f$$

$$= -v_{kl}\nabla^k f \nabla^l f + 2g^{ij}\nabla_i\left(\frac{\partial f}{\partial t}\right)\nabla_j f$$

*We also have*

$$\frac{\partial}{\partial t}d\mu = \frac{v}{2}d\mu$$

*Moreover, since $\int_{M^n} f^2 d\mu = 1$, take $\frac{d}{dt}$ on both side gives*

$$0 = \int_{M^n} 2f\frac{\partial f}{\partial t} + \frac{v}{2}f^2 d\mu$$

*So*

$$2\int_{M^n} f\frac{\partial f}{\partial t}d\mu = \frac{-1}{2}\int_{M^n} v f^2 d\mu \qquad (2.5)$$

*Now we compute the $\frac{d\lambda}{dt}$.*

$$\frac{d\lambda}{dt} = \int_{M^n} \frac{\partial}{\partial t}(|\nabla f|^2 d\mu)$$

$$= \int_{M^n} \frac{\partial}{\partial t}|\nabla f|^2 + \frac{v}{2}|\nabla f|^2 d\mu$$

$$= -\int_{M^n} v_{ij}\nabla^i f \nabla^j f d\mu + 2\int_{M^n} g^{ij}\nabla_i\left(\frac{\partial}{\partial t}f\right)\nabla_j f d\mu + \frac{1}{2}\int_{M^n} v|\nabla f|^2 d\mu$$

$$= -\int_{M^n} v_{ij}\nabla^i f \nabla^j f d\mu - 2\int_{M^n}\left(\frac{\partial}{\partial t}f\right)\Delta f d\mu + \frac{1}{2}\int_{M^n} v|\nabla f|^2 d\mu$$

$$= -\int_{M^n} v_{ij} \nabla^i f \nabla^j f d\mu + 2\lambda \int_{M^n} f \frac{\partial}{\partial t} f d\mu + \frac{1}{2} \int_{M^n} v |\nabla f|^2 d\mu$$

*by plugging (2.5) in last equation, proof will be complete.*
*So we obtain the following lemma directly*

*Lemma2.5. If $\frac{\partial g_{ij}}{\partial t} = -2R_{ij}$, Then the evolution of $\lambda_\alpha(g(t))$ is given by*

$$\frac{d\lambda_\alpha}{dt} = 2\int_{M^n} R_{ij} \nabla^i f_\alpha \nabla^j f_\alpha d\mu + \lambda_\alpha \int_{M^n} R f_\alpha^2 d\mu - \int_{M^n} R |\nabla f_\alpha|^2 d\mu$$

*Where $f_\alpha(x,t)$ is an Eigen function of $\lambda_\alpha(g(t))$ with $\int_{M^n} f_\alpha(x,t)^2 d\mu = 1$ for all $t \in [0,T)$.*

## Cigar Soliton

*In Hamilton's program for the Ricci flow on 3-manifold, via dimension reduction, the cigar soliton is a potential singularity model.*

*Definition2.3. Hamilton's cigar soliton is the complete Riemannian surface $(R^2, g_\Sigma)$ where*

$$g_\Sigma = \frac{dx^2 + dy^2}{1 + x^2 + y^2}$$

*Where $x^2 := dx \otimes dx$. As a solution to the Ricci flow, its time-dependent version is:*

$$g_\Sigma(t) = \frac{dx^2 + dy^2}{e^{4t} + x^2 + y^2}$$

*From the change of variables $\tilde{x} = e^{-2t}x$ and $\tilde{y} = e^{-2t}y$, we see that $g_\Sigma = \frac{d\tilde{x}^2 + d\tilde{y}^2}{1 + \tilde{x}^2 + \tilde{y}^2}$ is isometric to $g_\Sigma = g_\Sigma(0)$. That is, if we define the 1- parameter group diffeomorphisms $\varphi_t: R^2 \to R^2$ by $\varphi_t(x,y) = (e^{-2t}x, e^{-2t}y)$, Then $g_\Sigma(t) = \varphi_t^* g_\Sigma(0)$. So $g_\Sigma(t)$ is a Ricci soliton. Thus, for each two times $t_1, t_2 \in (-\infty, \infty), g(t_1)$ is isometric to $g(t_2)$.*

*In polar coordinates, we may rewrite to cigar metric as*

$$g_\Sigma = \frac{dr^2 + r^2 d\theta^2}{1 + r^2}$$

*If we define the new radial distance variable*

$$s := \operatorname{arcsinh} r = \log\left(r + \sqrt{1 + r^2}\right)$$

*Then we may rewrite $g_\Sigma$ as*

$$g_\Sigma = ds^2 + \tanh^2 s\, d\theta^2$$

*Also we can write*

$$g_\Sigma = \left(1 - \frac{M}{\rho}\right) d\theta^2 + \left(1 - \frac{M}{\rho}\right)^{-1} \frac{d\rho^2}{4\rho^2}$$

$$= (e^{-2z} + 1)^{-1}(dz^2 + d\theta^2)$$

*Where $\rho = M\cosh^2 s$ and $z = \frac{1}{2}\log\left(\frac{r}{M} - 1\right)$ (where $r = \sqrt{x^2 + y^2}$)*

*And $\theta \in S^1(1) = R/2\pi Z$ where $S^1(1)$ denotes the circle of radius 1.*

*Corollary 2.2 (uniqueness of the cigar) if $(M^2, g(t))$ is a complete gradient Ricci soliton with positive curvature, then $(M^2, g(t))$ is the cigar soliton.*

**Ricci soliton**

*Definition 2.4. Let $(M, g)$ be a Riemannian manifold. A Ricci soliton structure on $M$ is a smooth vector field $X$ satisfying the Ricci equation*

$$Ric + \frac{1}{2} L_X g = \lambda g$$

*For some constant $\in R$.*

*Here*

*Ric = Ricci curvature of $M$*

*$L_X$ = Lie derivative in the direction $X$ (a sort of directional derivative in the direction of $X$ for tensor fields. Built using the local 1-parameter group of diffeomorphisms generated by )*

*We say the Ricci soliton is*

*Contractive or shrinking if $\lambda > 0$*

*Steady if $\lambda = 0$*

*Expansive if $\lambda < 0$*

*If $X = \nabla f$ for some smooth function $f: M \to R$. we say that $(M, g, \nabla f)$ is a gradient Ricci soliton with potential f. In this situation, the soliton equation reads*

$$Ric + Hess(f) = \lambda g$$

*Ricci solitons are a generalization of Einstein manifolds .they give rise to self similar solutions of the Ricci flow and arise as the blow up of some of the singularities of the Hamilton Ricci flow.*

# OVERVIEW ON POINCARE CONJECTURE

## Chapter 3

# Chapter 3

# The Ricci flow approach to Poincare conjecture

*We will spend this section giving a high-level overview of Perelman's Ricci flow –based proof of the Poincare conjecture ,and in particular how that conjecture is reduced to verifying a number of (highly non-trivial) facts about Ricci flow .Our exposition is based on [13],[14],[15],[16],[17] and [18]*

*We start with a question.*

*"If M is a closed* $3-$*manifold with trivial fundamental group, then is M diffeomorphic to* $S^3$.

*The Poincare conjecture is that the answer to this question is "yes". In 1980's Thurston developed another approach to 3-manifolds with Riemannian metrics of constant negative curvature -1 .These manifolds, which are locally isometric to Hyperbolic 3-space ,are called Hyperbolic manifold .There are fairly obvious obstructions showing that not every 3-manifold can admit such a metric . Thurston formulated a dereral conjecture that roughly says that the obvious obstructions are the only ones, Should they vanish for a particular 3-manifold, then that manifold admits such a metric .The aspect of Thurston's Geometrization conjecture that is most relevant for us is that the conjectural existence of especially nice metrics on 3-manifolds suggests a more analytic*

*approach to the problem of classifying 3-manifolds.Hamilton's formalized one such approach by introducing the Ricci flow on the space of Riemannian metrics .*

*Hamilton made significant progress on the program he initiated by establishing many crucial analytic estimates for understanding the evolving metric and its curvatures. After Perelman, on this work has claimed to surmount all of the various technical, geometric and analytic difficulties of Hamilton's program .In this way he claims to have established Thurston's Geometrization conjecture and hence the Poincare Conjecture .*

*Definition3.1. A homogeneous Riemannian manifold $(M,g)$ is one whose group of isometries acts transitively on the manifold.*

*Examples of homogeneous manifolds are the round sphere $S^n$ ,Euclidean space $R^n$ and Hyperbolic space $H^n$ .*

*We say that a Riemannian manifold is modeled on a given homogeneous manifold $(M,g)$ if every point of the manifold has a neighborhood isometric to an open set in $(M,g)$ .Such manifolds are called locally homogeneous manifolds provided that they are complete .*

*Remark3.1. In dimension2 there are four simply connected homogeneous manifolds up to isometry : $S^2$ , $R^2$ , $H^2$ ,and G ,with a left invariant metric ,where G is the group $R \ltimes R^*$ with the natural action of $R^*$ on R.*

*Remark 3.2 .the geometry and topology of the surface are related by the Gauss –Bonnet formula $2\pi\chi(M) = \int_M k dV$ ,where $\chi(M)$ is the Euler characteristic of M ,Which is related with the genus g of M by $\chi(M) = 2 - 2g$ ,and every oriented closed surface has a genus g and can be described as a sphere with g handles glued to it ,where a handle is $I \times S^1$ .*

*Theorem3.1 (uniformization in dimension 2) let X be a compact surface then X admits a locally homogeneous metric locally isometric to one of constant curvature models above .The model will be positively curved if $\chi(X) > 0$ ,flat if $\chi(X) = 0$ ,and negatively curved or hyperbolic if $\chi(X) < 0$ .*

*The Sphere (or prime) Decomposition:*

*The first two main steps in the way to classification of 3-manifolds are the sphere Decomposition Theorem.*

*Definition 3.2.(connected sum) if a closed 3-manifold M contains an embedded sphere $S^2$ separating M into two components, we can split M along this $S^2$ into manifolds $M_1$ and $M_2$ with boundary $S^2$. We can then fill in these boundary spheres with 3-balls to produce two closed manifolds $N_1$ and $N_2$. One says that M is the connected sum of $N_1$ and $N_2$, and one writes $M = N_1 \# N_2$. This splitting operation is commutative and associative.*

*There is also a strong relationship between the topology of a connected sum and that of its components.*

*Remark.3.3 Let M and M' be connected manifolds of the same dimension*

*1) M#M' is compact if and only if M and M' are both compact*

*2) M#M' is Orientable if and only if M and M' are both Orientale*

*3) M#M' is simply connected if and only if M and M' are both simply connected*

*Remark3.4 .on rather trivial possibility for the decomposition of M a connected sum is $= M \# S^3$.*

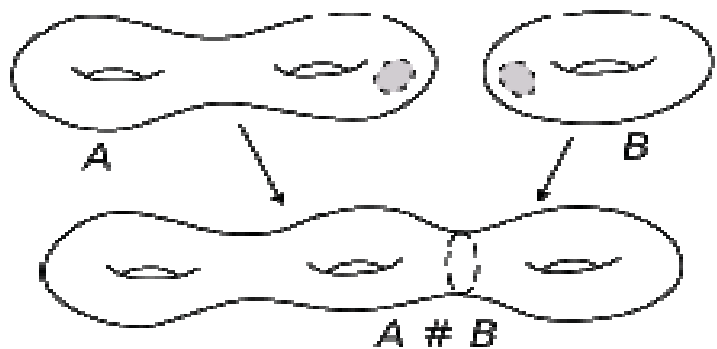

Illustration of connected sum.

*Now we give a theorem from Grushko about connected sum.*

*Theorem3.2 (Grushko) Let M be a compact, connected 3-manifold, and $\pi_1(M) = G_1 * G_2$ Then $M = M_1 \# M_2$ where $\pi_i(M_i) = G_i$ $(i = 1,2)$.*

*Definition3.3 (prime manifolds) a closed 3-manifold is called prime if every separating embedded 2-sphere bounds a 3-ball.*

*Or a another simple definition: a 3-manifold X is said to be prime if it is not diffeomorphic to $S^3$ and if every $S^2 \subset X$ that separates X into two pieces has the Property that one of the two pieces is diffeomorphic to a 3-ball.*

*Definition 3.4 : an embedded 2-sphere $S^2 \hookrightarrow M^3$ is essential if it does not bound a ball in M or is not parallel to a sphere in $\partial M^3$.*

*An orientable 3-manifold M is irreducible if any embedding of the 2-sphere into M extends to an embedding of the 3-ball into M.*

*One of the first theorems in the topology of three manifolds in due to Keneser in 1929.*

*Theorem3.3 .every closed, oriented 3-manifold admit a decomposition as a connected sum of oriented prime 3-manifolds, called its prime factors. This decomposition is unique up to the order of the prime factors (and orientation-preserving diffeomorphism of the factors).*

*There are a countably infinite number of prime 3-manifolds up to diffeomorphism.*

*There are three main classes of prime manifolds*

*Type I: With finite fundamental group .All the known examples are the spherical 3-manifolds, of the form, $M = \frac{S^3}{\Gamma}$, where $\Gamma$ is a finite subgroup of $SO(4)$ acting freely on $S^3$ by rotations .Thus $\Gamma = \pi_1(M)$ .It is and old conjecture that spherical 3-manifolds are the only closed 3-manifolds with finite fundamental group (it is ,in fact ,Poincare conjecture ) .*

*Type II) with infinite cyclic fundamental group .there is only one prime 3-manifold satisfying this condition: $S^1 \times S^3$ .This also the only orientable 3-manifold that is prime but not irreducible .it is also the only prime orientable 3-manifold with non-trivial $\pi_2$ .*

*Type III) with infinite non cyclic fundamental group .these are $K(\pi,1)$ Manifolds (also called aspherical); i.e , manifolds with contractible universal cover. Any irreducible 3-manifold M, with $\pi_1$ infinite is a $K(\pi,1)$.*

*Suppose that we start with a 3-manifold M which is connected and not prime.*

*Then we can decompose $M = N_1 \# N_2$ ,where no $N_i$ is a sphere .Now ,either each $N_i$ is irreducible ,or we can iterate this procedure. The theorem of Keneser (1929) states that this procedure always stops after a finite number of steps, yielding a manifold M such that each connected component of M is irreducible .*

*Theorem3.4 (prime Decomposition) Let M be a orientable closed 3-manifold.Then M admits a finite connected sum decomposition*

$$M = (K_1 \# K_2 \dots \# K_P) \# (L_1 \# L_2 \dots \# L_q) \# (\#_1^r S^2 \times S^1)$$

*The K and L factors here are closed and irreducible 3-manifolds.The K factors have infinite fundamental group and are aspherical 3-manifold (are of type III), while the L factors have finite fundamental group and have universal cover a homotopy 3-sphere (are of type I)*

### The Thurston Geometrization conjecture

*In dimension three, every finite volume, locally homogeneous manifold is modeled on one of the eight homogeneous manifolds listed below .First; we have the constant (sectional) curvature examples:*

1) *$S^3$,of constant curvature +1*
2) *$R^3$,which is flat*
3) *$H^3$ of constant curvature -1*

*Next we have the homogeneous 3-manifold with product metrics:*

*4) $S^2 \times R$*

*5) $H^2 \times R$*

*Finite volume locally homogeneous manifolds modeled on $S^2 \times R$ are automatically compact and either are isometric $S^2 \times S^1$ or $RP^3 \# RP^3$.*

*Finite volume locally homogeneous manifolds modeled on $H^2 \times R$ either are of the form $\sum \times S^1$, where $\sum$ is a finite area hyperbolic surface, or are finitely covered by such manifolds.*

*Lastly, we have the homogeneous manifolds $(M, g)$ where $M$ is a simply connected lie group and $g$ is a left-invariant metric.*

*Three-dimensional examples of this type admitting locally homogeneous examples of finite volume are:*

*Definition 3.5. Heisenberg group. Is a group of $3 \times 3$ upper triangular matrices of the form $\begin{pmatrix} 1 & a & c \\ 0 & 1 & b \\ 0 & 0 & 1 \end{pmatrix}$. Elements $a, b, c$ can be taken from some arbitrary commutative ring.*

*6) **The unipotent group** (Heisenberg group) Locally homogeneous manifolds modeled on this group are called Nil-manifolds. Note that this is the only 3-dimensional nilpotent but not abelian connected and simply connected Lie group; This explains the term Nil geometry topologically, Nil is diffeomorphic to $R^3$ under the map*

$$Nil \ni \gamma = \begin{pmatrix} 1 & a & c \\ 0 & 1 & b \\ 0 & 0 & 1 \end{pmatrix} \to (a, b, c) \in R^3$$

*Under this identification, left multiplication by $\gamma$ corresponds to the map*

$$L_\gamma(x, y, z) = (x + a, y + b, z + ay + c)$$

*In other words, from this point of view, $R^3$ has the multiplication*

$$(x_0, y_0, z_0)(x, y, z) = (x + x_0, y + y_0, z + z_0 + x_0 y)$$

*It is easy to prove that Nil is a Lie group, so it admits a metric invariant under left multiplication; we shall take $ds^2 = dx^2 + dy^2 + dz^2$ at (0,0,0) -the unit of the Heisenberg group and the extend $ds^2$ at all other points of X as a left invariant metric. The result is*

$$ds^2 = dx^2 + dy^2 + (dz - xdy)^2$$

*Next we check that $ds^2$ is invariant under multiplications*

$$L_\gamma^*(ds^2) = d(L_\gamma^* x)^2 + d(L_\gamma^* y)^2 + \left(d(L_\gamma^* z) - L_\gamma^* x d(L_\gamma^* y)\right)^2$$

$$= d(xoL_\gamma)^2 + d(yoL_\gamma)^2 + \left(d(zoL_\gamma) - (xoL_\gamma) d(yoL_\gamma)\right)^2$$

$$= d(x + a)^2 + d(y + a)^2 + \left(d(z + ay + c) - (x + a)d(y + b)\right)^2$$
$$= dx^2 + dy^2 + (dz + ady - xdy - ady)^2 = ds^2$$

*If we identify $S^1$ with the interval $[0,2\pi]$ with the ends identified then a point $\theta \in S^1$ acts on Nil by $\mathcal{A}: S^1 \times Nil \to Nil$ such that*

$$(\theta, (x, y, z)) \xrightarrow{\mathcal{A}} (\theta x, \theta y, \theta z) \text{ Where}$$

$$f(x) = \begin{cases} \theta x = x\cos\theta - y\sin\theta \\ \theta y = x\sin\theta + y\cos\theta \\ \theta z = z + \frac{1}{2}\sin\theta(\cos\theta(x^2 - y^2) - 2\sin\theta xy) \end{cases}$$

*Where $\mathcal{A}$ is an action of $S^1$ on Nil which is a group of automorphisms of Nil preserving the above metric.*

*Finite volume manifolds with this geometry are compact and orientable and have the structure of Seifert fiber space.*

*Also we can identify Nil with the subset $G$ of $\mathbb{C}^2$ defined as $\{(u,v) \in \mathbb{C}^2 : Imv = |u|^2\}$ and with multiplication on $G$ defined by*

$$(u,v).(u',v') = (u+u', v+v'+2iu\bar{u}')$$

*An isomorphism from Nil to $G$ can be given by the formula*

$$(x,y,z) = \left(\frac{1}{2}(x+iy), z - \frac{1}{2}xy + \frac{1}{4}i(x^2+y^2)\right)$$

### 7) Sol geometry

*At first we give two points .*

*1) A Lie group $G$ is said to be solvable if it is connected and its Lie algebra is solvable.*

*2) Let $G$ and $H$ two Lie groups and consider a homomorphism from $G$ to the abstract group of automorphisms of $H$, that is, $\rho: G \to Aut(H)$. The semi direct product $H \times_\rho G$ of $H$ and $G$ with respect to $\rho$ is the product manifold $H \times G$ endowed with the Lie group structure given by*

$(h,g).(h',g') = (h\rho(g)h', gg')$   ,  $(h,g)^{-1} = (\rho(g^{-1})h^{-1}, g^{-1})$   $\forall h, h' \in M$ *and* $\forall g, g' \in G$

*Topologically, we can identify Sol with $R^3$ so that the multiplication is given by*

$$(x_0, y_0, z_0)(x,y,z) = (x + e^{-z_0}x_0, y + e^z y_0, z + z_0)$$

*Clearly $(0,0,0)$ is identity and $xy$−plane is a normal subgroup isomorphic to $R^2$. In fact,*

$$(x,y,z)^{-1}.(a,b,0).(x,y,z) = (-xe^{-z}, -ye^z, -z)(a,b,0).(x,y,z)$$

$$= (a - xe^{-z}, b - ye^z, -z)(x,y,z) = (*,*, 0)$$

*Metrically, Sol is just $R^3$ but endowed with the left invariant Riemannian metric which at $(x, y, z)$ is*

$$ds^2 = e^{2z}dx^2 + e^{-2z}dy^2 + dz^2 \quad (3.1)$$

*Another equivalent approach to the algebraic definition of Sol will be given along the proof of Thurston theorem. In this way; we can say that Sol is the Unimodular Lie group completely determined by*

$$[e_1, e_2] = 0, \quad [e_2, e_3] = e_1, \quad [e_3, e_1] = -e_2$$

*Where $\{e_1, e_2, e_3\}$ is an orthonormal basis of eigenvectors of the Lie algebra of $X \equiv Sol$. As the generators $e_1$ and $e_2$ are commute. Sol contains a copy of $R^2$ which is a normal subgroup and the quotient group is R. The group is therefore a semi direct product of $R^2$ with R. Therefore the transformations in this basis are of the form $t \to \begin{pmatrix} e^t & 0 \\ 0 & e^{-t} \end{pmatrix}$.*

*The metric (1.3) is preserved by the group G of transformations of X of the form*

$$(x, y, z) \to (\varepsilon e^{-c}x + a, \varepsilon' e^c y + b, z + c) \text{ or } (\varepsilon e^{-c}y + a, \varepsilon' e^c x + b, z + c)$$

*Where $a, b, c \in \mathbb{R}$ and $\varepsilon, \varepsilon' = \pm 1$.*

**8) $G = \widetilde{PSL_2}(\mathbb{R})$** *The universal covering group of $PSL_2(\mathbb{R})$. This manifold can also be viewed by as the universal covering of the unit tangent bundle to $H^2$ with its induced metric. Finite volume locally homogeneous manifolds modeled on this example are circle bundles over hyperbolic surfaces.*

*Theorem 3.5 suppose that $X^3$ is connected and orientable and admits a locally homogeneous Riemannian metric of finite volume. Then, X is diffeomorphic to the interior of a compact 3-manifold with boundary, all of whose boundary components are Tori. Furthermore, each of these tori has fundamental group which injects into the fundamental group of X. if X is non-compact, then it is modeled on $H^3, H^2 \times \mathbb{R}$ or $\widetilde{PSL_2}(\mathbb{R})$, and hence, X either is a hyperbolic 3-manifold or is Seifert –fibered with hyperbolic two –dimensional orbifold base.*

*The Geometrization conjecture reads as follows:*

**Geometrization conjecture:** *Every closed orientable $3$-manifold $M$ is a connected sum of closed $3$-manifold $M_i$ such that for every $i$ there is a finite collection of pair wise disjoint embedded Tori $T_{ij} \subset M_i$ such that*

1) *The Tori $\{T_{ij}\}_j$ are incompressible in $M_i$ (i.e. $\pi_1 T_{ij} \hookrightarrow \pi_1 M_i$)*

2) *The components of $M_i \setminus \cup_j T_{ij}$ are diffeomorphic to metric quotients with finite volume of one of the following 8 homogeneous geometries*

$S^3, S^2 \times S^1, H^3, R^3, H^2 \times R, \widetilde{PSL}(2,R), Nil$ *or* $Sol$

*Note that the Poincare conjecture follows easily from the geometrization conjecture*

**Actions on geometric manifolds**

*Definition 3.6 $((X,\Gamma) - structure)$. Let $X$ be a topological manifold and $\Gamma$ a group acting on $X$. An $(X,\Gamma)$-structure for a manifold $M$ is a maximal $(X,\Gamma)$-compatible collection of charts $\varphi_i: U_i \to X$ covering $M$. Two charts $\varphi_i, \varphi_j$ are $(X,\Gamma)$-compatible, if on each component $V$ of $U_i \cap U_j$ the coordinate change $\varphi_j \circ \varphi_i^{-1}|_{\varphi_i^{-1}(V)}$ is the restriction of some element $g \in \Gamma$.*

*Definition 3.7 (model geometry, geometric manifold) a model geometry is a smooth, simply connected manifold $X$ with a Lie-group $\Gamma$ acting transitively on $X$ such that,*

1) *$\Gamma$ has compact point stabilizer*
2) *$\Gamma$ Is maximal in the sense that it is not contained in any larger group of diffeomorphisms of $X$ with compact point stabilizer.*
3) *There exists at least one compact manifold with an $(X,\Gamma)$-structure if $(X,\Gamma)$ is a model geometry, then an $(X,\Gamma)$-structure is called geometric structure and a manifold with a geometric structure is called geometric manifold.*

*We consider smooth actions $\rho: G \curvearrowright M$ of a finite group $G$ on a smooth manifold $M$.*

*Definition 3.8 (standard action) let $(X, \Gamma)$ be a model geometry and $M$ an $(X, \Gamma)-$ manifold. We say, the action $\rho: G \curvearrowright M$ is standard, if there exists a $\rho(G)$-invariant complete locally homogeneous metric on $M$.*

*Theorem 3.6: Let $M = S^2 \times S^1$, and $\rho: G \curvearrowright M$ a smooth finite group action. Then the action $\rho$ is standard.*

*Any locally homogeneous manifold modeled on the $S^3$ of constant curvature +1 is a Riemannian manifold of constant positive sectional curvature.*

*These are of the from $S^3/\Gamma$ where $\Gamma$ is a finite subgroup of $SO(4)$ acting freely on $S^3$, $RP^3$, lens space, as well as the quotients by the symmetry groups of the exceptional regular solids. These manifolds are called spherical space-form.*

*Definition 3.9. A Riemannian four-manifold is said to have positive isometric curvature if for every orthonormal four-frame the curvature tensor satisfies*

$$R_{1313} + R_{1414} + R_{2323} + R_{2424} > 2R_{1234}$$

*Definition 3.10. An incompressible space from $N^3$ in a four-manifold $M^4$ is a three dimensional sub manifold diffeomorphic to $S^3/\Gamma$, (the quotient of the three sphere by a group of isometries without fixed point) such that the fundamental group $\pi_1(N^3)$ injects into $\pi_1(M^4)$.*

*The space form is said to be essential unless $\Gamma = \{1\}$, or $\Gamma = Z_2$ and the normal bundle is non-orientable.*

*Theorem 3.7 (Bing-Long chen Xi-Ping Zhu) let $M^4$ be a compact four-manifold with no essential incompressible space-form and with a metric $g_{ij}$ of positive isotropic curvature. Then we have a finite collection of smooth solutions ${g^{(k)}}_{ij}(t)$, $k = 0,1,2,\ldots,m$. to the Ricci flow, defined on $M_k^4 \times [t_k, t_{k+1})$, $(0 = t_k < \cdots < t_{m+1})$ with $M_0^4 = M^4$ and ${g^{(0)}}_{ij}(t_0) = g_{ij}$, which go singular as $t \to t_{k+1}$, such that the following properties hold:*

i) *For each $k = 0,1,\ldots,m-1$, the compact (possible disconnected) four-manifold $M_k^4$ contains an open set $\Omega_k$ such that the solution $g^{(k)}{}_{ij}(t)$ can be smoothly extended to $t = t_{k+1}$ over $\Omega_k$ ;*

ii) *For each $k = 0,1,\ldots,m-1, \left(\Omega_k, g^{(k)}{}_{ij}(t_{k+1})\right)$ and $\left(M^4{}_{k+1}, g^{(k+1)}{}_{ij}(t_{k+1})\right)$ contain compact (possible disconnected) four-dimensional sub manifolds with smooth boundary, which are isometric and then can be denoted by $N_k^4$,*

iii) *For each $k = 0,1,\ldots,m-1$, $M^4{}_k \setminus N^4{}_k$ consists of a finite number of disjoint pieces diffeomorphic to $S^3 \times I$, $B^4$ or $RP^4 \setminus B^4$, while $M^4{}_{k+1} \setminus N^4{}_k$ consists of a finite number of disjoint pieces diffeomorphic to $B^4$ ;*

iv) *For $k = m$, $M_m^4$ is diffeomorphic to the disjoint union of a finite number of $S^4$, or $RP^4$, or $S^3 \times S^1$, or $S^3 \widetilde{\times} S^1$, or $RP^4 \# RP^4$. (note $S^3 \widetilde{\times} S^1$, is equal with $S^3 \times S^1/Z_2$ where $Z_2$ flips $S^3$ antipodally and rotates $S^1$ by $180°$ )*

*As a direct consequence we have the following classification result of Hamilton.*

*Corollary (Hamilton): a compact four manifold with no essential incompressible space form and with a metric of positive isometric curvature is diffeomorphic to $S^4$, or $RP^4$, or $S^3 \times S^1$, or $S^3 \widetilde{\times} S^1$, or a connected sum of them.*

*Hanilton introduced a program to study all 3-manifolds using the Ricci flow.*

- **Short-time existence and uniqueness** :*if $g_0$ is a smooth metric on a compact manifold, then there is some $\varepsilon > 0$ depending on $g_0$ and a unique solution to the Ricci flow equation defined for $t \in [0, \varepsilon)$ with $g = g_0$.*

- **Curvature characterization of singularity formulation** :*if the solution exists on the time interval $[0, T)$ but does not extend to any strictly larger time interval, Then there is a point $x$ in the manifold for which curvature tensor $Rm(x,t)$ of the metric $g(t)$ is unbounded as $t$ approaches $T$.*

*Hamilton was discovered quite early that the Ricci flow may develop singularities even in the case of a sphere if the Ricci curvature is not positive. an example is the so-called*

*neck pinch singularity .Perelman's work describes what happens to the Ricci flow nears a singularity and also how to Perform the surgery .The new flow is called Ricci flow with surgery.*

*Theorem 3.8(Hamilton) Let $X^3$ be a compact connected 3-manifold with non-negative Ricci curvature .Then one of the following happens:*

1) *The Ricci curvature becomes strictly positive for all $t > 0$ sufficiently small. In this case, the Ricci flow develops a singularity in finite time.*

   *As the singularity develops ,the diameter of the manifolds goes to zero .Rescaling the evolving family of metrics so that their diameters are one leads to a family of metrics converging smoothly to a metric of constant positive curvature .In particular ,the manifold is diffeomorphic to a spherical space-form (note that if $T < \infty$ and the curvature becomes unbounded as $t$ tends to $T$ ,we say the maximal solution develops singularities as $t$ tends to $T$ and $T$ is the singular time )*

2) *There is a finite cover of the Riemannian manifold which ,with the induced metric ,is a metric product of a compact surface of positive curvature and $S^1$ .this remains true for all the Riemannian metrics in the Ricci flow .The Ricci flow develops a singularity in finite time ,and the manifold in question is diffeomorphic to $S^2 \times S^1$ ,or $RP^3 \# RP^3$ (the case of $RP^3 \# RP^3$ is interesting in that it is apparently the only non prime 3-manifold which admit a geometric structure)*

3) *The metric is flat and the evolution equation is constant. In this case, of course, the manifold is covered by $T^3$.*

*Theorem3.9 (Hamilton) if the Ricci flow with initial conditions $(M, g_0)$ be a connected, compact and Ricci flow exists for all $t \in [0, \infty)$ and also the normalized curvature . $Rm(x,t)$ is bounded as $t \to +\infty$ .Then M satisfies Thurston's Geometrization conjecture .*

**Perelman's claims**

*Regions of High curvature in the flow :in order to do surgery we need to understand regions of the flow where the scalar curvature is large .of course ,since it is possible to*

*rescale the metric and time in any flow ,we must normalize in some fashion to have an invariant notion of large curvature .Thus ,we arrange that our flow has normalized initial conditions in the sense that at $t = 0$ the absolute value of the Riemannian curvature at each point is at most one, and the volume of any metric ball of radius one is at least half that of the unit ball in $R^3$.*

*From now on we implicity assume that all flows under consideration have normalized initial conditions.*

*For all the following definitions we fix $0 < \varepsilon < \frac{1}{2}$ .Set k equal to the greatest integer less than or equal to $\varepsilon^{-1}$ .in particular, $k \geq 2$.*

*Definition 3.11 suppose that we have a fixed metric $g_0$ on a manifold M and an open sub manifold $X \subset M$ .We say that another metric g on X is within $\varepsilon$ of $g_0|_X$ in the $C^{\left[\frac{1}{\varepsilon}\right]}$-Topology if, setting $k = \left[\frac{1}{\varepsilon}\right]$ we have*

$$\sup_{x \in X}\left(|g(x) - g_0(x)|_{g_0}^2 + \sum_{l=1}^{k}\left|\nabla_{g_0}^l g(x)\right|_{g_0}^2\right) < \varepsilon^2$$

*Where the covariant derivative $\nabla_{g_0}^l$ is the Levi-Civita connection of $g_0$ and norms are the point wise $g_0$ −norms on*

$$Sym^2 T^*M \otimes \underbrace{T^*M \dots T^*M}_{l-times}$$

*Definition 3.12 Let $(N, g)$ be a Riemannian manifold and $x \in N$ a point .Then an $\varepsilon$-neck structure on $(N, g)$ centered at x consists of a diffeomorphism*

$$\varphi: S^2 \times (-\varepsilon^{-1}, \varepsilon^{-1}) \to N$$

*With $x \in \varphi(S^2 \times \{0\})$ ,such that the metric $R(x)\varphi^*g$ is within $\varepsilon$ in the $C^{\left[\frac{1}{\varepsilon}\right]} - Topology$ of the product of the usual Euclidean metric on the open interval with the metric of constant Gaussian curvature $\frac{1}{2}$ on $S^2$ .*

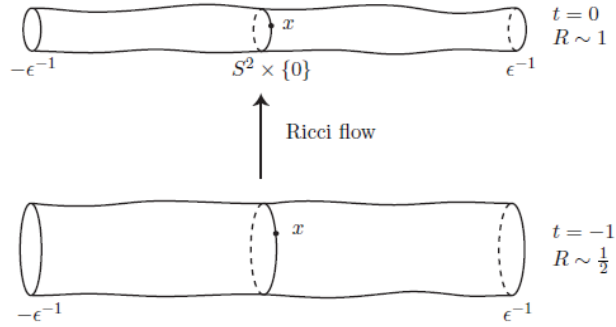

*Now we give a claim from Perelman.*

*Claim ( Perelman) there is $r > 0$ such that the following holds .Let $(M, g(t))$ be a Ricci flow with normalized initial conditions defined for $0 \leq t < T$ with M a closed ,orientable 3-manifold.Then , for any point $(x, t)$ in the flow with $R(x, t) \geq r^{-2}$ one of the following holds*

1) *The components of $M_t$ containing $(x, t)$ is diffeomorphic to a spherical space-form.*
2) *$(x, t)$ is the center of an $\varepsilon - neck$ in $M_t$*
3) *$(x, t)$ is contained in a sub manifold of $M_t$ diffeomorphic to $D^3$ with every point in the boundary being the center of an $\varepsilon -neck$ in $M_t$ .*

*Note that: such a region is called an $\varepsilon - cap$*

*Definition3.13 supposes that as $t \to T$ the Ricci flow becomes singular .there is an extension of the above theorem to this time as well. There is an open subset $\Omega_T \subset M$ consisting of all the points where the Riemannian Curvature tensor remains bounded .and according to papers of Shi and Hamilton, There is a limiting metric on this open subset as $t \to T$ .we denote it $(\Omega_T, g(T))$.*

*Definition3.14 in doing surgery there are three parameters that are fixed .The first is the coarse control parameter $\varepsilon > 0$ ,sufficiently small, a universal constant fixed once and for all .The other two parameters are non-increasing functions of t limiting to Zero as $t \to \infty$ .They are $\delta(t) > 0$,The fine control parameter, and $h(t)$,The scale parameter .we also fix an auxiliary parameter $\rho(t) = \varepsilon\delta(t)$ .*

*Definition 3.15 we fix $r > 0$ and a surgery parameter $\delta > 0$ .set $\rho = r.\delta$. Let $\Omega_T(\rho)$ be the subset of points $(x, T)$ for which $R(x, T) \leq \rho^{-2}$*

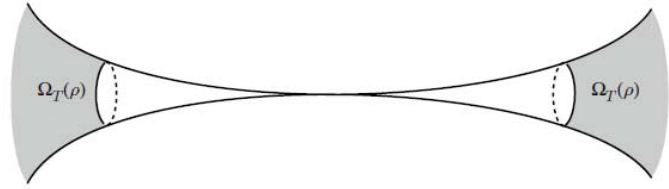

1) *An $\varepsilon-$Tube in $\Omega_T$ is a sub manifold diffeomorphic to $S^2 \times I$ such that each point is the center of an $\varepsilon-$neck in $\Omega_T$.*

2) *An $\varepsilon-$Circuit in $\Omega_T$ is a component of $\Omega_T$ which is a closed manifold and each one of its points is the center of an $\varepsilon-$neck. It is diffeomorphic to $S^2 \times S^1$.*

3) *an $\varepsilon-$horn is a closed subset $H \subset \Omega_T$ diffeomorphic to $S^2 \times [0,1)$ with boundary contained in $\Omega_T(\rho)$ such that every point of $H$ is the center of an $\varepsilon-$neck in $\Omega_T$. of course, the scalar curvature goes to infinity at the other end of H.*

4) *A double $\varepsilon-$horn is a component of $\Omega_T$ diffeomorphic to $S^2 \times (0,1)$ such that every point is the center of an $\varepsilon-$neck in $\Omega_T$. The scalar curvature goes to infinity at each end of this component.*

5) *A capped $\varepsilon-$horn is a component of $\Omega_T$ diffeomorphic to $\text{int}(D^3)$ such that each is either the center of an $\varepsilon-$neck or is contained in an $\varepsilon-$cap. The scalar curvature goes to infinity near the end of a capped $\varepsilon-$horn*

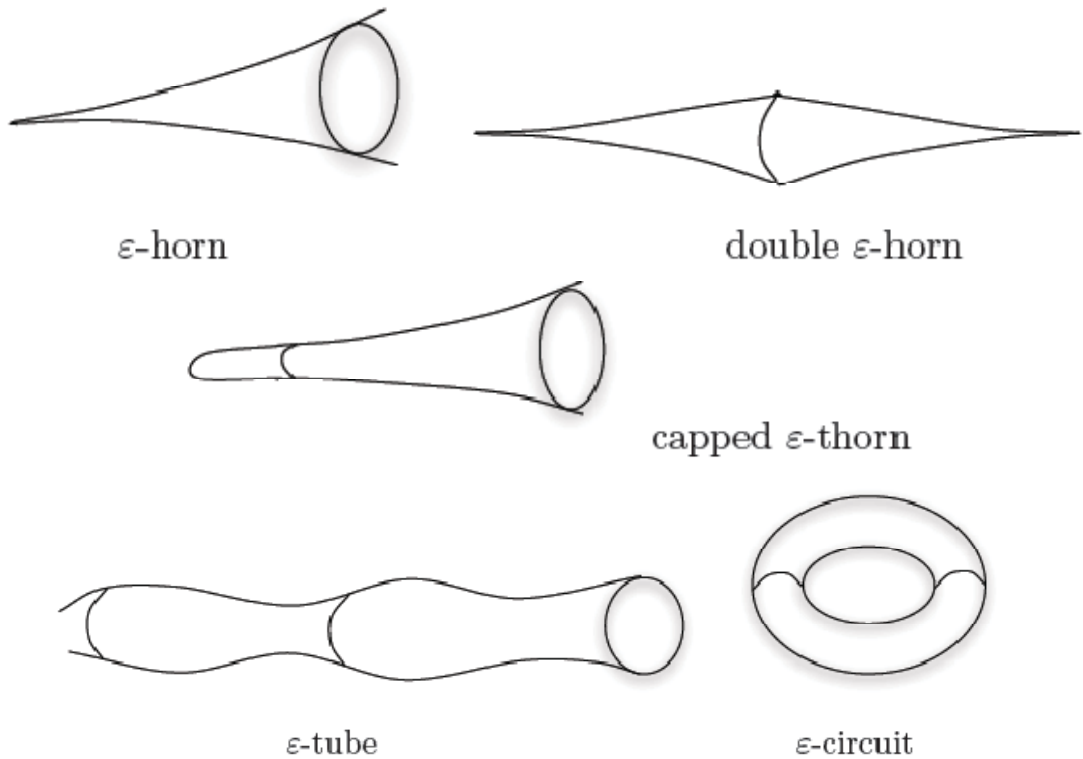

*Corollary At time T every $x \in \Omega_T - \Omega_T(\rho)$, is contained in one of the following :*

1) *A component of $\Omega$ containing x is diffeomorphic to a quotient of a sphere ,*
2) *An $\varepsilon -$circuit diffeomorphic to $S^2 \times S^1$,*
3) *An $\varepsilon -$Tube with boundary components in $\Omega_T(\rho)$*
4) *an $\varepsilon -$Cap with boundary in $\Omega_T(\rho)$*
5) *an $\varepsilon -$horn with boundary in $\Omega_T(\rho)$*
6) *a capped $\varepsilon -$horn*
7) *a double $\varepsilon -$horn*

**Perelman's functional $\mathcal{F}$ and $\mathcal{W}$**

*One of perelmans claims is that, he proved that Ricci flow is a gradient flow up to diffeomorphism .Of course note that Ricci flow is not gradient flow .for explain the trick of Perelman for it we at first define a function that called $\mathcal{F}$-functional*

$$\mathcal{F}(g,f) = \int_M (R + |\nabla f|^2) e^{-f} \, dVol.$$

*Where $g$ is a Riemannian metric on $M$ and $f$ is a smooth function.*

*Restrict this function to the subspace of $\mathcal{M} \times C^\infty(M)$ given by $(g,f)$ such that the volume form $e^{-f} dVol$ is constant, equal to a fixed one denoted by $dm$. The gradient flow equation for $\mathcal{F}$ is then*

$$\frac{\partial g}{\partial t} = -2(Ric + D^2 f), \qquad \frac{\partial f}{\partial t} = -\Delta f - R \qquad (3.1)$$

*Now the idea is to deform a solution $(g(t), f(t))$ of (3.1) to kill the term $D^2 f$ and obtain a solution $\tilde{g}(t)$ of Ricci flow. Now we recall a classical result of Riemannian geometry. if $\varphi_s$ is a family of diffeomorphisms such that $\frac{d}{ds}\varphi_s = \nabla_g f$, then $\frac{d}{ds}(\varphi_s^* g) = 2D^2 f$. Given a solution of (3.1) one can define diffeomorphisms $\varphi_t$ such that $\frac{d}{dt}\varphi = \nabla f$ at time $t$. Then $\tilde{g}(t) = \varphi_t^* g(t)$ and $\tilde{f}(t) = f(t) \circ \varphi_t$ solve the system :*

$$\frac{\partial \tilde{g}}{\partial t} = -2Ric, \quad \frac{\partial \tilde{f}}{\partial t} = -\Delta \tilde{f} - R + |\nabla \tilde{f}|^2$$

*The first equation gives a solution to the Ricci flow. Moreover, $\mathcal{F}(g(t), f(t)) = \mathcal{F}(\tilde{g}(t), \tilde{f}(t))$.*

**Elliptization and Hyperbolization conjecture**

*Definition 3.16 A Seifert fibration is a partition by circles (locally a product) except for a finite number of singular fibres, that have the following local model. We consider a cylinder, the product $D^2 \times [0,1]$ of a disc with an interval, and glue $D^2 \times \{0\}$ with $D^2 \times \{1\}$ by a rotation of finite order. The fibration by horizontal intervals induces a partition by circles of the solid torus, that is a fibration except for the singular fibre $\{0\} \times S^1$, which is shorter.*

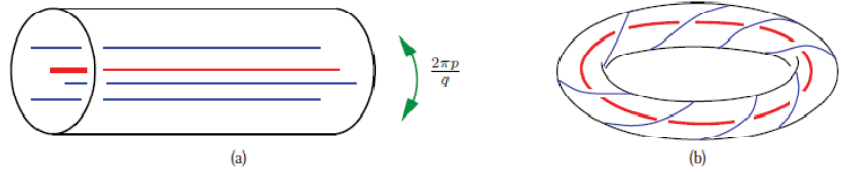

*We know a manifold said to be hyperbolic if its interior has a complete metric of constant curvature -1 .it is important to notice that manifolds with finite fundamental group cannot be hyperbolic.*

**Hyperbolization conjecture**: *Let $M^3$ be a prime, closed three manifold, with $\pi_1(M^3)$ infinite and such that every subgroup $Z \oplus Z < \pi_1(M^3)$ comes from $\partial(M^3)$ .Then $M^3$ is hyperbolic.*

*Also a 3-manifold is said to be elliptic if it admits a metric of constant curvature +1.*

*If $M^3$ is a closed three manifold and $\pi_1(M^3)$ is finite, then equipping $M^3$ with an elliptic metric is equivalent to admit a Seifert fibration, using the fact both families elliptic manifolds and Seifert fibered ones are classified.*

**Elliptization conjecture**: *let $M^3$ be a prime, closed three-manifold, with $\pi_1(M^3)$ finite .Then $M^3$ is elliptic.*

*Now we explain some partial result about Geometrization before Perelman*

*Definition 3.17 A 3-manifold is called a graph manifold if it is a union of Seifert manifolds along the boundary, consisting of Tori.*

*Definition3.18 on orientable 3-manifold is called sufficiently large if it contains an property embedded incompressible surface.*

*The next result is Thurston's Geometrization for sufficiently large manifolds, that was one of the main evidences to support this conjecture*

*A surface $F^2 \subset M^3$ is called properly embedded if it is embedded and $\partial F = F \cap \partial M^3$*

*Theorem3.10 (Thurston) a sufficiently large three manifold satisfies the Geometrization conjecture.*

*The following theorem is involves convergence groups*

*Theorem3.11 let $M^3$ be a compact irreducible 3- manifold .if $\pi_1(M^3)$ has an infinite cyclic normal sub-group, then $M^3$ is Seifert fibered*

**Consequences of Geometrization :**

*Some consequences of Geometrization in 3-dimensional topology are listed in this section:*

*Theorem3.12 (Borel conjecture in dimension three) if two aspherical compact three manifolds are homotopically equivalent, then they are homeomorphic .*

*Theorem3.13: compact aspherical three manifolds are classified by its fundamental group.*

*Theorem3.14 (Gromov-Lawson) if a compact three manifold admits a metric of non-negative scalar curvature then it is either flat or a connected sum of elliptic manifolds and $S^2 \times S^1$ or its quotients.*

**Hamilton claims on Geometrization conjecture**

*In the three dimensional case, the curvature operator (acting on 2-forms) diagonalizes*

$$\begin{pmatrix} \alpha_1 & 0 & 0 \\ 0 & \alpha_2 & 0 \\ 0 & 0 & \alpha_3 \end{pmatrix}$$

*So that the $\alpha_i$ are functions on M .Then the $\frac{\alpha_i}{2}$ are sectional curvatures,*

$$\begin{pmatrix} \dfrac{\alpha_2 + \alpha_3}{2} & 0 & 0 \\ 0 & \dfrac{\alpha_1 + \alpha_2}{2} & 0 \\ 0 & 0 & \dfrac{\alpha_1 + \alpha_3}{2} \end{pmatrix}$$

And $R = \alpha_1 + \alpha_2 + \alpha_3$

The evolution equations for the $\alpha_i$ are

$$\alpha_1' = \Delta\alpha_1 + \alpha_1^{\,2} + \alpha_2\alpha_3$$

$$\alpha_2' = \Delta\alpha_2 + \alpha_2^{\,2} + \alpha_3\alpha_1$$

$$\alpha_3' = \Delta\alpha_3 + \alpha_3^{\,2} + \alpha_1\alpha_2$$

*Hamilton proved, in dimension three, $Ric \geq 0$ , $Ric > 0$ , $Sec \geq 0$ and , $Sec > 0$ are conditions invariant under the Ricci flow .*

*Maximum principles where also used by Hamilton to control $\alpha_i$ in the following Theorem.*

*Theorem 3.15(Hamilton) if a compact three manifold $M^3$ admits a metric with $Ric > 0$ ,then the Ricci flow, after rescaling ,converges to a metric with constant positive sectional curvature .In Particular $M^3$ is elliptic .*

*Also Hamilton developed a strong maximum principle for tensors, and used it to show that if $Ric \geq 0$ ,Then one of the three possibilities happen*

1) the metric is flat
2) $Ric > 0$ at $t > 0$ , hence M is elliptic
3) The metric is locally a product $g = g_1 \oplus dx^3$ .In this case the manifold is diffeomorphic to $S^2 \times S^1$ or it is a quotient $RP^3 \# RP^3$ .

*Hence we conclude the following Theorem from Hamilton*

*Theorem3.16 a closed manifold with a Riemannian metric with $Ric \geq 0$ satisfies Thurston's Geometrization conjecture.*

*Hamilton –Ivey estimates are another example of cleaver application of maximum principles for tensors in dimension three .Let $\varphi: [-1, +\infty) \to [1, \infty)$ be the inverse map of $x \to \log x - x$ .Then we have the following Theorem*

*Theorem 3.17 (Hamilton-Ivey pinching) The inequalities*

*$R \geq -1$ and $\alpha_1, \alpha_2, \alpha_3 \geq -\varphi(R)$*

*Are invariant under the Ricci flow.*

**Hamilton's work on formation of singularities and the Zoom technique**

*The Ricci Tensor is invariant by homoteties ,Thus the rescaled metric still satisfies the Ricci flow equation .Now according to geometric tools we consider the Hamilton's work on singularities .One of geometric tools is zoom ,or parabolic rescaling .The idea is to dilate the metric and the time .*

*Definition 3.19 (parabolic rescaling) Let $g(t)$ be a Ricci flow on $M \times [0, T)$, $x_0 \in M, t_0 \in [0, T)$ Such that $R(x, t) \leq \lambda_0 \equiv R(x_0, t_0)$ for all $x \in M$ and $t \leq t_0$. Then*

$$g_0(t) = \lambda_0 g\left(t_0 + \frac{t}{\lambda_0}\right)$$

*Is a Ricci flow on $[-t_0 \lambda_0, (T - t_0)\lambda_0)$ and $R_0(x, t) \leq 1$ for every $x \in M$ and $t \leq 0$ .This is called a parabolic rescaling of $g(t)$ at $(x_0, t_0)$ .Perelman showed if The Ricci flow on a closed three manifold M encounters a singularity then an entire connected component disappear or there are nearby 2-Spheres on which to do surgery .To attack this ,Hamilton initiated a blowup analysis for Ricci flow .It is known that singularities arise from curvature blowups .That is ,if a Ricci flow solution exists on a maximal time interval $[0, T)$ ,with $T < \infty$ ,then $\lim_{t \to T^{-1}} \operatorname{Sup}_{x \in M} |Riem(x, t)| = \infty$ , where Riem denotes the sectional curvature .*

*Assume we have a singularity developing at time T under the Ricci flow $(M, g)$ .we take a sequence $x_i \in M$ and $t_i \to T$ so that $R(x_i) = \max_M R$ at time $t_i$ ,and we Parabolic rescale to have $R(x_i) = 1$ .we also move the time by a translation ,So that the initial $t_i$ becomes zero .In order to analyze the singularity ,The idea is to look at the limit of* *pointed Ricci flows (with base point $x_i$ at time 0, after the translation ).There is a compactness theorem for pointed Ricci flows ,proved we have a positive lower bound on the injectivity radius of the base point .if we had This lower bound ,Then there would be a convergent subsequence to a flow ,and Hamilton work would yield that the limit has the following properties :*

1) *It is an ancient solution ,i.e. defined on time $(-\infty, 0]$*
2) *The metric is complete*
3) *The sectional curvature is non-negative ,because we rescale by R and apply Hamilton-Ivey Pinching*
4) *$R \geq 0$,because $R_{min}$ is non decreasing .*

*Now suppose that $T < \infty$ and is maximal and therefore $\sup|Riem(x,t)| \to \infty$ as $t \to T$ .We can pick a sequence of points $(x_k, t_k)$ such that $t_k \to T$ and $\lambda_k \equiv R(x_k, t_k) \geq R(x,t)$ for any $x \in M$ and $t < t_k$ .Now consider the sequence $g_k(t)$ of parabolic rescaling at $(x_k, t_k)$ .it is defined at least on intervals $[-t_k \lambda_k, 0]$ ,where ,thanks to the Hamilton-Ivey pinching ,we have uniform bounds on the sectional curvature .Note that $-t_k \lambda_k \to -\infty$ .and by additional hypotheses we can take a limit and get a Ricci flow on $(-\infty, 0]$.note that using theorem of M.gromov we can define pointed convergence of the Ricci flow .*

*Let $(a,b)$ be an interval such that $-\infty \leq 0 \leq b \leq +\infty$ .Let $(M_k, g_k(t), x_k)$ be a sequence of pointed Ricci flows on $(a,b)$ with $x_k \in M_k$ .One says that the sequence converges in the pointed topology to the pointed Ricci flow $(M, g(t), x)$, if there exist, for any k ,embeddings $\mathcal{F}_k$ from $B(x, 0, k)$ into $M_k$ ,taking $x$ to $x_k$ ,such that the pull-back metrics $\mathcal{F}_k^* g_k(t)$ converges to $g(t)$ in the $C^\infty - Topology$ uniformly on any compact set of $M \times (a,b)$.*

# Maximum Principle

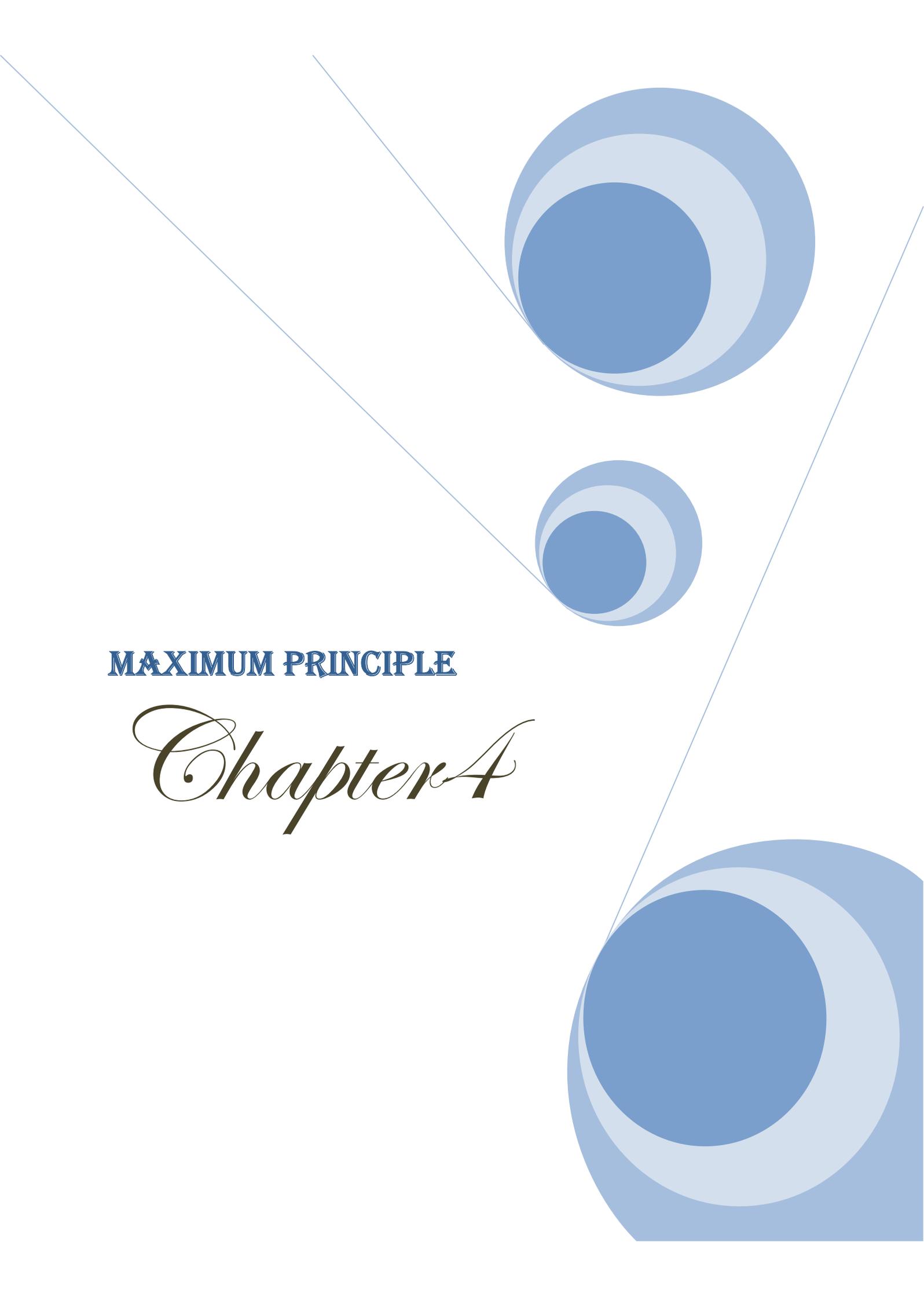

## Chapter 4

# Chapter 4

## Maximum principle

*In this section we will look at some basic PDE techniques and apply them to the Ricci flow to obtain some important results about preservation and pinching of Curvature quantities. The important fact is that the curvatures satisfy certain Reaction - diffusion equations which can be studied with the maximum principle. The essence of parabolic maximum principles is to compare the solution of a heat- type PDE by an ODE obtained by dropping the Laplacian and gradient term. The Solution to the ODE, which is easier to obtain, will act as a barrier to the PDE solution.* .Our exposition is based on [20-32]                    .

*Recall that if a smooth function* $u : U \rightarrow R$ *where* $U \subseteq R^n$ *has a local minimum at* $x_0$ *in the interior of U, Then*

$$\frac{\partial u}{\partial x^i}(x_0) = 0$$

$$\frac{\partial^2 u}{\partial x^i \partial x^j}(x_0) \geq 0$$

*Where the second statement is that the Hessian is nonnegative definite (has all nonnegative eigenvalues). The same is true on a Riemannian manifold, replacing regular derivatives with covariant derivatives.*

*Lemma4.1 Let $(M, g)$ be a Riemanian manifold and $u : M \to R$ be a smooth (or at least $C^2$) function that has a local minimum at $x_0 \in M$. Then*

$$\nabla_i u(x_0) = 0$$

$$\nabla_i \nabla_j u(x_0) \geq 0$$

$$\Delta u(x_0) = g^{ij}(x_0) \nabla_i \nabla_j u(x_0) \geq 0$$

*Proof. In a coordinate patch, the first statement is clear, since $\nabla_i u = \frac{\partial u}{\partial x^i}$ the second statement is that the Hessian is positive definite, Recall that in coordinates, the Hessian is*

$$\nabla_i \nabla_j u = \frac{\partial^2 u}{\partial x^i \partial x^j} - \Gamma_{ij}^k \frac{\partial u}{\partial x^k}$$

*But at a minimum, the second term is zero and the positive definiteness follows from the case in $R^n$. The last statement is true since both $g$ and the Hessian are positive definite.*

*Note: there is a similar statement for maxima.*

*Lemma4.2. Let $(M, g(t))$ be a smooth family of compact Riemannian manifolds for $t \in [0, T]$. let $u : [0, T] \times M \to R$ be a $C^2$ function such that*

$$u(0, x) \geq 0$$

*For all $x \in M$. Then exactly one of the following is true*

1. $u(t, x) \geq 0$ for all $(t, x) \in [0, T] \times M$, or
2. there exists a $(t_0, x_0) \in (0, T] \times M$ such that all of the following are true:
a) $u(t_0, x_0) < 0$
b) $\nabla_i u(t_0, x_0) = 0$
c) $\Delta_{g(t_0)} u(t_0, x_0) \geq 0$

d) $\frac{\partial u}{\partial t}(t_0, x_0) \leq 0$

*Proof. Certainly both cannot hold. Now suppose 1) fails. Then there must exists $(t_0, x_0)$ such that $u(t_0, x_0) < 0$. We may move this to the minimum point, at which all of the first three must hold. If we take this to be the first time that such a point occurs, the last must hold as well.*

*Theorem4.1 (pichotomy) Let $(M, g(t))$ be a smooth family of compact Riemannian Manifolds for $t \in [0, T]$. Let $u, v \in [0, T] \times M \to R$ by $C^2$ functions such that*

$$u(0, x) \geq v(0, x)$$

*For all $x \in M$. also let $A \in R$. Then exactly one of the following is true:*

1. $u(t, x) \geq v(t, x)$ *for all* $(t, x) \in [0, T] \times M$, *or*
2. *there exists a* $(t_0, x_0) \in (0, T) \subset (0, T] \times M$ *sush that all of the following are true*

a) $u(t_0, x_0) < v(t_0, x_0)$
b) $\nabla_i u(t_0, x_0) = \nabla_i v(t_0, x_0)$
c) $\Delta_{g(t_0)} u(t_0, x_0) \geq \Delta_{g(t_0)} v(t_0, x_0)$
d) $\frac{\partial u}{\partial t}(t_0, x_0) \leq \frac{\partial v}{\partial t}(t_0, x_0) + A[u(t_0, x_0) - v(t_0, x_0)]$

*Proof .replace $u$ with $e^{-At}(u - v)$.*

*Note: For elliptic equations on a manifold, the facts we use are that if a function $f$ : $M \to R$ attains its minimum at a point $x_0 \in M$, then*
$\nabla f(x_0) = 0$ *and* $\Delta f(x_0) \geq 0$

*For equations of parabolic type, a simple version says the following.*

*Theorem4.2 (weak maximum principle for super solutions of the heat equation). Let $g(t)$ be a family of metrics on a close manifold $M^n$ and let $u : M^n \times [0, T) \to R$ satisfy*

$$\frac{\partial u}{\partial t} \geq \Delta_{g(t)} u.$$

*Then if $u \geq c$ at $t = 0$ for some $c \in R$, then $u \geq c$ for all $t \geq 0$.*

*Proof. The idea is simply that given a time $t_0 \geq 0$, if the spatial minimum of $u$ is attained at a point $x_0 \in M$, then*

$$\frac{\partial u}{\partial t}(t_0, x_0) \geq \Delta_{g(t)} u(t_0, x_0) \geq 0$$

*so that the minimum should be non decreasing. Note that at $(t_0, x_0)$ we actually have $(\nabla_i \nabla_j u) \geq 0$. More rigorously, we proceed as follows. Given any $\varepsilon > 0$, define $u_\varepsilon : M \times [0, T) \to R$,*

$$u_\varepsilon = u + \varepsilon(1 + t)$$

*since $u \geq c$ at $t = 0$, we have $u_\varepsilon > c$ at $t = 0$. Now suppose for some $\varepsilon > 0$ we have $u_\varepsilon \leq c$ somewhere in $M \times (0, T)$. Then since $M$ is closed, there exists $(t_1, x_1)$ such that $u_\varepsilon(t_1, x_1) = c$ and $u_\varepsilon(t_1, x_1) > c$ for all $x \in M$ and $t \in M$ and $t \in [0, t_1)$. we then have at $(t_1, x_1)$*

$$0 \geq \frac{\partial u_\varepsilon}{\partial t} \geq \Delta_{g(t)} u_\varepsilon + \varepsilon > 0$$

*Which is a contradiction. Hence $u_\varepsilon > c$ on $M \times [0, T)$ for all $\varepsilon > 0$ and by taking the limit as $\varepsilon \to 0$ we get $u \geq c$ on $\times [0, T)$.*

*Let us work on a closed manifold $M$ with a Riemannian metric $g(t)$ that varies with time. in this section we will consider PDE s of the form*

$$\frac{\partial u}{\partial t} = \Delta_{g(t)} u + \langle X(t), \nabla u \rangle + F(u) \quad (4.1)$$

*where $u : M \times [0, T) \to R$ is a time - dependent real - valued function on $M$, $X(t)$ is a time - dependent vector field on $M$ and $F : R \to R$. we will see many PDEs of this broad type - they consist of a Laplacian term $\Delta_{g(t)} u$ and the reaction terms $\langle X(t), \nabla u \rangle + F(u)$. we call such PDEs heat - type equations, because of the analogy with the heat equation(4.1)*

*Lemma4.3 let $(M, g(t))$ be a closed manifold with a time - dependent Riemannian metric $g(t)$. suppose that $u : M \times [0, T) \to R$ is initially non-positive (i.e. $u(x, 0) \leq 0$ for all $x \in M$) and that is satisfies the differential inequality*

$$\frac{\partial u}{\partial t} \leq \Delta_{g(t)} u + \langle X(t), \nabla u \rangle$$

*at all points $(x, t) \in M \times [0, T)$ where $u(x, t) > 0$. Then $u(x, t) \leq 0$ for all $x \in M$ and $t \in [0, T)$.*

*Proof. By applying $u_\varepsilon = u - \varepsilon(1 + t)$ in proof of previous theorem, proof will be complete.*

*The following theorem essentially tells us that out upper bound grows no faster that we would expect from the reaction term $F(u)$ in (4.1)*

## The Scalar maximum principle

*Proposition4.1 Let $(M, g(t))$ be a closed manifold with a time - dependent Riemannian metric $g(t)$. suppose that $u : M \times [0, T) \to R$ satisfies*

$$\frac{\partial u}{\partial t} \leq \Delta_{g(t)} u + \langle X(t), \nabla u \rangle + F(u)$$

$u(x, 0) \leq c$ for all $x \in M$

*For some constant $c$, where $X(t)$ is a time - dependent vector field on $M$ and $F : R \to R$ is locally Lipschitz. Suppose that $\phi : R \to R$ is the solution of the associated ODE. Which is formed by neglecting the Laplacian and gradient terms:*

$$\frac{d\phi}{dt} = F(\phi)$$

$$\phi(0) = c$$

*Then $u(x, t) \leq \phi(t)$ for all $x \in M$ and $t \in [0, T)$ such that $\phi(t)$ exists.*

*Proof. Let us set $v = u - \phi$ we know that $v(x, 0) \leq 0$ for all $x \in M$, and we desire to show that $v(x, t) \leq 0$ for all $(x, t) \in M \times [0, T)$. to do this, we fix an arbitrary $\tau \in [0, T)$ and show that $v \leq 0$ on $[0, \tau]$, for any $\tau \in [0, T)$.*

*First note that*

$$\frac{\partial u}{\partial t} = \frac{\partial(u-\phi)}{\partial t}$$

$$\leq \Delta u + \langle X, \nabla u \rangle + F(u) - F(\phi)$$

$$= \Delta v + \langle X, \nabla v \rangle + (F(v) - F(\phi)) \quad (4.2)$$

*Note that $\nabla u = \nabla v$ and $\Delta u = \Delta v$ because $\phi$ depends only on t. we now want to deal with the last term on RHS.*

*Because $M \times [0,\tau]$ is compact, there exists a constant c (dependent on $\tau$) such that $|u(x,t)| \leq c$ and $|\phi(t)| \leq c$ on $M \times [0,\tau]$. Because F is locally Lipschitz and the interval $[-c,c]$ is compact, there exists $c_1$ (also dependent on  ) such that $|F(x) - F(y)| \leq c_1 |x-y|$ for all $x, y \in [-c,c]$. Therefore, because $u, \phi \in [-c,c]$, $|F(u) - F(\phi)| \leq c_1|u - \phi| = c_1 |v|$ on $M \times [0,\tau]$. Plugging this into the evolution equation (4.2) for v , we obtain*

$$\frac{\partial v}{\partial t} \leq \Delta v + \langle X, \nabla v \rangle + c_1 |v|$$

*Now let $= e^{-c_1 t} v$ . Then we have*

$$\frac{\partial \omega}{\partial t} \leq e^{-c_1 t}(\Delta v + \langle X, \nabla v \rangle + c_1 |v| - c_1 v)$$

$$= \Delta \omega + \langle X, \nabla \omega \rangle + c_1 (|\omega| - \omega) \quad (4.3)$$

*We are going to apply previous lemma to the function $\omega$. Because $v(x,0) \leq 0$ for all $x \in M$, we have $\omega(x,0) \leq 0$ for all $x \in M$. Furthermore, if $\omega > 0$ then $|\omega| = \omega$, so the differential inequality (4.3) gives us*

$$\frac{\partial \omega}{\partial t} \leq \Delta \omega + \langle X, \nabla \omega \rangle$$

*at any point $(x,t) \in M \times [0,\tau]$ such that $\omega(x,t) \geq 0$, Hence, by previous lemma $\omega(x,t) \leq 0$ for all $(x,t) \in M \times [0,\tau]$. It follows that $v(x,t) \leq 0$, and hence*

$v(x,t) \leq \phi(t)$ for all $(x,t) \in M \times [0,\tau]$, for any $\tau \in [0,T]$. Therefore $u(x,t) \leq \phi(t)$ for any $(x,t) \in M \times [0,T)$. So proof is complete.

**The maximum principle on non-compact manifolds**

*We start with following general results of Karp - Li and Ni – Tam*

*Definition4.1 we say that $u \in H^1_{loc}(M^n \times [0,T])$ is a weak sub-solution of the heat equation $(\frac{\partial}{\partial t} - \Delta)u(x,t) = 0$ if for every non-negative $C^\infty$ function $\phi$ with compact support in $M^n \times (0,T)$ we have*

$$\int_0^T \int_{M^n} \left( u \frac{\partial \phi}{\partial t} - \nabla u \nabla \phi \right) d\mu(x) dt \geq 0$$

Theorem 4.3(see [11-30]) *Assume that the curvature of $(M^n, g(t))$, $t \in [0,T)$, are uniformly bounded. if $u$ is a weak sub-solution of the heat equation on $M^n \times [0,T]$ with $u(0,0) \leq 0$ and if*

$$\int_0^T \int_{M^n} \exp\left(-\alpha d_{g(0)}^2(x,o) u_+^2(x,t)\right) d\mu_{g(t)}(x) dt < \infty$$

*For some $\alpha > 0$, then $u \leq 0$ on $M^n \times [0,T]$ ($u_+ := max\{0,u\}$ and $d(x,o)$ denote the distance function of $x$ to a fixed point $o \in M^n$)*

*Theorem4.4(see[31]) Let $(M^n, g(t)), t \in [0,T)$, be a complete solution to the Ricci flow with bounded curvature and let $\alpha_0$, be a $(p,q)$ –tensor with*

$$|\alpha_0(x)|_{g(0)} \leq e^{A(d(x,p)+1)}$$

*For some $A < \infty$. Let $E^{p,q} := (\otimes^p T^*M) \otimes (\otimes^q T^*M)$ and suppose that*

$$F_t: E^{p,q} \to E^{p,q}$$

*Is a fiber –wise linear map with*

$$\|F_t\|_\infty = \sup_{\beta(x) \in E^{p,q}} \frac{|F(\beta(x))|_{g(t)}}{|\beta(x)|_{g(t)}} < \infty$$

*Then there exists $\beta < \infty$ and a solution $\alpha(t), t \in [0,T)$, of*

$$\frac{\partial}{\partial t}\alpha = \Delta_{g(t)}\alpha + F_t(\alpha)$$

*With $\alpha(0) = \alpha_0$ and $|\alpha|_{g(t)} \leq e^{\beta(d(x,p)+1)}$. This solution is unique among all solutions with $|\alpha|_{g(t)} \leq e^{c(d(x,p)+1)}$ for all $c < \infty$.*

*Li - Yau proved the uniqueness of solutions bounded from below under a certain lower bound assumption on the Ricci curvature.*

*Theorem4.5(see [31]) Assume that the curvature and their first derivatives of $(M^n, g(t))$, $t \in [0, T)$, are uniformly bounded. For any $a > 0$ and $A > 0$, there exists a positive function $\phi(x, t)$ and $b > 0$ such that*

$$\left(\frac{\partial}{\partial t} - \Delta\right)\phi \geq A\phi$$

*On $M^n \times [0, T)$ and*

$$\exp(a. d(o, x)) \leq \phi(x, t) \leq \exp(b. d(o, x)).$$

**The maximum Principle for tensors**

*Let M be an n - dimensional complete manifold. Consider a family of smooth metrics $g_{ij}(t)$ evolving by the Ricci flow with uniformly bounded curvature for $t \in [0, T]$ with $T < \infty$. Denote by $d_t(x, y)$ the distance between two points $x, y \in M$ with respect to the metric $g_{ij}(t)$.*

*Now we give a lemma that will be useful in next theorems*

*Lemma4.4 There exists a smooth function f on M such that $f \geq 1$ everywhere, $f(x) \to +\infty$ as $d_0(x, x_0) \to +\infty$ (for some fixed $x_0 \in M$ )*

$$|\nabla f|_{g_{ij}(t)} \leq c \text{ and } |\nabla^2 f|_{g_{ij}(t)} \leq c$$

*On $M^n \times [0, T]$ for some positive constant c.*

*Here we give one of applications of weak maximum principle.*

*Proposition4.2 if the scalar curvature R of the solution $g_{ij}(t)$, $0 \leq t \leq T$, to the Ricci flow is non-negative at $t = 0$, Then it remains so on $0 \leq t \leq T$.*

*Proof Let f be the function constructed in previous lemma, recall*

$$\frac{\partial R}{\partial t} = \Delta R + 2|Ric|^2$$

*For any small constant $\varepsilon > 0$ and large constant $A > 0$, we have*

$$\frac{\partial}{\partial t}(R + \varepsilon e^{At} f) = \frac{\partial R}{\partial t} + \varepsilon A e^{At} f$$

$$= \Delta(R + \varepsilon e^{At} f) + 2|Ric|^2 + \varepsilon e^{At}(Af - \Delta f)$$

*By choosing A large enough.*

*We claim that*

$$R + \varepsilon e^{At} f > 0 \text{ on } M^n \times [0, T]$$

*Suppose not, then there exists a first time $t_0 > 0$ and a point $x_0 \in M$ such that*

$$(R + \varepsilon e^{At} f)(x_0, t_0) = 0$$

$$\nabla(R + \varepsilon e^{At} f)(x_0, t_0) = 0$$

$$\Delta(R + \varepsilon e^{At} f)(x_0, t_0) \geq 0$$

$$\frac{\partial}{\partial t}(R + \varepsilon e^{At} f)(x_0, t_0) \leq 0$$

*Then $0 \geq \frac{\partial}{\partial t}(R + \varepsilon e^{At} f)(x_0, t_0) > \Delta(R + \varepsilon e^{At} f)(x_0, t_0) \geq 0$*

*So we have $0 > 0$ that is contradiction. So we have proved that*
*$(R + \varepsilon e^{At} f) > 0$ on $M \times [0, T]$*
*Letting $\varepsilon \to 0$, we get $R \geq 0$ on $M \times [0, T]$. so proof is complete.*

### Hamilton's maximum principle

*To introduce Hamilton's maximum principle, let us start with some basic set-up. We assume $(M, g_{ij}(x, t)), t \in [0, T]$, is a smooth complete solution to the Ricci flow with bounded curvature. Let V be an abstract vector bundle over M with a metric $h_{\alpha\beta}$, and connection $\nabla = \Gamma_{ij}^k$ compatible with h. Now we may form the Laplace $\Delta \sigma = g^{ij} \nabla_i \nabla_j \sigma$ which acts on the sections $\sigma \in \Gamma(V)$ of . Suppose $M_{\alpha\beta}(x, t)$ is a family of bounded symmetric bilinear forms on V satisfying the equation*

$$\frac{\partial}{\partial t} M_{\alpha\beta} = \Delta M_{\alpha\beta} + u^i \nabla_i M_{\alpha\beta} + N_{\alpha\beta} \quad (4.4)$$

*where $u^i(t)$ is a time - dependent uniform bounded vector field on the manifold M, and $N_{\alpha\beta} = P(M_{\alpha\beta}, h_{\alpha\beta})$ is a polynomial in $M_{\alpha\beta}$ formed by contracting products of $M_{\alpha\beta}$ with itself using the metric $h = \{h_{\alpha\beta}\}$. Hamilton established the following weak maximum principle:*

*Let $M_{\alpha\beta}$ be a bounded solution to (4.4) and suppose $N_{\alpha\beta}$ satisfies the condition that*

$$N_{\alpha\beta}v^\alpha v^\beta \geq 0 \text{ whenever } N_{\alpha\beta}v^\beta = 0$$

*For $0 \leq t \leq T$. If $M_{\alpha\beta} \geq 0$ at $t = 0$, then it remains so for $0 \leq t \leq T$. Hamilton also established a strong maximum principle for solutions to equation (4.4): Let $M_{\alpha\beta}$ be a bounded solution to (4.4) with $u^i = 0$, and $N_{\alpha\beta}$ satisfies "$N_{\alpha\beta} > 0$ whenever $M_{\alpha\beta} \geq 0$". suppose $M_{\alpha\beta} \geq 0$ at $t = 0$. then there exists an interval $0 < t < \delta$ on which the rank of $M_{\alpha\beta}$ is constant and the null space of $M_{\alpha\beta}$ is invariant under parallel translation and invariant in time and also lies in the null space of $N_{\alpha\beta}$.*

## The maximum principle for systems

*Recall that the Riemannian curvature tensor may be considered as an operator : $\Omega^2 T^*M \to \Omega^2 T^*M$. As*

$$Rm(a) := g^{kr}g^{ls}R_{ijkl}a_{sr}dx^i \wedge dx^j$$

*to simply the evolution equation of Rm we need to introduce the notion of Li algebra square $Rm^\#: \Omega^2 T^*M \to \Omega^2 T^*M$.*

*First we introduce the lie algebra structure on $\Omega^2 T^*M$ by defining its lie bracket $[.,.]: \Omega^2 T^*M \times \Omega^2 T^*M \to \Omega^2 T^*M$. For any two forms and $\beta$, we define*

$$[\alpha,\beta]_{ij} = g^{kl}(\alpha_{ik}\beta_{lj} - \beta_{ik}\alpha_{lj})$$

*Let $\{\phi^\alpha, \phi^\beta, \ldots\}$ be a basis of $\Omega^2 T^*M$, and $C_\gamma^{\alpha\beta}$ be the structure symbols defined by $[\phi^\alpha, \phi^\beta] = C_\gamma^{\alpha\beta}\phi^\gamma$. We define the lie algebra square by*

$$(Rm^\#)_{\alpha\beta} := C_\alpha^{\gamma\delta}C_\beta^{\varepsilon\xi}Rm_{\gamma\varepsilon}Rm_{\delta\xi}$$

*Lemma 4.5(Evolution of the curvature operator )*

$$\frac{\partial}{\partial t}Rm = \Delta Rm + Rm^2 + Rm^\# \quad (4.5)$$

*Where $Rm^2$ is the composition $Rm \circ Rm$ and $Rm$ is the lie algebra square.*

*Let $E_x = \pi^{-1}(x)$ ($\pi: E \to M$ where $E := \Omega^2 M^n \otimes_s \Omega^2 M^n$) be the fiber over $x$. For each $\in M$, consider the system of ODE on $E_x$ corresponding to the PDE (4.5) obtained by dropping the laplacian term :*

$$\frac{d}{dt}M = M^2 + M^\# \quad (4.6)$$

*Where $M \in E_x$ is a symmetric $N \times N$ matrix, where $N = \frac{n(n-1)}{2} = \dim So(n)$ the maximum principle for systems says the following. A set K in a vector space is said to* *be convex, if for any $, Y \in K$, we have $sX + (1 - s)Y \in K$ for all $s \in [0,1]$. A subset K of the vector bundle E is said to be invariant under parallel translation, if for every path $\gamma: [a, b] \to M$ and vector $\in K \cap E_{\gamma(a)}$, the unique parallel section $X(s) \in E_{\gamma(s)}$, $s \in [a, b]$, along $\gamma(s)$ with $X(a) = X$ is contained in K.*

*Theorem4.6 (see [20-31]): Let $g(t)$ $t \in [0, T)$ be a solution to the Ricci flow on a closed manifold $M^n$. Let $K \subseteq E$ be a subset which is invariant under parallel translation and whose intersection $K_x := K \cap E_x$ with each fiber is closed and convex. Suppose the ODE (4.6) has the property that for any $M(0) \in K$, we have $M(t) \in K$ for all $\in [0, T)$, if $Rm(0) \in K$, then $Rm(t) \in K$ for all $t \in [0, T)$.*

**Strong Maximum principle**

*In order to give an idea of the proof of the strong maximum principle we need to make some remarks on functions which are not quite differentiable*

*Definition4.2 let $f$ be a real valued function defined for all $x$ in an interval containing , we define*

$$\overline{\lim_{x \to y}} f(x) \equiv \text{limsup}_{x \to y} f(x) = \inf_{\delta > 0} \sup_{|x-y| < \delta} f(x) = \lim_{\delta \to 0} \sup_{|x-y| < \delta} f(x)$$

$$\underline{\lim}_{x \to y} f(x) \equiv \text{liminf}_{x \to y} f(x) = \sup_{\delta > 0} \inf_{|x-y| < \delta} f(x) = \lim_{\delta \to 0} \inf_{|x-y| < \delta} f(x)$$

$$\underline{\lim}_{x \to y^+} f(x) \equiv \text{liminf}_{x \to y} f(x) = \sup_{\delta > 0} \inf_{0 < x-y < \delta} f(x) = \lim_{\delta \to 0} \inf_{0 < x-y < \delta} f(x)$$

$$\overline{\lim_{x \to y^+}} f(x) \equiv \text{limsup}_{x \to y^+} f(x) = \inf_{\delta > 0} \sup_{0 < x-y < \delta} f(x) = \lim_{\delta \to 0} \sup_{0 < x-y < \delta} f(x)$$

*Definition4.3 Let $: R \to R$. The upper right and lower right Dini derivatives of $f$ at $t \in R$ are, respectively, defined by*

$$D^+ f(t) = \limsup_{h \to 0^+} \frac{f(t+h) - f(t)}{h}$$

$$D_+ f(t) = \liminf_{h \to 0^+} \frac{f(t+h) - f(t)}{h}$$

*Lemma4.6 Let $f: [a, b] \to R$ be a Lipschitz function such that $f(a) \leq 0$ and $D^+ f(t) \leq 0$ when $f \geq 0$ for $\leq t \leq b$, then $f(b) \leq 0$*

*Corollary 4.1 if $f(a) \leq g(a)$ and $D^+ f(t) \leq D_+ g(t)$ for all $a \leq t \leq b$. Then $(b) \leq g(b)$.*

*Theorem4.7 (Hamilton) Let $g(t,y)$ be a smooth function of $t \in R$ and $y \in R^k$. Let $f(t) := \sup_{y \in Y} g(t,y)$, where $Y \subseteq R^k$ is a compact set. Then $f$ is a Lipschitz function and its upper right derivative satisfies*

$$D^+ f(t) \leq \sup_{y \in Y(t)} \frac{\partial}{\partial t} g(t,y)$$

*Being $Y(t) = \{y \in Y : f(t) = g(t,y)\}$.*

*Proof: choose an arbitrary $t_0 \in R$ and a sequence $\{t_j\}_{j=1}^{\infty}$ decreasing to $t_0$ for which $\lim_{t_j \to t_0} \frac{f(t_j) - f(t_0)}{t_j - t_0}$ equals the limsup. since Y is a compact set, the maximum is attained; so for each index j we can take $y_j \in Y$ such that $f(t_j) = g(t_j, y_j)$, therefore, $\{y_i\}$ is a sequence in Y and, because of the compactness of Y, there is a subsequence convergent to some $y_0 \in Y$. We can assume (for simplicity of the notation) $y_j \to y_0$ taking limits in $f(t_j) = g(t_j, y_j)$ and using the continuity of f and g, we have $f(t_0) = g(t_0, y_0)$; so $y_0 \in Y(t_0)$. By definition of f, $g(t_0, y_*) \leq g(t_0, y_0)$ $\forall y_* \in Y$; then*

$$f(t_j) - f(t_0) = g(t_j, y_j) - g(t_0, y_0) \leq g(t_j, y_j) - g(t_0, y_j)$$

*Dividing by $t_j - t_0$ and using the mean value theorem, we obtain*

$$\frac{f(t_j) - f(t_0)}{t_j - t_0} \leq \frac{g(t_j, y_j) - g(t_0, y_j)}{t_j - t_0} = \frac{\partial}{\partial t} g(T_j, y_j)$$

*With $t_0 < T_j < t_j$. Taking limits when $t_j \to t_0$, we have*

$$\lim_{t_j \to t_0} \frac{f(t_j) - f(t_0)}{t_j - t_0} \leq \frac{\partial}{\partial t} g(t_0, y_i) \leq \sup_{y \in Y(t)} \frac{\partial}{\partial t} g(t,y)$$

*Since $t_0$ is arbitrary, we have proved the estimate on the upper right derivative of f. Moreover, since $Y(t)$ is a compact set, this supremum is attained and the above inequality shows that f has bounded first derivative and so it is a Lipschitz function.*

*Next, we state the analog result for lower right derivatives.*

*Theorem4.8 (Hamilton) Let $g(t,y)$ be a smooth function of $t \in R$ and $y \in R^k$. Let $h(t) := \inf_{y \in Y} g(t,y)$, where $Y \subseteq R^k$ is a compact set. Then h is a Lipschitz function.*

$$D_+ h(t) \geq \sup_{y \in Y(t)} \frac{\partial}{\partial t} g(t,y)$$

*Being $Y(t) = \{y \in Y : h(t) = g(t,y)\}$.*

*Corollary 4.2 if $f(a) \leq 0$ and $D^+f(t) \leq cf$ for some $c \in R$ and for $\leq t \leq b$, then $f(b) \leq 0$.*

*Proof take $= e^{-ct}f$, then*

$$D^+g(t) = \limsup_{s \to t^+} \frac{e^{-cs}f(s) - e^{-ct}f(t)}{s-t}$$

$$= \limsup_{s \to t^+} \left(\frac{e^{-cs}f(s) - e^{-ct}f(s) + e^{-ct}f(s) - e^{-ct}f(t)}{s-t}\right)$$

$$= \limsup_{\delta \to 0 \; 0<s-t<\delta} \left(f(s)\frac{e^{-cs} - e^{-ct}}{s-t} + e^{-ct}\frac{f(s) - f(t)}{s-t}\right)$$

$$\leq \lim_{\delta \to 0}\left[\sup_{0<s-t<\delta}\left(f(s)\frac{e^{-cs}-e^{-ct}}{s-t}\right) + \sup_{0<s-t<\delta}\left(e^{-ct}\frac{f(s)-f(t)}{s-t}\right)\right]$$

$$= \limsup_{s \to t^+}\left(f(s)\frac{e^{-cs}-e^{-ct}}{s-t}\right) + e^{-ct}\limsup_{s \to t^+}\left(\frac{f(s)-f(t)}{s-t}\right)$$

$$= \limsup_{s \to t^+}\left(f(s)\frac{e^{-cs}-e^{-ct}}{s-t}\right) + e^{-ct}D^+f(t)$$

*So the upper limit in the first addend above is actually a limit, and the limit of a product is the product of limits. In the second addend, we can apply the hypothesis of the corollary about $D^+f(t)$. So $D^+g(t) \leq f(t)(e^{-cs})'_{s=t} + e^{-ct}cf(t) = -f(t)ce^{-ct} + e^{-ct}cf(t) = 0$. As a result of applying lemma (4.6) to $g$, we have $g(b) \leq 0$; but $g(b) = e^{-cb}f(b)$ so $f(b) \leq 0$.*

*Let $(M,g)$ be a compact Riemannian manifold and let $f = (f^1, f^2, \dots, f^k): M \to R^k$ be a system of k functions on M. Let $U \subset R^k$ be an open subset and let $\phi: U \to R^k$ be a smooth vector field on. We let $f, g$ and $\phi$ depend on time, also consider the nonlinear heat equation*

$$\frac{\partial f}{\partial t} = -\Delta_g f + \phi \circ f \quad (4.7)$$

*With $f(0) = f_0$, and suppose that it has a solution for some time interval $0 \leq t \leq T$. before dealing with this, we need some definitions.*

*Definition 4.4 we define the tangent cone $T_z(X)$ to a closed convex set $X \subset R^k$ at a point $z \in \partial X$ as the smallest closed convex con with vertex at z which contains X. It is the* *intersection of all the closed half spaces containing X with z on the boundary of the half space.*

*Definition 4.5 we say that a linear function $l: R^k \to R$ is a support function for $X \subseteq R^k$ at $z \in \partial X$ (and write $l \in S_z X$) if*

*1. $|l| = 1$*

*2. $l(z) \geq l(x)$ for all $x \in X$ (i.e. $l(z - x) \geq 0$)*

*Remark 4.2 from the view point of the support functions, we can write*

$$x \in T_z(X) \text{ if and only if } l(x - z) \leq 0 \text{ for every } l \in S_z X$$

*Lemma 4.7 the solution of the PDE equation $\frac{df}{dt} = \phi o f$ which are in the closed convex set $X \subseteq R^k$ at $t = 0$ will remain in X if and only if $\phi(z) \in T_z(X)$ for every $z \in \partial X$.*

*Theorem 4.9 if the solution of the ODE equation $\frac{df}{dt} = \phi o f$ with $f(0) \in X$ stays in X, then the solution of the PDE (4.7) with $f(0) \in X$ stays in X (suppose X is a compact set)*

*Proof first we introduce a notation*

$$s(z) := d(z, X) = \sup\{l(z - x) : x \in \partial X, l \in S_x X\}$$

*Let $s(z) = d(z, X)$, $z \in R^k$. Given a solution $f: M \times R \to R^k$ of $\frac{\partial f}{\partial t} = -\Delta_g f + \phi o f$, we define $s(t) := \sup_{x \in M} s(f(x, t))$, so, by notation we have*

$$s(t) = \sup_{(x,q,l) \in Y} l(f(x,t) - q) \quad (4.8)$$

*Being $Y = \{(x, q, l) : x \in M, q \in \partial X, l \in S_q X\}$ a compact set. So previous Hamilton's Theorem 4.10 assures that*

$$D^+ s(t) \leq \sup_{(x,q,l) \in Y} \frac{\partial}{\partial t} l(f(x,t) - q)$$

*Where $Y(t) = \{(x, q, l) \in Y : (f(x, t) - q) = s(t)\}$ from this definition of $Y(t)$ and (4.8) ; its follows that if $(x, q, l) \in Y(t)$, then $\max_{(x,q,l) \in Y(t)} l(f(x,t) - q) = l(f(x,t) - q)$, Since l is a linear function independent of t, we have*

$$D^+ s(t) \leq \sup_{(x,q,l) \in Y(t)} \frac{\partial}{\partial t} l(f(x,t) - q) = \sup_{(x,q,l) \in Y(t)} l\left(\frac{\partial f(x,t)}{\partial t}\right) =$$

$$\sup_{(x,q,l) \in Y(t)} \{-l(\Delta_g f(x,t)) + l(\phi(f(x,t)))\} \quad (4.9)$$

*Note that the last equality is true because f is a solution of the PDE, $\frac{\partial f}{\partial t} = -\Delta_g f + \phi o f$ .By definition of (t) ,l(f(x,t)) has its maximum at x ;so*

$$l(-\Delta_g f) = -\Delta_g l(f) \leq 0 \quad (4.10)$$

*On the other hand, by hypothesis, the solution of the ODE $\frac{\partial f}{\partial t} = \phi o f$ stays in X and to easily we can prove that this means $l(\phi(z)) \leq 0$ for every $l \in S_z X$ and $\in \partial X$ . Then, in particular, $l(\phi(q)) \leq 0$. So we get*

$$l(\phi(f(x,t))) \leq l(\phi(f(x,t))) - l(\phi(q)) = l(\phi(f(x,t)) - \phi(q))$$
$$\leq |l||\phi(f(x,t)) - \phi(q)| \leq c|f(x,t) - q| = cl(f(x,t) - q)$$

*Where c is the Lipschitz constant of $\phi$ and the last equality follows by the definition of Y(t) .finally, by sub situation of the above inequality and (4.9)in (4.10) ,we obtain $D^+ s(t) \leq cs(t)$ , and , since $f(0) \in X$ , $s(t) = 0$. Then applying lemma (4.6), we conclude that $s(t) = \sup_{x \in M} s(f(x,t)) = 0$ for all time in which the solution is defined. But this shows that $f(x,t)$ remains in X*

**Applications of maximum principle**

*By applying the maximum principle to various equations and inequalities governing the evolution of curvature, we will get some preliminaries control on how R and Rm evolve .*

*Theorem4.11: Under the Ricci flow, the scalar curvature , satisfy*

$$\frac{\partial R}{\partial t} \geq \Delta R + \frac{2}{n} R^2$$

*Proof since  $|Rc|^2 \geq \frac{1}{n} R^2$ (more generally, $|a|_g^2 \geq \frac{1}{n} trace_g(a)^2$ for any 2- tensor a ) also we know that $\frac{\partial R}{\partial t} = \Delta R + \frac{2}{n} |Ric|^2$, so we get $\frac{\partial R}{\partial t} \geq \Delta R + \frac{2}{n} R^2$.*

*Now we recall if $g(t)$ be a family of metrics on a closed manifold $M^n$ and suppose $v: M^n \times [0,T) \to R$ satisfies $\frac{\partial}{\partial t} v \geq \Delta_{g(t)} v + X_t(v)$ where $X_t(v)$ denotes the directional derivative of v in the direction of some time-dependent vector field $X_t$. If $v \geq c$ at $t = 0$ for some $\in R$ , then $v \geq c$ for all $t \geq 0$.*

*Now if we set $v(x,t) = e^{-ct} u(x,t)$ we get the following lemma.*

*Lemma4.8: suppose*

$$\frac{\partial}{\partial t}u \geq \Delta_{g(t)}u + X_t(u) + cu$$

*Where $u \leq 0$. Then if $u \geq 0$ at $t = 0$ so $u \geq 0$ for all time.*

*Theorem4.12 we have*

$$R(x,t) \geq \left(\frac{1}{\min_{x \in M} R(x,0)} - \frac{2}{n}t\right)^{-1}$$

*If $\min_{x \in M} R(x,0) \neq 0$. If $\min_{x \in M} R(x,0) = 0$, then we have $R(x,t) \geq 0$.*

*Proof we know that*

$$\frac{\partial R}{\partial t} \geq \Delta R + \frac{2}{n}R^2 \quad (4.11)$$

*Now we apply a trick which is useful for obtaining sharp estimate. Let $\rho(t)$ be a solution of the ODE*

$$\frac{\partial \rho}{\partial t} = \frac{2}{n}\rho^2$$

*This ODE is obtained by replacing $\geq$ by $=$ and dropping the $\Delta$ term in(4.11). we then have*

$$\frac{\partial}{\partial t}(R - \rho) \geq \Delta(R - \rho) + \frac{2}{n}(R - \rho)(R + \rho)$$

*Let $[0,T)$ be the time interval of existence of the solution of Ricci flow. For any $0 < t < T$ we have $R + \rho \leq c$ for any constant $< \infty$. Hence*

$$\frac{\partial}{\partial t}(R - \rho) \geq \Delta(R - \rho) + \frac{2c}{n}(R - \rho)$$

*Whenever $-\rho \leq 0$. Now we choose $\rho$ so that $\rho(0) = \min_{x \in M} R(x,0)$, by exercise, since $R - \rho \geq 0$ at $t = 0$, we have $R - \rho \geq 0$ on $M \times [0,\tau]$, since $\tau < T$ is arbitrary ,it follows that $R - \rho \geq 0$ on $M \times [0,\tau]$. we also easily see that the solution to $\frac{d\rho}{dt} = \frac{2}{n}\rho^2$ is $\rho(t) = \left(\frac{1}{\rho(0)} - \frac{2}{n}t\right)^{-1}$ unless $\rho(0) = 0$ in which case $\rho(t) \equiv 0$ so proof is complete.*

*Remark4.1 in particular, if $\min_{x \in M} R(x,0) > 0$, then the maximal time interval of existence $[0,T)$ satisfies $< \frac{n}{2 \min_{x \in M} R(x,0)}$.*

*We mention some useful applications of the weak and strong maximum principle:*

*Corollary 4.3 let $M \times [0,T)$ be a Ricci flow on a compact n-dimensional manifold M and suppose that $S_{min} \leq S(.,0) \leq S_{max}$ on M .Then*

$$\frac{1}{\frac{1}{S_{min}} - \frac{2}{n}t} \leq S(.,t)$$

*And if $Ric \geq 0$ everywhere, $S(.,t) \leq \frac{1}{\frac{1}{S_{max}} - 2t}$.*

*Corollary4.4(see [32] ) let $M \times [0,T]$ be a Ricci flow on an n-dimensional manifold. Assume that the curvature operator $\hat{R}$ is everywhere non-negative definite and $M(T) \cong N \times R$, where $(N, g_N)$ is an $n-1$ dimensional Riemannian manifold . Then the splitting already existed before T and the Ricci flow $M \times [0,T]$ is of product form $(N \times R) \times [0,T]$ where $N \times [0,T]$ denotes a Ricci flow on N.*

*Corollary 4.5(see [32] ) consider a Ricci flow $M \times [0,T)$ on a compact manifold M .If the curvature operator at time 0 is everywhere non-negative definite ,then this property is preserved under the Ricci flow .In particular ,in dimension 3 non-negative sectional curvature is preserved .*

*Now we give some another applications of maximum principle.*

*Theorem4.13(see [32] )  let $(U, g(t))$ ,$0 \leq t \leq T$ be a a 3-dimensional Ricci flow with non-negative sectional curvature  with U connected but not necessarily complete and $T > 0$.If $R(p,T) = 0$ for some $p \in U$, Then $(U, g(t))$ is flat for every $t \in [0,T]$.*

*Proposition4.3 (see [32] ) Let $(U, g(t))$ ,$0 \leq t \leq T$ , is a 3-dimensional Ricci flow with non-negative sectional curvature ,with U being connected but not necessarily complete and $T > 0$ .Suppose that $(U, g(T))$ is isometric to a non-empty open subset of a cone over a Riemannian manifold. Then $(U, g(t))$ is flat for every $t \in [0,T]$.*

*Also extending flows is one of applications of the maximum principle that is important.*

*Proposition 4.4(see [32]and[20-29] ) let $(M, g(t))$ $0 \leq t < T < \infty$ , be a Ricci flow with M a compact manifold. Then either the flow extends to an interval $[0,T')$ for some $T' > T$  or $|Rm|$ is unbounded on $\times [0,T)$ .*

*Now as a last proposition we give another application of maximum principle.*

*Hamilton –Ivey pinching*

*Note that we cannot completely establish a unilateral lower bound on the individual curvatures $\lambda, \mu, \upsilon$ , but one can at least show that if one of these curvatures is large and negative , then one of the others must be extremely large and positive ,and so in regions of high curvature, the positive curvature components dominate. This is important phenomenon for Ricci flow that Ivey and Hamilton formalized it.*

*Proposition4.5 (see[20-29]) Let $(M, g(t))$ be a Ricci flow on a compact 3-dimensional manifold on some time interval $[0, T]$ . Suppose that the least Eigen value $\upsilon(t, x)$ of the Riemann curvature tensor is bounded below by -1 at times $t = 0$ and all $x \in M$.then , at all space time points $(t, x) \in [0, T] \times M$, we have the scalar curvature bounded*

$$R \geq \frac{-6}{4t + 1}$$

*And furthermore whenever one has negative curvature in the sense that $\upsilon(t, x) < 0$ , then one also has the pinching bound*

$$R \geq 2|\upsilon|(log|\upsilon| + \log(1 + t) - 3)$$

# Li-Yau-Hamilton Estimates

## Chapter 5

# Chapter 5

*Li-Yau-Hamilton Estimates*

*The classical Harnack inequality from parabolic PDE theory states that for $0 < t_1 < t_2 \leq T$ a non-negative smooth solution $u \in C^\infty(M \times [0, T])$ of the linear heat equation $\partial_t u = \Delta u$ on a closed, connected manifold M satisfies*

$$\sup_M u(.,t_1) \leq c \inf_M u(.,t_2) \quad (5.1)$$

*Where c depends on $t_1, t_2$ and the geometry of M.*

*In 1986, peter Li and Shing Tung Yau found a completely new Harnack type result, namely a point wise gradient estimate that can be integrated along a path to find a classical Harnack inequality of the form (5.1) they proved that on a manifold with $Ric \geq 0$ and convex boundary ,The differential Harnack expression*

$$H(u,t) := \frac{\partial_t u}{u} - \frac{|\nabla u|^2}{u^2} + \frac{n}{2t}$$

*Is non-negative for any positive solution u of the linear heat equation.*

*Also Richard Hamilton proved a matrix version of the Li-Yau inequality under slightly different assumptions. Hamilton then found a nonlinear analog for the Ricci flow case. In this section we will introduce Hamilton's Harnack inequality for the scalar curvature under the Ricci flow.*

*Theorem5.1 (classical Harnack inequality) let $(M, g)$ be a compact $n-$dimensional Riemannian manifold with non-negative Ricci curvature .Let f be a positive solution of the heat equation*

$$\frac{\partial f}{\partial t} = -\Delta f \quad , \quad for\ 0 < t < T$$

*Then for any two points $(\xi_1, t_1), (\xi_2, t_2) \in M \times (0, T)$ with $t_1 < t_2$, we have*

$$t_1^{\frac{n}{2}} f(\xi_1, t_1) \leq e^{\frac{\psi}{4}} t_2^{\frac{n}{2}} f(\xi_2, t_2)$$

*Where $\psi = \frac{d(\xi_1, \xi_2)^2}{t_2 - t_1}$ and $d(.,.)$ denotes the distance in $(M, g)$.*

*Proof introduce $= \log f$, then*

$$\frac{\partial L}{\partial t} = -\Delta L + |dL|^2$$

*In fact, using the expression $\Delta f = -\frac{1}{\sqrt{g}} \partial_i(\sqrt{g} g^{jk} \partial_k f)$ for the Laplacian, we have*

$$\Delta L = -\frac{1}{\sqrt{g}} \partial_j(\sqrt{g} g^{jk} \partial_k L) = -\frac{1}{\sqrt{g}} \partial_j \left(\sqrt{g} g^{jk} \frac{\partial_k f}{f}\right)$$

$$= -\frac{1}{\sqrt{g}} \partial_j(\sqrt{g} g^{jk} \partial_k f) \frac{1}{f} - \frac{1}{\sqrt{g}} (\sqrt{g} g^{jk} \partial_k f) \left(\frac{-\partial_j f}{f^2}\right)$$

$$= \frac{\Delta f}{f} + g^{jk} \frac{\partial_k f}{f} \frac{\partial_j f}{f} = \frac{\Delta f}{f} + g^{jk} \partial_k L \partial_j L = \frac{\Delta f}{f} + |dL|^2$$

*Derivating respect of t in the definition of L, we get*

$$\frac{\partial L}{\partial t} = \frac{1}{f} \frac{\partial f}{\partial t} = -\frac{\Delta f}{f} \quad \left(\text{Because by heat equation } \frac{\partial f}{\partial t} = -\Delta f\right)$$

$$= -\Delta L + |dL|^2$$

*Next define $:= \frac{\partial L}{\partial t} - |dL|^2 = -\Delta L$.*

*So we have*

$$\frac{\partial L}{\partial t} = Q + |dL|^2$$

*Now we recall Boehner's formula*

*Boehner's formula: for any smooth function on $(M, g)$, we have*

$$-\frac{1}{2} \Delta |grad f|^2 = |\nabla^2 f|^2 - \langle grad f, grad(\Delta f) \rangle + Ric(grad f, grad f).$$

*And using Bochner's formula, let us compute*

$$\frac{\partial Q}{\partial t} = \frac{\partial}{\partial t}(-\Delta L) = -\Delta\left(\frac{\partial L}{\partial t}\right) = -\Delta Q - \Delta(|dL|^2)$$

$$= -\Delta Q - 2\langle grad(\Delta L), gradL\rangle + 2|\nabla^2 L|^2 + 2Ric(gradL, gradL)$$

$$= -\Delta Q + 2\langle gradQ, gradL\rangle + 2|\nabla^2 L|^2 + 2Ric(gradL, gradL)$$

*By hypothesis, we have $ic(gradL, gradL) \geq 0$ . in order to find a lower bound for $\frac{\partial Q}{\partial t}$ we are going to use the well known inequality between the square of the norm and the trace of a symmetric tensor :*

$$|\nabla^2 L|^2 \geq \frac{1}{n}(tr\nabla^2 L)^2 = \frac{1}{n}(-\Delta L)^2 = \frac{1}{n}Q^2$$

*By this remarks together with $|\nabla^2 L| \geq 0$ we get*

$$\frac{\partial Q}{\partial t} \geq -\Delta Q + 2\langle gradL, gradQ\rangle + \frac{2}{n}Q^2 \quad (5.2)$$

*by applying the scalar maximum principle ,*

*scalar maximum principle: let M be a compact manifold, and let $g_t$ a 1- parametric family of smooth metrics on M depending smoothly on t. let $X_t$ be a family of smooth vector fields on M depending smoothly on t. Let us consider the partial differential inequations*

$$\frac{\partial f_t}{\partial t} \geq -\Delta_{g_t} f_t + \langle grad_t f_t, X_t\rangle + \phi o f_t \qquad (1)$$

*And the associated ordinary differential equation*

$$\frac{\partial h(t)}{\partial t} = \phi o h \quad \text{with} \quad h(0) = \min\{f_0(x): x \in M\} \quad (2)$$

*Then $f_x(t) \geq h(t)$ for every t in an interval $[0,T]$ where there is a solution of (1) and(2).*

*So by the scalar maximum principle, we obtain the following inequality from ()*

$$Q(x,t) \geq Q_{min}(t) \geq \varphi(t) \qquad (5.3)$$

*Where $\varphi(t)$ is the solution of the ODE*

$$\frac{\partial \varphi}{\partial t} = \frac{2}{n}\varphi^2 \quad \text{with} \quad \varphi(0) = Q_{min}(0) \quad (5.4)$$

*Since M is compact, $\int_M \Delta L = 0$, then, if f is not constant at the start, $Q = -\Delta L$ must be negative at some point at the start, then $Q(0) = Q_{min}(0) < 0$. From (5.3) and (5.4) we obtain*

$$\varphi(t) \geq -\frac{n}{2t - \frac{n}{\varphi(0)}} \geq -\frac{n}{2t} \text{ and } Q \geq -\frac{n}{2t}$$

*By of Hopf –Rinow theorem, because M is compact manifold, it is also complete. So it follows that there exists a minimal geodesic jointing any pair of points in the manifold. Hence we can take a minimal geodesic parameterized by time t joining $\xi_1$ and $\xi_2$, that is $\gamma: [t_1, t_2] \to M$, such that $\gamma(t_1) = \xi_1$, $\gamma(t_2) = \xi_2$ and $d(\xi_1, \xi_2) = L(\gamma|_{[t_1,t_2]})$.*

*Denoting $\gamma^i \equiv \gamma^i(t)$ the component functions of the geodesic, we compute (using the chain rule)*

$$\frac{dL}{dt}(\gamma(t), t) = \frac{\partial L}{\partial t}(\gamma(t), t) + \frac{\partial L}{\partial x^i}\frac{d\gamma^i(t)}{dt} = \frac{\partial L}{\partial t}(\gamma(t), t) + dL(\gamma'(t))$$

$$= \frac{\partial L}{\partial t}(\gamma(t), t) + \langle gradL, \gamma'(t) \rangle \geq -\frac{n}{2t} + |dL|^2 + \langle gradL, \gamma'(t) \rangle$$

$$\geq -\frac{n}{2t} + |dL|^2 - |\langle gradL, \gamma'(t) \rangle|$$

$$\geq -\frac{n}{2t} + |dL|^2 - |gradL||\gamma'(t)| \qquad \text{(By Cauchy –Schwarz)}$$

$$= -\frac{n}{2t} + \left(|dL| - \frac{1}{2}\left|\frac{d\gamma}{dt}\right|\right)^2 - \frac{1}{4}\left|\frac{d\gamma}{dt}\right|^2 \geq -\frac{n}{2t} - \frac{1}{4}\left|\frac{d\gamma}{dt}\right|^2$$

*Integrating along the geodesic and using the fundamental theorem of calculus.*

$$L(\gamma(t_2), t_2) - L(\gamma(t_1), t_1) = \int_{t_1}^{t_2} \frac{dL}{dt}(\gamma(t), t)dt \geq -\int_{t_1}^{t_2} \frac{n}{2t} dt - \frac{1}{4}\int_{t_1}^{t_2} \left|\frac{d\gamma}{dt}\right|^2 dt$$

$$= -\frac{n}{2}\ln\left(\frac{t_2}{t_1}\right) - \frac{1}{4}\int_{t_1}^{t_2} \left|\frac{d\gamma}{dt}\right|^2 dt$$

*So we have*

$$d(\xi_1, \xi_2) = L(\gamma|_{[t_1,t_2]}) = \int_{t_1}^{t_2} \left|\frac{d\gamma}{dt}\right| dt = \left|\frac{d\gamma}{dt}\right|(t_2 - t_1) \implies \left|\frac{d\gamma}{dt}\right| = \frac{d(\xi_1, \xi_2)}{t_2 - t_1}$$

*(Because all Riemannian geodesics are constant speed curves)*

*In conclusion, we get*

$$\int_{t_1}^{t_2} \left|\frac{d\gamma}{dt}\right|^2 dt = \int_{t_1}^{t_2} \frac{d(\xi_1,\xi_2)^2}{(t_2-t_1)^2} dt = \frac{d(\xi_1,\xi_2)^2}{(t_2-t_1)^2} \int_{t_1}^{t_2} dt = \frac{d(\xi_1,\xi_2)^2}{t_2-t_1} = \psi$$

*So we reach*

$$L(\gamma(t_2),t_2) - L(\gamma(t_1),t_1) \geq \frac{-n}{2} \ln\left(\frac{t_2}{t_1}\right) - \frac{\psi}{4}$$

*and the definition of L given at the beginning of the proof yields have* $\ln\frac{f(\gamma(t_2),t_2)}{f(\gamma(t_1),t_1)} \geq \frac{-n}{2}\ln\left(\frac{t_2}{t_1}\right) - \frac{\psi}{4}$ *, and, taking exponentials,*

$$\frac{f(\xi_2,t_2)}{f(\xi_1,t_1)} \geq \left(\frac{t_2}{t_1}\right)^{\frac{-n}{2}} e^{\frac{-\psi}{4}}$$

*So proof is complete.*

*We have seen that in the Ricci flow the curvature tensor satisfies a nonlinear heat equation, and no-negatively of the curvature operator is preserved by the Ricci flow. Roughly speaking, the Li-Yau-Hamilton estimate say the non-negativity of a certain combination of the derivatives of the curvature up to second order is also preserved by the Ricci flow.*

*We start Li-Yau's gradient estimate for the linear heat equation*

*Theorem5.2: let $(M,g)$ be a compact Riemannian manifold. Assume that on the ball $B(0,2R), Ric(M) \geq -k$. Then for any $> 1$ , we have that*

$$\text{Sup}_{B(0,2R)} \left(\frac{|\nabla u|^2}{u^2} - \alpha\frac{u_t}{u}\right) \leq \frac{c\alpha^2}{R^2}\left(\frac{\alpha^2}{\alpha^2-1} + \sqrt{k}R\right) + \frac{n\alpha^2 k}{2(\alpha-1)} + \frac{n\alpha^2}{2t} \quad (5.5)$$

*Remark (see[22])if $(M,g)$ has non-negative Ricci, letting $R \to +\infty$ in (5.5) gives the clean estimate (a Hamilton Jacobi inequality):*

$$\frac{|\nabla u|^2}{u^2} - \frac{u_t}{u} \leq \frac{n}{2t} \quad (5.6)$$

*This estimate is sharp in the sense that the equality satisfied for some $(x_0,t_0)$ implies that $(M,g)$ is isometric to $R^n$.it can also be easily checked that if $u$ is the fundamental solution on $R^n$ given by the formula $\frac{1}{(4\pi t)^{\frac{n}{2}}} exp\left(-\frac{|x|^2}{4t}\right)$then the equality holds in (5.6).Because if we set $u(x,t) = (4\pi t)^{\frac{-n}{2}} exp\left(-\frac{|x|^2}{4t}\right)$ by differentiating the function $u$ ,we get*

$$\nabla_j u = -u \frac{x_j}{2t} \text{ or } \nabla_j u + u v_j = 0, \quad (5.7)$$

where $v_j = \frac{x_j}{2t} = -\frac{\nabla_j u}{u}$

*Differentiating (5.7) we get*

$$\nabla_i \nabla_j u + \nabla_i u v_j + \frac{u}{2t} \delta_{ij} = 0 \quad (5.8)$$

*To make the expression in (5.8) symmetric in i,j , we multiply $v_i$ to (5.7) and add to (5.8) and obtain*

$$\nabla_i \nabla_j u + \nabla_i u v_j + \nabla_j u v_i + u v_i v_j + \frac{u}{2t} \delta_{ij} = 0 \quad (5.9)$$

*Taking the trace in (5.9) and using the equation $\frac{\partial u}{\partial t} = \Delta u$, we arrive at*

$$\frac{\partial u}{\partial t} + 2 \nabla u . v + u |v|^2 + \frac{n}{2t} u = 0$$

*Now if we take the optimal vector field $v = -\frac{\nabla u}{u}$, we recover the equality*

$$\frac{\partial u}{\partial t} - \frac{|\nabla u|^2}{u} + \frac{n}{2t} u = 0$$

*Hamilton obtained the following Li-Yau estimate for scalar curvature $R(x,t)$.*

*Recall that under the Ricci flow on a Riemann surface the scalar curvature satisfies the following heat type equation*

$$\left( \frac{\partial}{\partial t} - \Delta \right) R = R^2$$

*By the maximum principle , the positive of the curvature is preserved by the Ricci flow , Hamilton considered the quantity $Q = \frac{\partial}{\partial t} \log R - |\nabla \log R|^2$ and computed*

$$\frac{\partial}{\partial t} \left( Q + \frac{1}{t} \right) \geq \Delta \left( Q + \frac{1}{t} \right) + 2 \nabla \log R . \nabla \left( Q + \frac{1}{t} \right) + \left( Q - \frac{1}{t} \right) \left( Q + \frac{1}{t} \right)$$

*From the maximum principle, it follows.*

*Theorem5.3 (Hamilton[33]) let $g_{ij}(x,t)$ be a complete solution to the Ricci flow with bounded curvature on a surface M . Assume the scalar curvature of the initial metric is positive . Then*

$$\frac{\partial R}{\partial t} - \frac{|\nabla R|^2}{R} + \frac{R}{t} \geq 0$$

In n-dimensions, we have the following generalization of Hamilton's differential Harnack estimate for surfaces

*Theorem5.4 (Hamilton 1993-matrix Harnack for the Ricci flow [33]) if $(M^n, g(t))$ is a solution to the Ricci flow with non-negative curvature operator, so $R_{ijkl}U_{ij}U_{kl} \geq 0$ for all 2-forms, and either $(M^n, g(t))$ is compact or complete non-compact with bounded curvature, Then for any 1-form $W \in C^\infty(\Omega^1 M)$ and 2-form $U \in C^\infty(\Delta^2 M)$ we have*

$$M_{ij}W_i W_j + 2P_{pij}U_{pi}W_j + R_{pijq}U_{pi}U_{qj} \geq 0$$

Here the 3-tensor P is defined by

$$P_{kij} := \nabla_k R_{ij} - \nabla_i R_{kj}$$

And the symmetric two tensor M is defined by

$$M_{ij} := \Delta R_{ij} - \frac{1}{2}\nabla_i \nabla_j R + 2R_{kijl}R_{kl} - R_{ip}R_{pj} + \frac{1}{2t}R_{ij}$$

We call this Hamilton's matrix Harnack estimate for Ricci flow.

Consequently For any one-form $V_i$, we have

$$\frac{\partial R}{\partial t} + \frac{R}{t} + 2\nabla_i R.V_i + 2R_{ij}V_i V_j \geq 0$$

Because by taking $U_{ij} = \frac{1}{2}(V_i W_j - V_j W_i)$ and tracing over $W_i$ we immediately get it (where $V_i = \nabla_i f$ for some function.)

*Remark in particular by taking $V \equiv 0$ we see that the function $tR(x, t)$ is point wise non-decreasing in time*

*Theorem5.5 (Li-Yau Hamilton for the conjugate heat equation[35]): assume that $(M, g(t))$ is a solution to Ricci flow on $M \times [0, T]$. Note we consider the conjugate heat equation*

$$\left(\frac{\partial}{\partial \tau} - \Delta + R\right) u(x, t) = 0 \quad (7.10)$$

Here $= T - t$. This equation is the adjoint of the heat equation $\left(\frac{\partial}{\partial t} - \Delta\right)v = 0$, then we have

$$2f_\tau + |\nabla f|^2 - R + \frac{f}{\tau} \leq 0$$

*A fundamental application of the linear trace Harnack estimate is the following classification of eternal solutions. That is the important application of Hamilton's inequality.*

*Theorem 5.6 if $(M^2, g(t)), t \in (-\infty, \infty)$, is a complete solution to the Ricci flow with positive curvature and such that $\text{Sup}_{M^2 \times (-\infty, \infty)} R$ is attained at some point in space and time, Then $(M^2, g(t))$ is a gradient Ricci Soliton. By the classification theorem, it must be the cigar soliton*

*Proof one can compute that*

$$\frac{\partial}{\partial t} Q = \Delta Q + 2\langle \nabla \log R, \nabla Q \rangle + 2 \left| \nabla \nabla \log R + \frac{1}{2} R g \right|^2 \quad (7.11)$$

*Hence $\frac{\partial}{\partial t} Q \geq \Delta Q + \langle \nabla \log R, \nabla Q \rangle + Q^2$, where we applied the inequality $|a_{ij}|^2 \geq \frac{1}{n}(tr a)^2$ to $a = \nabla \nabla \log R + \frac{1}{2} R g$ with $n = 2$. Now we used our assumption that the solution exists on all of $(-\infty, +\infty)$. in particular, for any $\alpha \in R$, the solution exists on the interval $(\alpha, +\infty)$. Since the solution PDE*

$$\frac{dq}{dt} = q^2$$

*With $\lim_{t \to \alpha} q(t) = \infty$ is $(t) = -\frac{1}{t-\alpha}$, the maximum principle says*

$$Q = \nabla \log R + R \geq -\frac{1}{t - \alpha}$$

*For all $> \alpha$. Hence, on all of $M^2 \times (-\infty, +\infty)$, by taking $\alpha \to -\infty$, we get*

$$\nabla \log R + R \geq 0$$

*By our hypothesis, there's a point $(x_0, t_0) \in M^2 \times (-\infty, +\infty)$, such that $R(x_0, t_0) = \text{Sup}_{M^2 \times (-\infty, +\infty)} R$.*

*At $(x_0, t_0)$ we have*

$$\frac{\partial R}{\partial t} = 0 \text{ and } |\nabla R| = 0$$

*And hence $(x_0, t_0) = 0$. Since $Q \geq 0$, applying the strong maximum principle to (7.11) we see that*

$$\nabla logR + R = Q \equiv 0$$

Plugging this back into (7.11), we get

$$\nabla\nabla logR + \frac{1}{2}Rg \equiv 0$$

On $M^2 \times (-\infty, +\infty)$. This says that $g(t)$ is a gradient Ricci soliton following along $\nabla logR$ :

$$\frac{\partial g}{\partial t} = -Rg = 2\nabla\nabla logR$$

As a last theorem we extend it for each $n \geq 1$.

Theorem 5.7 (see [33],[22] )

If $(M^n, g(t)), t \in (-\infty, \infty)$, is a complete solution to the Ricci flow with non-negative curvature operator ,positive Ricci curvature and such that $\text{Sup}_{M^2 \times (-\infty, \omega)} R$ is attained at some point in space and time , Then $(M^n, g(t))$ is a steady gradient Ricci Soliton.

# Two Functional F and W of Perelman

## Chapter 6

## Chapter 6

## *Two functionals $\mathcal{F}$ and $\mathcal{W}$ of Perelman*

*In this chapter, we introduce two functionals of Perelman $\mathcal{F}$ and $\mathcal{W}$, and discuss their relations with the Ricci flow. The functional $\mathcal{F}$ can be found in the literature on the String Theory, where it describes the flow energy effective action. It was not known whether the Ricci flow is a gradient flow until Perlman showed that the Ricci flow is, in a certain sense, the gradient flow of the functional $\mathcal{F}$. In the theory of dynamical systems we can take the two functionals $\mathcal{F}$ and $\mathcal{W}$ as a Lyapunov type.*

### The $\mathcal{F}$ −functional

*Let $\mathcal{M}$ denote the space of smooth Riemannian metrics $g$ on M. We think of $\mathcal{M}$ formally as an infinite-dimensional manifold. Also M is closed manifold and the tangent space $T_g M$ consists of the symmetric covariant 2-tensors $v_{ij}$ on M. f is a function $f : M \to R$ . Let $dv$ denote the Riemannian volume density associated to a metric $g$. The functional $\mathcal{F} : \mathcal{M} \times C^\infty(M) \to R$ is given by*

$$\mathcal{F}(g,f) = \int_M (R + |\nabla f|^2)e^{-f} dV$$

*Given $v_{ij} \in T_g M$ and $h = \delta f$. The evaluation of the differential $d\mathcal{F}$ on $(v_{ij}, h)$ is written as $\delta \mathcal{F}(v_{ij}, h)$. put $v = g^{ij} v_{ij}$.*

*Theorem 6.1 (Perelman) we have*

$$\delta \mathcal{F}(v_{ij}, h) = \int_M e^{-f} \left[ -v_{ij}(R_{ij} + \nabla_i \nabla_j f) + \left(\frac{v}{2} - h\right)(2\Delta f - |\nabla f|^2 + R)\right] dV \quad (6.1)$$

*Proof first we prove* $\delta R = -\Delta v + \nabla_i \nabla_j v_{ij} - R_{ij} v_{ij}$ *.In any normal coordinates at a fixed point, we have*

$$\delta R^h_{ijl} = \frac{\partial}{\partial x^i}(\delta \Gamma^h_{jl}) - \frac{\partial}{\partial x^j}(\delta \Gamma^h_{il})$$

$$= \frac{\partial}{\partial x^i}\left[\frac{1}{2}g^{hm}(\nabla_j v_{lm} + \nabla_l v_{jm} - \nabla_m v_{jl})\right]$$

$$- \frac{\partial}{\partial x^j}\left[\frac{1}{2}g^{hm}(\nabla_i v_{lm} + \nabla_l v_{im} - \nabla_m v_{il})\right]$$

$$\delta R_{jl} = \frac{\partial}{\partial x^i}\left[\frac{1}{2}g^{im}(\nabla_j v_{lm} + \nabla_l v_{jm} - \nabla_m v_{jl})\right] - \frac{\partial}{\partial x^j}\left[\frac{1}{2}g^{im}(\nabla_i v_{lm} + \nabla_l v_{im} - \nabla_m v_{il})\right]$$

$$= \frac{1}{2}\frac{\partial}{\partial x^i}\left[\nabla_j v^i_l + \nabla_l v^i_j - \nabla^i v_{jl}\right] - \frac{1}{2}\frac{\partial}{\partial x^j}[\nabla_l v]$$

*Therefore*

$$\delta R = \delta(g^{il} R_{jl}) = -v_{jl} R_{jl} + g^{jl} \delta R_{jl}$$

$$= -v_{jl} R_{jl} + \frac{1}{2}\frac{\partial}{\partial x^i}\left[\nabla^l v^i_l + \nabla_l v^{il} - \nabla^i v\right] - \frac{1}{2}\frac{\partial}{\partial x^j}[\nabla^j v]$$

$$= -v_{jl} R_{jl} + \nabla_i \nabla_l v_{il} - \Delta v$$

As $|\nabla f|^2 = g^{ij} \nabla_i f \nabla_j f$

We have $\delta |\nabla f|^2 = -v^{ij} \nabla_i f \nabla_j f + 2\langle \nabla f, \nabla h \rangle$, as $dV = \sqrt{\det(g)}\, dx_1 \ldots dx_n$, we have $\delta(dV) = \frac{v}{2} dV$, so

$$\delta(e^{-f} dV) = \left(\frac{v}{2} - h\right) e^{-f} dV \quad (6.2)$$

*Putting this together gives*

$$\delta \mathcal{F} = \int_M e^{-f}\left[-\Delta v + \nabla_i \nabla_j v_{ij} - R_{ij} v_{ij} - v_{ij} \nabla_i f \nabla_j f + 2\langle \nabla f, \nabla h \rangle + (R + |\nabla f|^2)\left(\frac{v}{2} - h\right)\right] dV \quad (6.3)$$

*The goal now is to rewrite the right-hand side of (6.3) so that* $v_{ij}$ *and h appear algebraically, i.e. without derivatives .as*

$$\Delta e^{-f} = (|\nabla f|^2 - \Delta f) e^{-f}$$

*We have*

$$\int_M e^{-f}[-\Delta v] dV = -\int_M (\Delta e^{-f}) v\, dV = \int_M e^{-f}(\Delta f - |\nabla f|^2) v\, dV$$

*Next*

$$\int_M e^{-f}\nabla_i\nabla_j v_{ij}dV = \int_M (\nabla_i\nabla_j e^{-f})v_{ij}dV$$

$$= -\int_M \nabla_i(e^{-f}\nabla_j f)v_{ij}dV = \int_M e^{-f}(\nabla_i f\nabla_j f - \nabla_i\nabla_j f)v_{ij}dV$$

*Finally*

$$2\int_M e^{-f}\langle\nabla f,\nabla h\rangle\, dV = -2\int_M \langle\nabla e^{-f},\nabla h\rangle\, dV$$

$$= 2\int_M (\Delta e^{-f})h\,dV = 2\int_M e^{-f}(|\nabla f|^2 - \Delta f)h\,dV$$

*Then*

$$\delta\mathcal{F} = \int_M e^{-f}\left[\left(\frac{v}{2}-h\right)(2\Delta f - 2|\nabla f|^2) - v_{ij}(R_{ij}+\nabla_i\nabla_j f)\right.$$
$$\left.+\left(\frac{v}{2}-h\right)(R+|\nabla f|^2)\right]dV$$

$$= \int_M e^{-f}\left[-v_{ij}(R_{ij}+\nabla_i\nabla_j f) + \left(\frac{v}{2}-h\right)(2\Delta f - |\nabla f|^2 + R)\right]dV$$

*So proof is complete.*

*Remark 6.1*(see [38-39]) *notice that*

1) *the functional $\mathcal{F}$ is invariant under diffeomorphism i.e. $\mathcal{F}(M,\phi^*g, f\circ\phi) = \mathcal{F}(M,g,\phi)$ for any diffeomorphisms $\phi$*

2) *Also for any $c > 0$ and $b$, $\mathcal{F}(M,c^2g, f+b) = c^{n-2}e^{-b}\mathcal{F}(M,g,f)$.*

   *Example (fundamental important example) let $(M,g)$ be Euclidean space $M = R^n$ and let*

$$f(t,x) = \frac{|x|^2}{4T} + \frac{n}{2}\log 4\pi T = -\log\left[(4\pi T)^{\frac{-n}{2}}e^{\frac{-|x|^2}{4T}}\right]$$

*Where $\tau = t_0 - t$, notice that $e^{-f}dx$ is the Gaussian measure, which solves the backward heat equation .If $t < t_0$ .This choice of $g$ and $f$ satisfy the equations*

$$\frac{\partial g}{\partial t} = -2Rc(g)$$

$$\frac{\partial f}{\partial t} = -\Delta f - R + |\nabla f|^2$$

*We can check*

$$\frac{\partial f}{\partial t} = \frac{\partial}{\partial t}\left[\frac{|x|^2}{4T} + \frac{n}{2}\log 4\pi T\right] = \frac{|x|^2}{4T^2} - \frac{n}{2T}$$

*And*

$$\nabla f = \frac{x}{2T}$$

*So $|\nabla f|^2 = \frac{|x|^2}{4T^2}$ and $\Delta f = \frac{n}{2T}$.*

*Also because we know $\int_{R^n} e^{\frac{-|x|^2}{4T}} dV = (4\pi T)^{\frac{n}{2}}$. By differentiating with respect to $\tau$ gives*

$$\int_{R^n} \frac{|x|^2}{4T^2} e^{\frac{-|x|^2}{4T}} dV = (4\pi T)^{\frac{n}{2}} \frac{n}{2T}.$$

*So*

$$\int_{R^n} |\nabla f|^2 e^{-f} dv = \frac{n}{2T}$$

*Then $\mathcal{F}(t) = \frac{n}{2T} = \frac{n}{2(t_0-t)}$. In particular, this is non-decreasing as a function of $t \in [0, t_0)$.*

*Theorem 6.2 let $g_{ij}$ and $f(t)$ evolve according to the coupled flow*

$$\frac{\partial g_{ij}}{\partial t} = -2R_{ij}$$

$$\frac{\partial f}{\partial t} = -\Delta f - R + |\nabla f|^2$$

*Then $\int_M e^{-f} dV$ is constant.*

*Proof by the chain rule we have*

$$\frac{\partial}{\partial t} dV = \frac{\partial}{\partial t}\sqrt{\det g_{ij}} dx^1 \wedge dx^2 \wedge \ldots \wedge dx^n$$

$$= \frac{1}{2\sqrt{\det g_{ij}}} g^{ij}(-2R_{ij}) \det g dx^1 \wedge dx^2 \wedge \ldots \wedge dx^n = -RdV$$

*Because we have $\frac{d}{dt}\det A = (A^{-1})^{ij}\left(\frac{dA_{ij}}{dt}\right)\det A$*

*Hence*

$$\frac{d}{dt}(e^{-f} dV) = e^{-f}\left(-\frac{\partial f}{\partial t} - R\right) dV = (\Delta f - |\nabla f|^2) e^{-f} dV = -\Delta(e^{-f}) dV$$

*Because $\Delta(e^{-f}) = (|\nabla f|^2 - \Delta f)e^{-f}$.*

*So it then follows that*

$$\frac{d}{dt}\int_M e^{-f}dV = -\int_M \Delta(e^{-f})dV$$

*Because M is closed manifold according to divergence theorem*

$$\int_M \Delta(e^{-f})dV = 0$$

*So this finishes the proof of the theorem.*

*we would like to get rid of the $\left(\frac{v}{2} - h\right)(2\Delta f - |\nabla f|^2 + R)$ term in (6.1). we can do this by restricting our variations so that $\frac{v}{2} - h = 0$. From (6.2), this amounts to assuming that assuming $e^{-f}dV$ is fixed. We now fix a smooth measure $dm$ on $M$ and relate $f$ to $g$ by requiring that $e^{-f}dV = dm$. Equivalently, we define a section $s: \mathcal{M} \to \mathcal{M} \times C^\infty(M)$*

*by $s(g) = \left(g, \ln\left(\frac{dV}{dm}\right)\right)$. Then the composition $\mathcal{F}^m = \mathcal{F} \circ s$ is a function on $\mathcal{M}$ and its deferential is given by*

$$d\mathcal{F}^m(v_{ij}) = \int_M e^{-f}\left[-v_{ij}(R_{ij} + \nabla_i\nabla_j f)\right]dV \quad (6.4)$$

*Defining a formal Riemannian metric on $\mathcal{M}$ by*

$$\langle v_{ij}, v_{ij}\rangle_g = \frac{1}{2}\int_M v^{ij}v_{ij}dm$$

*The gradient flow of $\mathcal{F}^m$ on $\mathcal{M}$ is given by*

$$\frac{\partial}{\partial t}g_{ij} = -2(R_{ij} + \nabla_i\nabla_j f) \quad (6.5)$$

*The induced flow equation for $f$ is*

$$\frac{\partial f}{\partial t} = \frac{1}{2}g^{ij}\frac{\partial}{\partial t}g_{ij} = -\Delta f - R$$

*As with any gradient flow, the function $\mathcal{F}^m$ is non-decreasing along the flow line with its derivative being given by the length squared of the gradient, i.e.*

$$\frac{\partial}{\partial t}\mathcal{F}^m = 2\int_M |R_{ij} + \nabla_i\nabla_j f|^2 dm$$

*We can check that if $g_{ij}(t)$ and $f(t)$ evolve according to the coupled flow*

$$\begin{cases} \frac{\partial g_{ij}}{\partial t} = -2R_{ij} \\ \frac{\partial f}{\partial t} = -\Delta f + |\nabla f|^2 - R \end{cases} \quad (6.7)$$

Then $\frac{d}{dt}\mathcal{F}(g_{ij}(t), f(t)) = 2\int_M |R_{ij} + \nabla_i\nabla_j f|^2 e^{-f}\, dV$. Now we prove that to obtaining a solution of (6.5) and (6.6) is to show that it is some how equivalent to the decoupled system of (6.7).

Proposition 6.1 (see [39]) defining $\hat{g}(t) = \sigma(t)\psi_t^*(g(t))$, we have

$$\frac{\partial \hat{g}}{\partial t} = \sigma'(t)\psi_t^*(g) + \sigma(t)\psi_t^*\left(\frac{\partial g}{\partial t}\right) + \sigma(t)\psi_t^*(\mathcal{L}_X g) \quad (6.9)$$

Also for some function $f: M \to \mathbb{R}$, we then have

$$\mathcal{L}_{(\nabla f)} g = 2 Hess(f)$$

Theorem 6.3 the solutions of (6.5) and (6.6) may be generated by pulling back solutions of (6.7) by an appropriate time-dependent diffeomorphism.

Proof By previous proposition we know $\mathcal{L}_X g = -2Hess(f)$ where $X(t) = -\nabla f$. We fix $\sigma(t) \equiv 1$. Now we define $\hat{g}(t)$ by (6.8), it will evolve, by (6.9), according to

$$\frac{\partial \hat{g}}{\partial t} = \psi_t^*(t)\left[-2Ric(g) - 2Hess_g(f)\right]$$

Where $Hess_g(f)$ is the Hessian of $f$ with respect to the metric $g$. Keeping in mind that $\psi_t: (M, \hat{g}) \to (M, g)$ is an isometry, we may then write $\hat{f} := f \circ \psi_t$ to give

$$\frac{\partial \hat{g}}{\partial t} = -2\big(Ric(\hat{g}) + Hess_{\hat{g}}(\hat{f})\big)$$

The evolution of $\hat{f}$ is then found by the chain rule

$$\frac{\partial \hat{f}}{\partial t}(x,t) = \frac{\partial f}{\partial t}(\psi_t(x), t) + X(f)(\psi_t(x), t)$$
$$= \left[(-\Delta_g f + |\nabla f|^2 - R_g) - |\nabla f|^2\right](\psi_t(x), t)$$
$$= \left[-\Delta_{\hat{g}} \hat{f} - R_{\hat{g}}\right](x, t)$$

(Here we have used of this fact that $\mathcal{L}_{\nabla f} f = |\nabla f|^2$)

Thus we have a solution to

$$\frac{\partial \hat{f}}{\partial t} = -\Delta_{\hat{g}} \hat{f} - R_{\hat{g}}$$

We would like to use this to develop a controlled quantity for Ricci flow, but we need to eliminate $f$. This can be accomplished by taking an infimum, defining

$$\lambda(M,g) = \inf_{f:\int_M e^{-f}dV=1} \mathcal{F}(M,g,f)$$

*Lemma6.1(see[39]) let $g_{ij}(t)$ and $f(t)$ evolve according to the coupled flow*

$$\begin{cases} \dfrac{\partial g_{ij}}{\partial t} = -2R_{ij} \\ \dfrac{\partial f}{\partial t} = -\Delta f + |\nabla f|^2 - R \end{cases}$$

*Then $\mathcal{F}(g_{ij}(t),f(t))$ is non-decreasing in time and monotonicity is strict unless we are on a steady gradient solution.*

*Proof according to previous computations we can show that*

$$\frac{\partial}{\partial t}\mathcal{F}(g_{ij}(t),f(t)) = 2\int_M |R_{ij} + \nabla_i\nabla_j f|^2 e^{-f}\, dV$$

*So proof is complete.*

*So by previous lemma we obtain*

$$\lambda\left(g_{ij}(t)\right) \leq \mathcal{F}(g_{ij}(t),f(t)) \leq \mathcal{F}(g_{ij}(t_0),f(t_0)) = \lambda\left(g_{ij}(t_0)\right) \quad (6.10)$$

*For $t < t_0$ and $\int_M e^{-f}dV = 1$.*

*Definition6.1 a steady breather is a Ricci flow solution on an interval $[t_1,t_2]$ that satisfies the equation $g(t_2) = \phi^* g(t_1)$ for some $\phi \in Diff(M)$.*

*Now we show that a steady breather on a compact manifold is a gradient steady soliton.*

*Theorem6.4 (Perelman[39]) a steady breather is a gradient steady soliton.*

*Proof .We have $\lambda(g(t_2)) = \lambda(\phi^* g(t_1)) = \lambda(g(t_1))$. Because $\lambda$ is invariant under diffeomorphism. So by (6.10), $\mathcal{F}(g(t),f(t))$ must be constant in t. From*

$$\frac{\partial}{\partial t}\mathcal{F}(g_{ij}(t),f(t)) = 2\int_M |R_{ij} + \nabla_i\nabla_j f|^2 e^{-f}\, dV$$

*We conclude, $R_{ij} + \nabla_i\nabla_j f = 0$. Then $R + \Delta f = 0$ and so the system*

$$\begin{cases} \dfrac{\partial g_{ij}}{\partial t} = -2R_{ij} \\ \dfrac{\partial f}{\partial t} = -\Delta f + |\nabla f|^2 - R \end{cases}$$

*Change to*

$$\frac{\partial g_{ij}}{\partial t} = -2R_{ij}$$

$$\frac{\partial f}{\partial t} = |\nabla f|^2$$

*This is a gradient expanding soliton.*

*Lemma 6.2(see [38-39])we have*

$$\frac{\partial \lambda}{\partial t} \geq \frac{2}{n}\lambda^2(t)$$

*Definition6.2 an expanding breather is a Ricci flow solution on an interval $[t_1, t_2]$ that satisfies the equation $g(t_2) = c\phi^* g(t_1)$ for some $c > 1$ and $\phi \in Diff(M)$.*

*Proposition6.2 an expanding breather is a gradient expanding solution.*

*Lemma6.3 (see [43])$\lambda(M,g)$ is finite*

*Lemma 6.4(see [43])$\lambda(M,g)$ is the least number for which one has the inequality*

$$\int_M 4|\nabla u|_g^2 + R|u|^2 dV \geq \lambda(M,g) \int_M |u|^2 dV$$

*For all u in the sobolev space $H^1(M)$.(note $\|f\|_{H^1} = \int (|\nabla f|_g^2 + f^2)dV$ for $C^1$ functions)*

### *The $\mathcal{W}$ − Functional*

*we know that the metric $g_{ij}(t)$ evolving by Ricci flow is called a breather, if for some $t_1 < t_2$ and $\alpha > 0$ the metrics $\alpha g_{ij}(t_1)$ and $g_{ij}(t_2)$ differ only by 0 diffeomorphism; the cases $\alpha = 1, \alpha < 1, \alpha > 1$ correspond to steady, shrinking and expanding breathers, respectively.*

*In order to handle the shrinking case when $\lambda(M,g) > 0$ we need to replace our functional $\mathcal{F}$ by its generalization, which contains explicit insertions of the scale parameter, to be denoted by $\tau$. Thus consider the functional.*

$$\mathcal{W}(g_{ij}, f, \tau) = \int_M [\tau(|\nabla f|^2 + R) + f - n](4\pi\tau)^{\frac{-n}{2}} e^{-f} dV$$

*Restricted to f satisfying*

$$\int_M (4\pi\tau)^{\frac{-n}{2}} e^{-f} dV = 1$$

*$\tau > 0$. Where $g_{ij}$ is a Riemannian metric, $f$ is a smooth function on M and $\tau$ is a positive scale parameter. Also for any positive number $a$ and any diffeomorphism $\varphi$*

$$\mathcal{W}(a\varphi^* g_{ij}, \varphi^* f, a\tau) = \mathcal{W}(g_{ij}, f, \tau).$$

*Now we start with first variation for $\mathcal{W}$.*

*Theorem6.5 (Perelman) assume that $\delta g_{ij} = v_{ij}$ and $\delta f = h$. Put $\sigma = \delta\tau$ then we have*

$$\delta \mathcal{W}(v_{ij}, h, \sigma) = \int_M \left[ \sigma(R + |\nabla f|^2) - \tau v_{ij}(R_{ij} + \nabla_i \nabla_j f) + h + [\tau(2\Delta f - |\nabla f|^2 + R) + f - n]\left(\frac{v}{2} - h - \frac{n\sigma}{2\tau}\right)\right](4\pi\tau)^{\frac{-n}{2}} e^{-f} dV$$

*Proof By $\frac{\partial}{\partial t} dV = \frac{v}{2} dV$. We see*

$$\delta\left((4\pi\tau)^{\frac{-n}{2}} e^{-f} dV\right) = \left(\frac{v}{2} - h - \frac{n\sigma}{2\tau}\right)(4\pi\tau)^{\frac{-n}{2}} e^{-f} dV$$

*But we know*

$$\delta \mathcal{F}(v_{ij}, h) = \int_M e^{-f}\left[-v_{ij}(R_{ij} + \nabla_i \nabla_j f) + \left(\frac{v}{2} - h\right)(2\Delta f - |\nabla f|^2 + R)\right] dV$$

*We obtain*

$$\delta \mathcal{W} = \int_M \left[\sigma(R + |\nabla f|^2) + \tau\left(\frac{v}{2} - h\right)(2\Delta f - 2|\nabla f|^2) - \tau v_{ij}(R_{ij} + \nabla_i \nabla_j f) + h \right.$$
$$\left. + [\tau(R + |\nabla f|^2) + f - n]\left(\frac{v}{2} - h - \frac{n\sigma}{2\tau}\right)\right](4\pi\tau)^{\frac{-n}{2}} e^{-f} dV$$

*But we know*

$$\Delta e^{-f} = (|\nabla f|^2 - \Delta f) e^{-f}$$

*Therefore*

$\delta \mathcal{W}(v_{ij}, h, \sigma)$

$$= \int_M \left[\sigma(R + |\nabla f|^2) - \tau v_{ij}(R_{ij} + \nabla_i \nabla_j f) + h \right.$$
$$\left. + [\tau(2\Delta f - |\nabla f|^2 + R) + f - n]\left(\frac{v}{2} - h - \frac{n\sigma}{2\tau}\right)\right](4\pi\tau)^{\frac{-n}{2}} e^{-f} dV$$

*So proof is complete.*

*Definition6.3 the arguments $g, f$ and $\tau$ are called compatible if*

$$\int_M \frac{e^{-f}}{(4\pi\tau)^{\frac{n}{2}}} dV = 1$$

*Lemma6.5(see[39]) under the transformation $(g,f,\tau) \to (\tau^{-1}g,f,1)$ compatibility is preserved, as is the functional :*

$$\mathcal{W}(g,f,\tau) = \mathcal{W}(\tau^{-1}g,f,1)$$

*Note that on $R^n$, the Gaussian measure $d\mu$ is defined in terms of the lebesgue measure $dx$ by*

$$d\mu = (2\pi)^{\frac{-n}{2}} e^{\frac{-|x|^2}{2}} dx$$

*The normalization being chosen so that*

$$\int_{R^n} d\mu = 1$$

*Lemma 6.6(see [42])(L.Gross) If $v: R^n \to R$ is, say, smooth and satisfies $v, |\nabla v| \in L^2(d\mu)$ then*

$$\int v^2 \ln|v| d\mu \leq \int |\nabla v|^2 d\mu + \left(\int v^2 d\mu\right) \ln\left(\int v^2 d\mu\right)^{\frac{1}{2}}$$

*so if we choose $v$ so that $\int v^2 d\mu = 1$, then inequality becomes*

$$\int v^2 \ln|v| d\mu \leq \int |\nabla v|^2 d\mu$$

*That called log-sobolev inequality.*

*Theorem6.6(see[39]) (Perelman) Let $g$ denote the flat metric on $R^n$. if $f$ and $\tau$ are compatible with $g$, then*

$$\mathcal{W}(g,f,\tau) \geq 0$$

*Proof let $f: R^n \to R$ be compatible with $g$ and $\tau$, which in this situation means that*

$$\int_{R^n} \frac{e^{-f}}{(4\pi\tau)^{\frac{n}{2}}} dx = 1$$

*If we set $\tau = \frac{1}{2}$ we obtain $\int_{R^n} \frac{e^{-f}}{(2\pi)^{\frac{n}{2}}} dx = 1$. If we define $v = e^{\frac{|x|^2}{4} - \frac{f}{2}}$, we have*

$$v^2 d\mu = e^{\frac{|x|^2}{2} - f} (2\pi)^{\frac{-n}{2}} e^{\frac{-|x|^2}{2}} dx = (2\pi)^{\frac{-n}{2}} e^{-f} dx$$

*So $\int_{R^n} v^2 dx = 1$. Therefore, by the log-sobolev inequality we obtain, $\int v^2 \ln|v| d\mu \leq \int |\nabla v|^2 d\mu$.*

*By computing the left-hand side and right hand side of this inequality we get*

$$\int v^2 \ln|v| d\mu = \int e^{\frac{|x|^2}{2}-f} \left(\frac{|x|^2}{4} - \frac{f}{2}\right) (2\pi)^{\frac{-n}{2}} e^{\frac{-|x|^2}{2}} dx$$

$$= \int \left(\frac{|x|^2}{4} - \frac{f}{2}\right) \frac{e^{-f}}{(2\pi)^{\frac{n}{2}}} dx$$

Also $\nabla v = \left(\frac{x}{2} - \frac{\nabla f}{2}\right) e^{\frac{|x|^2}{4} - \frac{f}{2}}$.

Which gives us

$$|\nabla v|^2 = \left(\frac{|x|^2}{4} - \frac{x.\nabla f}{2} + \frac{|\nabla f|^2}{4}\right) e^{\frac{|x|^2}{2} - f}$$

therefore

$$|\nabla v|^2 d\mu = \left(\frac{|x|^2}{4} - \frac{x.\nabla f}{2} + \frac{|\nabla f|^2}{4}\right) \frac{e^{-f}}{(2\pi)^{\frac{n}{2}}} dx$$

the integration − by − parts formula gives us

$$-\int \frac{x.\nabla f}{2} \frac{e^{-f}}{(2\pi)^{\frac{n}{2}}} dx = \frac{1}{2} \int x.\nabla(e^{-f}) \frac{dx}{(2\pi)^{\frac{n}{2}}}$$

But we can compute that $\nabla . x = n$ so

$$-\int \frac{x.\nabla f}{2} \frac{e^{-f}}{(2\pi)^{\frac{n}{2}}} dx = \frac{-n}{2} \int \frac{e^{-f}}{(2\pi)^{\frac{n}{2}}} dx$$

So

$$\int |\nabla v|^2 d\mu = \int \left(\frac{|x|^2}{4} - \frac{n}{2} + \frac{|\nabla f|^2}{4}\right) \frac{e^{-f}}{(2\pi)^{\frac{n}{2}}} dx$$

And the log-sobolev inequality gives us

$$\int \left(\frac{|x|^2}{4} - \frac{f}{2}\right) \frac{e^{-f}}{(2\pi)^{\frac{n}{2}}} dx \leq \int \left(\frac{|x|^2}{4} - \frac{n}{2} + \frac{|\nabla f|^2}{4}\right) \frac{e^{-f}}{(2\pi)^{\frac{n}{2}}} dx$$

So

$$\mathcal{W}\left(g, f, \frac{1}{2}\right) = \int \left[\frac{1}{2}|\nabla f|^2 + f - n\right] \frac{e^{-f}}{(2\pi)^{\frac{n}{2}}} dx \geq 0$$

by the scale invariance $\mathcal{W}(g,f,\tau) = \mathcal{W}\left(\frac{1}{2\tau}g, f, \frac{1}{2}\right)$ and because $(R^n, g)$ is preserved under the homothetic scaling ,we conclude

$$\mathcal{W}(g,f,\tau) \geq 0$$

and proof will be complete.

*Remark 6.2(see[38]) to easily we can check that for any $f$ , $\tau$ compatible with $g$ , $\mathcal{W}(g,f,\tau) = 0$ if and only if $f(x) \equiv \frac{x^2}{4\tau}$.*

*We have a analogous theorem for $\mathcal{W}$, like $\mathcal{F}$. We see that $\mathcal{W}$, is increasing under the Ricci flow when $f$ and $\tau$ are made to evolve appropriately.*

*Theorem6.7 (Perelman) if $g_{ij}(t), f(t)$ and $\tau(t)$ evolve according to the system*

$$\begin{cases} \frac{\partial g_{ij}}{\partial t} = -2R_{ij} \\ \frac{\partial f}{\partial t} = -\Delta f + |\nabla f|^2 - R + \frac{n}{2\tau} \\ \frac{\partial \tau}{\partial t} = -1 \end{cases}$$

*Then we have the identity*

$$\frac{d}{dt}\mathcal{W}\big(g_{ij}(t), f(t), \tau(t)\big) = \int_M 2\tau \left| R_{ij} + \nabla_i \nabla_j f - \frac{1}{2\tau}g_{ij}\right|^2 (4\pi\tau)^{\frac{-n}{2}} e^{-f} dV$$

*and $\int_M (4\pi\tau)^{\frac{-n}{2}} e^{-f} dV$ is constant .In Particular .in particular $\mathcal{W}\big(g_{ij}(t), f(t), \tau(t)\big)$ is non-decreasing in time and monotonicity is strict unless we are on a shrinking gradient soliton.*

*Now we give an example of a gradient shrinking soliton .*

*Example consider $R^n$ with the flat metric, constant in time $t \in (-\infty, 0)$ and let $\tau = -t$ and $f(\tau, x) = \frac{x^2}{4\tau}$.*

*proof because $f(\tau, x) = \frac{x^2}{4\tau}$ so we get $e^{-f} = e^{-\frac{x^2}{4\tau}}$ . To easily we can check $(g(t), f(t), \tau(t))$ satisfies the following system*

$$\begin{cases} \frac{\partial g_{ij}}{\partial t} = -2R_{ij} \\ \frac{\partial f}{\partial t} = -\Delta f + |\nabla f|^2 - R + \frac{n}{2\tau} \\ \frac{\partial \tau}{\partial t} = -1 \end{cases}$$

and $\int (4\pi\tau)^{\frac{-n}{2}} e^{-f} dV = 1$

now $\tau(|\nabla f|^2 + R) + f - n = \tau \frac{|x|^2}{4\tau^2} + \frac{|x|^2}{4\tau} - n = \frac{|x|^2}{2\tau} - n$

*so it follows from of this fact that we have*

$$\int_{R^n} e^{\frac{-|x|^2}{4\tau}} dV = (4\pi\tau)^{\frac{n}{2}}$$

*and*

$$\int e^{\frac{-|x|^2}{4\tau}} \frac{|x|^2}{4\tau^2} dV = (4\pi\tau)^{\frac{n}{2}} \frac{n}{2\tau}$$

*Therefore $\mathcal{W}(t) = 0$ for all t. So proof is complete.*

*Remark6.3 now if g is the Euclidean metric and we let $u = (4\pi\tau)^{\frac{-n}{2}} e^{-f}$, we see that*

$$\log u = \frac{-n}{2} \log(4\pi\tau) - f$$

$$|\nabla u|^2 = (4\pi\tau)^{-n} |\nabla f|^2 e^{-2f}$$

*So* $\quad |\nabla f|^2 = \frac{|\nabla u|^2}{u^2}$

*Therefore*

$$\mathcal{W}(M,g,f,\tau) = \int \left[\tau \frac{|\nabla u|^2}{u^2} - u \log u\right] dx - \frac{n}{2} \log(4\pi\tau) - n$$

*but we know that $\mathcal{W} \geq 0$ so it implies a log-sobolev inequality*

$$\tau \int \frac{|\nabla u|^2}{u^2} dx \geq \int u \log u \, dx + \frac{n}{2}(4\pi\tau) + n$$

*If we set $\phi^2 = u$*

$$4\tau \int |\nabla \phi|^2 dx \geq \frac{1}{\tau} \int \phi^2 \log \phi^2 dx + \frac{n}{2} \log(4\pi\tau) + n$$

*For the general case, we have*

$$\mathcal{W}(M,g,f,\tau) = \int \left[\tau \left(Ru + \frac{|\nabla u|^2}{u^2}\right) - u \log u\right] dV - \frac{n}{2} \log(4\pi\tau) - n \quad (6.11)$$

*one can show that*

$$\mathcal{W}(M,g,f,\tau) \geq -c(M,g,\tau)$$

*so according to (6.11) we get*

$$\tau \int R\phi^2 dV + \tau \int 4|\nabla \phi|^2 dV \geq -c + \int \phi^2 \log \phi^2 dV + \frac{n}{2} \log(4\pi\tau) + n \quad (6.12)$$

*Definition6.4 we set*

$$\mu(g_{ij},\tau) = \inf\left\{\mathcal{W}(g_{ij},f,\tau): f \in C^\infty(M), \frac{1}{(4\pi\tau)^{\frac{n}{2}}}\int_M e^{-f}dV = 1\right\}$$

*Which according to (6.12) is the best possible constant−c.*

*Remark 6.4(see [38]) μ is finite.*

*Definition6.5 a shrinking breather is a Ricci flow solution on $[t_1,t_2]$ that satisfies $g(t_2) = c\phi^*g(t_1)$ for some $c < 1$ and $\phi \in Diff(M)$.*

*Definition6.6 a shrinking soliton lives on a time interval $(-\infty, 0)$. a gradient shrinking soliton satisfies the equations*

$$\frac{\partial g_{ij}}{\partial t} = -2R_{ij} = 2\nabla_i\nabla_j f + \frac{g_{ij}}{t}$$

$$\frac{\partial f}{\partial t} = |\nabla f|^2$$

*Remark 6.5(see [39]) one can show that for any time $t \leq t_0$ we have*

$$\mu(M, g(t), \tau(t)) \leq \mathcal{W}(M, g(t), f(t), \tau(t)) \leq \mu(M, g(t_0), \tau(t_0)) \quad (6.13)$$

*Theorem6.8 A shrinking breather is a gradient shrinking soliton.*

*Proof put $t_0 = \frac{t_2 - ct_1}{1-c}$. Then if $\tau_1 = t_0 - t_1$ and $\tau_2 = t_0 - t_2$, we get $\tau_2 = c\tau_1$. Therefore because the functional $\mathcal{W}$ is invariant under simultaneous scaling of $\tau$ and $g_{ij}$ are invariant under diffeomorphism, so*

$$\mu(g(t_2),\tau_2) = \mu\left(\frac{\tau_2}{\tau_1}\phi^*g(t_1),\tau_2\right) = \mu(\phi^*g(t_1),\tau_1) = \mu(g(t_1),\tau_1)$$

*so by (6.13) and this fact that*

$$\frac{d}{dt}\mathcal{W}(g_{ij}(t),f(t),\tau(t)) = \int_M 2\tau\left|R_{ij} + \nabla_i\nabla_j f - \frac{1}{2\tau}g_{ij}\right|^2 (4\pi\tau)^{\frac{-n}{2}}e^{-f}dV$$

*it follows that the solution is a gradient shrinking soliton.*

*Remark6.6(see [43]): $\mu(g_{ij}(t),\tau - t)$ is non-decreasing along the Ricci flow.*

*Proposition6.3(see[43]) $\mu(g,\tau)$ is negative for small $\tau > 0$ and tends to zero as $\tau \to 0$.*

**Recent developments on Perelman's functional $\mathcal{F}$ and $\mathcal{W}$**

*Definition6.7 in [36] Jun-Fang Li introduced the following $\mathcal{F}$ −functional,*

$$\mathcal{F}_k(g,f) = \int_M (kR + |\nabla f|^2)e^{-f}d\mu$$

*Where $k \geq 1$. when $k = 1$, this is the $\mathcal{F}$ −functional.*

*The following theorem is analogous result, like, $\mathcal{F}$.*

*Theorem6.9 (see[36])suppose the Ricci flow of $g(t)$ exists for $[0,T)$, then all the functional $\mathcal{F}_k(g,f)$ will be monotone under the following coupled system, i.e.*

$$\begin{cases} \frac{\partial g_{ij}}{\partial t} = -2R_{ij} \\ \frac{\partial f}{\partial t} = -\Delta f + |\nabla f|^2 - R \end{cases}$$

$$\frac{d}{dt}\mathcal{F}_k(g_{ij},f) = 2(k-1)\int_M |Rc|^2 e^{-f} d\mu + 2\int_M |R_{ij} + \nabla_i\nabla_j f|^2 e^{-f} d\mu \geq 0$$

*define $\lambda_k(g) = \inf \mathcal{F}_k(g,f)$, where infimum is taken over all smooth $f$, satisfying $\int_M e^{-f} d\mu = 1$. and we assume $\lambda_1(g) = \lambda(g)$.*

*Theorem6.10(see [36]and[46]) $\lambda(g)$ is the lowest eigenvalue of the parameter $-4\Delta + R$ and the non-decreasing of the $\mathcal{F}$ functional implies the non-decreasing of $\lambda(g)$. as an application ,Perelman was able to show that there is no non-trivial steady or expanding Ricci breathers on closed manifolds.*

*Theorem6.11(see[46]) on a compact Riemannian manifold $(M,g(t))$, where $g(t)$ satisfies the Ricci flow equation for $t \in [0,T)$, the lowest eigenvalue $\lambda_k$ of the operator $-4\Delta + kR$ is non-decreasing under the Ricci flow. The monotonicity is strict unless the metric is Ricci-flat.*

*X. D. Cao considered the eigenvalues of the operator*

$$-\Delta + \frac{R}{2}$$

*On manifolds with nonnegative cuvvature operator. He showed that the eigenvalues of these manifolds are non-decreasing along the Ricci flow.*

*Corollary  On a compact Riemannian manifold, the lowest eigenvalues of the operator $-\Delta + \frac{R}{2}$ are non-decreasing under the Ricci flow .*

*proof let $k = 2$, then $\frac{1}{4}\lambda_2$ is the lowest eigenvalue of $-\Delta + \frac{R}{2}$ and the result will follows.*

Theorem6.12 (see [46])*Let $g(t), t \in [0, T)$, be a solution to the Ricci flow on a closed Riemannian manifold $M^n$. Assume that there is a $C^1$-family of smooth functions $f(t) > 0$, which satisfy*

$$\lambda(t)f(t) = -\Delta_{g(t)}f(t) + \frac{1}{2}R_{g(t)}f(t)$$

$$\int f^2(t)d\mu_{g(t)} = 1$$

*Where $\lambda(t)$ is a function of t only. Then*

$$2\frac{d}{dt}\lambda(t) = 4\int R_{ij}\nabla^i f \nabla^j f \, d\mu + 2\int |Rc|^2 f^2 \, d\mu$$

$$= \int |R_{ij} + \nabla_i \nabla_j \varphi|^2 e^{-\varphi} d\mu + \int |Rc|^2 e^{-\varphi} d\mu \geq 0$$

**Entropy functional for diffusion operator**

*Let $(M, g)$ be a compact Riemannian manifold, $\phi \in C^2(M)$. Let*
$$L = \Delta - \nabla\phi.\nabla, \quad d\mu = e^{-\phi}dV$$
*Let*
$$u = \frac{e^{-f}}{(4\pi t)^{\frac{m}{2}}}$$
*be a positive solution of*
$$(\partial_t - L)u = 0$$
*Inspired by the work of Perelman and Ni, we have the following results.*

*Theorem 6.13(X.-D. Li 2006) let*

$$H_m(u, t) = \int_M u \log u \, d\mu - \left(\frac{m}{2}\log(4\pi t) + \frac{m}{2}\right)$$

$$\mathcal{W}(u, t) = \int_M (t|\nabla f|^2 + f - m)\frac{e^{-f}}{(4\pi t)^{\frac{m}{2}}} d\mu$$

*Then*

$$\frac{d}{dt}H_m(u,t) = -\int_M \left(L\log u + \frac{m}{2t}\right)u\,d\mu$$

$$\mathcal{W}(u,t) = \frac{d}{dt}\left(tH_m(u,t)\right)$$

*Theorem 6.14 (X.-D. Li) let u be a positive solution of the heat equation*

$$\left(\frac{\partial}{\partial t} - L\right)u = 0$$

*Suppose that*

$$Ric_{m,n}(L) := Ric + \nabla^2\phi - \frac{\nabla\phi \otimes \nabla\phi}{m-n} \geq 0$$

*Then*

$$\frac{d\mathcal{W}(u,t)}{dt} = -2\int_M \tau\left(\left|\nabla^2 f - \frac{g}{2\tau}\right|^2 u\,d\mu + Ric_{m,n}(L)(\nabla f, \nabla f)u\right)d\mu$$

$$-\frac{2}{m-n}\int_M \tau\left(\nabla\phi.\nabla f + \frac{m-n}{2\tau}\right)^2 u\,d\mu$$

*Corollary (X.-D.Li 2006) supposes that $Ric_{m,n}(L) \geq 0$ then $\tau \to \mu(\tau)$ is decreasing along the heat diffusion $(\partial_\tau - L)u = 0$.*

## Perelman functional $\mathcal{F}$ and $\mathcal{W}$ for extend Ricci flow system

*We consider a system of evolution equations such that the stationary points satisfy*

$$Rc(g) = 2du \otimes du$$

$$\Delta^g u = 0$$

*to this end* Bernhard list *in [37] extended the Ricci flow to the system*

$$\frac{\partial g(t)}{\partial t} = -2Rc(g(t)) + 4du \otimes du \qquad (6.14)$$

$$\frac{\partial u(t)}{\partial t} = \Delta^{g(t)} u(t)$$

*For a Riemannian metric $g(t)$, a function $u(t)$, and given initial data $g(0)$ and $u(0)$. this is a quasilinear, weakly parabolic, coupled system of second order. Here $du \otimes du$*

is the tensor $\partial_i u \partial_j u dx^i \otimes dx^j$ and the laplacian of a function $u$ with respect to $g$ is given by $\Delta^g u = g^{ij}(\partial_i \partial_j u - \Gamma_{ij}^k \partial_k u)$.

*Definition 6.8 (see[37]) let $\tau \in R$ be a positive real number. Then the entropy $\mathcal{W}$ of a configuration*

$$(g, u, f, \tau) \in \mathcal{M}(M) \times C^\infty(M) \times C^\infty(M) \times R^+$$

*is defined to be*

$$\mathcal{W}(g, u, f, \tau) := \int_M [\tau(S + |df|^2) + f - n](4\pi\tau)^{\frac{-n}{2}} e^{-f} dV$$

Where $S := R - 2|du|^2$ (So the evolution of the metric can then be written as $\frac{\partial g_{ij}}{\partial t} = -2S_{ij}$) now we prove that $\mathcal{W}$ is scaling invariant.

*Lemma 6.7 (see[37]) let $\alpha > 0$ be a constant and $\varphi$ be a diffeomorphism of M. then the entropy $\mathcal{W}$ is invariant under simultaneous scaling of $g$ and $\tau$ by $\alpha$ in the sense that*

$$\mathcal{W}(\alpha g, u, f, \alpha \tau) = \mathcal{W}(g, u, f, \tau)$$

*and invariant under diffeomorphisms*

$$\mathcal{W}(\varphi^* g, \varphi^* u, \varphi^* f, \tau) = \mathcal{W}(g, u, f, \tau)$$

*Proof the invariance under diffeomorphisms is clear. also we have*

$\mathcal{W}(\alpha g, u, f, \alpha \tau)$

$$= \int_M \left[ \alpha\tau \big( R(\alpha g) - 2(\alpha g)^{ij} \partial_i u \partial_j u + (\alpha g)^{ij} \partial_i f \partial_j f \big) + f \right.$$
$$\left. - n \right] (4\pi\alpha\tau)^{\frac{-n}{2}} e^{-f} \sqrt{det(\alpha g)} dx$$

$$= \int_M [\alpha\tau(\alpha^{-1} R - 2\alpha^{-1}|du|^2 + \alpha^{-1}|df|^2) + f$$
$$- n] \alpha^{\frac{-n}{2}} (4\pi\tau)^{\frac{-n}{2}} e^{-f} \alpha^{\frac{n}{2}} dV = \mathcal{W}(g, u, f, \tau)$$

*Theorem 6.15 (see[37]) let M be a closed Riemannian manifold and assume that $g, u, f$ and $\tau$ satisfy on $[0, T) \times M$ the evolution equations*

$$\partial_t g = -2S_y$$
$$\partial_t u = \Delta u$$
$$\partial_t f = -\Delta f + |\nabla f|^2 - S + \frac{n}{2\tau}$$

$$\partial_t \tau = -1$$

*(S is the trace of symmetric tensor field $S_y$ )*

*then the following monotonicity formula holds:*

$$\partial_t \mathcal{W}(t) = \int_M \left( 2\tau \left| S_y + \nabla^2 f - \frac{1}{2\tau} g \right|^2 + 4\tau |\nabla u - du(\nabla f)|^2 \right) dm \geq 0$$

*Remark6.7(see[37]) the entropy $\mathcal{W}$ is non-decreasing and equality holds if and only if the solution is a homothetic shrinking gradient soliton. In this case $(g, u, f, \tau)(t)$ satisfies at every $t \in [0, T)$.*

$$S_y + \nabla^2 f - \frac{1}{2\tau} g = 0 \text{ and } \nabla u - du(\nabla f) = 0$$

*Theorem6.16(see[37]) if we vary $\mathcal{W}$ along the variation given by the following evolution equations*

$$\partial_t g = -2S_y - 2\nabla^2 f$$
$$\partial_t u = \Delta u - \langle du, df \rangle$$
$$\partial_t f = -\Delta f - S + \frac{n}{2\tau}$$
$$\partial_t \tau = -1$$

*then we have*

$$\partial_t \mathcal{W}(g, u, f, \tau)(t) = \int_M \left( 2\tau \left| S_y + \nabla^2 f - \frac{1}{2\tau} g \right|^2 + 4\tau |\Delta u - \langle du, df \rangle|^2 \right) dm \geq 0$$

*Definition 6.9(see[37]) let $(g, u, \tau) \in \mathcal{M}(M) \times C^\infty(M) \times R^+$ be a configuration .then we define*

$$\mu := \mu(g, u, \tau) = \inf_{f \in C^\infty(M)} \left\{ \mathcal{W}(g, u, f, \tau) \mid \int_M (4\pi\tau)^{\frac{-n}{2}} e^{-f} dV = 1 \right\}$$

*We investigate the remaining case of shrinking breathers.*

*Theorem6.17. (see[37]) suppose $(g, u)(t)$ is a solution to (6.14) on $[0, T) \times M$ where $M$ is closed. Fix a $\bar{\tau} \in [0, T)$ and define $\tau(t) := \bar{\tau} - t$. then $\mu(g, u, \tau)(t)$ is non-decreasing in $t$. If $\frac{\partial}{\partial t} \mu(t) = 0$ the solution is a gradient shrinking soliton.*

*Proposition6.4(see[37]) Let $(g, u)(t)$ be a shrinking breather on a closed manifold $M$. Then it necessarily is a gradient shrinking soliton.*

**Perelman's Functional $\mathcal{F}$ and $\mathcal{W}$ on Ricci Yang-Mills flow**

*The Yang-Mills heat flow was first used by Atiyah (in [45])and Bott and simon Donaldson. Donaldson used it to give an analytic proof of a Theorem of Narasimhan and Seshadri. also Atiyah and Bott used the Yang-Mills heat flow to study the topology of minimal Yang Mills connections. The Yang-Mills heat flow is a gauge-theoretic heat flow; that is, it is a differential equation for a field on a principal fiber bundle. In the study of geometric evolution equations monotonic quantities have always played an important role. Here we give a Review of the Ricci Yang- Mills flow using energy functional, but as for Ricci flow the resulting equations are not a-priori gradient equations. In other words we would like to follow the ideas of Perelman in order to write the Ricci Yang-Mills flow as a gradient flow. We claim that our coupled system is the gradient flow of some functional $\mathcal{F}(g, a, f)$ analogous to that of Perelman.*

*For definition of the Ricci Yang-Mills flow at first, we recall the $U(m)$ −vector bundle.*

*Let G be a unitary group $U(m) \subseteq GL(m, R) \subseteq GL(2m, R)$, a $G$ −vector bundle is a complex vector bundle of rank m together with Hermitian metric, a smooth function which assigns to each $p \in M$ a map*

$$\langle , \rangle_p : E_p \times E_p \to C$$

Which satisies the axioms

1. $\langle v, w \rangle_p$ is complex linear in $w$ and conjugate linear in $v$.
2. $\langle v, w \rangle_p = \overline{\langle w, v \rangle_p}$
3. $\langle v, w \rangle_p \geq 0$, with equality holding only if $v = 0$

*Definition6.10(see [44]) Let P be a $U(1)$ −bundle over a compact manifold. One can choose a metric on P such that the Ricci flow equations, with the additional hypothesis that size of the fiber remains fixed, yield the Ricci Yang-Mills flow:*

$$\frac{\partial g}{\partial t} = -2Rc + F^2$$

$$\frac{\partial A}{\partial t} = -D_A^* F(A)$$

here "A" is a connection on P and $D^*$ is adjoint of the exterior derivative D. also F is a two-form on M where $F(A) = D_A A$ that DA denote the exterior covariant derivative. Recall that if G is a lie sub group of $GL(m, R)$, a $G$-vector bundle is a rank m vector bundle whose transition functions take their values in G.

Definition6.11(see [44]) Let $(M, g)$ be a Riemannian manifold and let $E \to M$ denote a principal $K$-bundle (K is a lie group) over M with connectionA. In this section $\nabla$ will always refer to the Levi-Civita connection of g. Consider the functional

$$\mathcal{F}(g, A, f) = \int_M \left(R - \frac{1}{4}|F|^2 + |\nabla f|^2\right) e^{-f} dV$$

where R is the scalar curvature of the base metric, $f \in C^\infty(M)$ and dV denotes the volume form of g.

We use the notation $\delta$ to refer to the first variation at 0 of other quantities with respect to the parameter t.

Lemma6.8(see [44]) Let $\delta g_{ij} = v_{ij}$, $\delta A_i = \alpha_i$, $\delta f = h$, then

$\delta \mathcal{F}(v, \alpha, h)$

$$= \int_M e^{-f} \left[ -v_{ij}\left(Rc_{ij} - \frac{1}{2}\eta_{ij} + \nabla_i \nabla_j f\right) - \alpha_j(\alpha^* F_j - \nabla^i f F_{ij}) \right.$$
$$\left. + \left(\frac{v}{2} - h\right)\left(2\Delta f - |\nabla f|^2 + R - \frac{1}{4}|F|^2\right) \right] dV$$

where $\mu_{ij} = F_i^k F_{kj}$

Theorem6.18 (see [44])Given $(g(t), A(t), f(t))$ a solution to following system

$$\frac{d}{dt} g_{ij} = -2Rc_{ij} + \eta_{ij} - 2\nabla_i \nabla_j f$$

$$\frac{d}{dt} A_i = -d^* F + \nabla^i f F_{ij} \qquad (6.15)$$

$$\frac{d}{dt} f = -\Delta f - R + \frac{1}{2}|F|^2$$

Then the functional $\mathcal{F}$ is monotonically increasing in t. In particular

$$\frac{d}{dt}\mathcal{F} = \int_M \left(2\left|Rc_{ij} - \frac{1}{2}\eta_{ij} + \nabla_i\nabla_j f\right|^2 + \left|d^*F - \nabla^i f F_{ij}\right|^2\right) e^{-f} dV \geq 0$$

*Definition 6.12: (we define a metric on our configuration space to be*

$$\langle (g_1, A_1), (g_2, A_2) \rangle = \int \left(2(A_1, A_2) + 2(g_1, g_2)\right) e^{-f} dV$$

*then the gradient flow of $\mathcal{F}$ becomes (6.15)*

*Remark 6.8(see[44]) under equations (6.15) we know $\frac{d}{dt}\mathcal{F}(g(t), A(t), f(t)) \geq 0$. Equality is attained precisely when $R_{ij} - \frac{1}{2}F_i^k F_{kj} + \nabla_i\nabla_j f = 0$ and $d^*F_j - \nabla_i f F_{ij} = 0$ i.e. When $(g, A)$ is a steady gradient Ricci Yang-Mills soliton.*

*Remark 6.9* The solutions to equations (6.15) are equivalent to the system

$$\frac{d}{dt} g_{ij} = -2Rc_{ij} + F_{ik}F_{jk}$$

$$\frac{d}{dt} A_i = -d^*F_i \qquad (6.16)$$

$$\frac{d}{dt} f = -\Delta f + |\nabla f|^2 - R + \frac{1}{2}|F|^2$$

*corollary under equations (6,16)*

$$\frac{d}{dt}\mathcal{F}(g(t), A(t), f(t)) \geq 0$$

*Proposition 6.5 (see[44]) there exists a unique minimizer $\overline{f}$ of $\mathcal{F}(g, a, f)$ subject to the constraint*

$$\int e^{-f} dV = 1$$

*Definition 6.13 according to this proposition we can then define*

$$\lambda(g, A) = \inf\left\{\mathcal{F}(g, A, f): f \in C^\infty, \int e^{-f} dV = 1\right\}$$

*Proposition 6.6(see[44]) if $(g(.), A(.))$ is a solution to the Ricci Yang - Mills flow, then $\lambda(g(t), A(t))$ is non - decreasing in time.*

*Remark 6.10(see[44]) the minimum value of $\lambda(g, A)$ is equal to $\lambda_1(g, A)$, where $\lambda_1(g, A)$ is the smallest eigenvalue of the elliptic operator $-4\Delta + R - \frac{1}{4}|F|^2$.*

*Then the minimizer,$f_0$, of $\mathcal{F}$ satisfies the Euler - Lagrange equation*

$$\lambda(g, A) = 2\Delta f_0 - |\nabla f_0|^2 + R - \frac{1}{4}|F|^2$$

*Definition 6.14 A solution $(g(t), A(t))$ to the Ricci Yang -Mills flow is called a breather if there exist times $t_1 < t_2$, a constant $\alpha$, and a diffeomorphism $\varphi: M \to M$ such that $g(t_2) = \alpha\varphi^* g(t_1), A(t_2) = \alpha\varphi^* A(t_1)$*

*$\alpha > 1, \alpha < 1$ and $\alpha = 1$ correspond to $(g(t), A(t))$ being a expanding, shrinking or steady breather respectively.*

*Theorem 6.19(see[44]) let $(M^n, g(t), A(t))$ be a solution to the Ricci Yang-Mills flow on a closed manifold. if there exist $t_1 < t_2$ with $\lambda(g(t_1), A(t_1)) = \lambda(g(t_2), A(t_2))$ then $(g(t), A(t))$ is a steady gradient Ricci Yang-Mills soliton, which must have $|F|^2 = 0$ and be scalar flat.*

*Definition 6.15(see[44]) we define*

$$\mathcal{W}(g, A, f, \tau) = \int_M \left(\tau\left(|\nabla f|^2 + R + \frac{1}{4}|F|^2\right) + f - n\right)(4\pi\tau)^{\frac{-n}{2}} e^{-f} dV$$

*where $\tau = T - t$ and $(g(t), A(t))$ is a solution to Ricci Yang - Mills flow which exists on a maximal time interval of the form $[0, T]$ where $T < \infty$.*

*Theorem 6.20(see[44]) let $v_{ij} = \delta g_{ij}, \delta A_i = \alpha_i, \delta f = h$ and $\delta\tau = \sigma$. Then*

$$\delta\mathcal{W}(v, \alpha, h, \sigma)$$

$$= \int_M (4\pi\tau)^{\frac{-n}{2}} e^{-f} dV \left[\sigma\left(|\nabla f|^2 + R + \frac{1}{4}|F|^2\right)\right.$$

$$- \tau v_{ij}\left(Rc_{ij} + \frac{1}{2}\eta_{ij} + \nabla_i\nabla_j f\right) - \tau\alpha_j\left(d^* F_j - \nabla^i f F_{ij}\right) + h$$

$$\left. + \left[\tau\left(2\Delta f - |\nabla f|^2 + R - \frac{1}{4}|F|^2\right) + f - n\right]\left(\frac{v}{2} - h - \frac{n\sigma}{2\tau}\right)\right]$$

*Remark 6.11(see[44]) Consider the following system of equations*

$$\frac{d}{dt} g_{ij} = -2\left(Rc_{ij} - \frac{1}{2}\eta_{ij} + \nabla_i\nabla_j f\right)$$

$$\frac{d}{dt} A_i = -d^* F_i + \nabla^i f F_{ij}$$

$$\frac{d}{dt} f = -\Delta f - R + \frac{1}{2}|F|^2 + \frac{n}{2\tau}$$

$$\frac{d}{dt}\tau = -1$$

*for a solution to this system we have*

$$\frac{d\mathcal{W}}{dt} = \int_M \left[ 2\tau \left| Rc_{ij} - \frac{1}{2\tau} g_{ij} + \nabla_i \nabla_j f \right|^2 + \tau \left| d^* F_j + \nabla^i f F_{ij} \right|^2 + \frac{1}{4} |F|^2 \right.$$
$$\left. - \frac{1}{2} \tau |\eta|^2 \right] (4\pi\tau)^{\frac{-n}{2}} e^{-f} dV$$

*Definition 6.15(see[44]) let $(M, g(t), A(t))$ be a solution to RYM-flow which exists on a maximal time interval $T < \infty$. $(M, g(t), A(t))$ is a low - energy solution if*

$$\lim_{t \to T} (T - t) |F|^2_{C^0(M_t)} = 0$$

*So according to this Definition we get the following corollary.*

*Corollary(see[44]) let $(M, g(t), A(t))$ be a low energy solution to RYM flow on $[0, T)$ then there exists $t_0 < T$ such that for all $t_0 \leq t < T$, we have*

$$\frac{d\mathcal{W}}{dt} \geq 0$$

*Definition 6.16 Given $(M, g, A)$ at $t \in R$ let*

$$\mu(g, A, \tau) = \inf_f \left\{ \mathcal{W}(g, A, f, \tau) \mid \int_M (4\pi\tau)^{\frac{-n}{2}} e^{-f} dV = 1 \right\}$$

*Corollary (see[44]) let $(M, g(t), A(t))$ be a low - energy solution to RYM flow on $[0, T)$. Then there exists $t_0 < T$ such that for all $t_0 \leq t < T$ we have $\frac{d\mu}{dt} \geq 0$.*

# Perelman Reduced Volume and Reduced Length

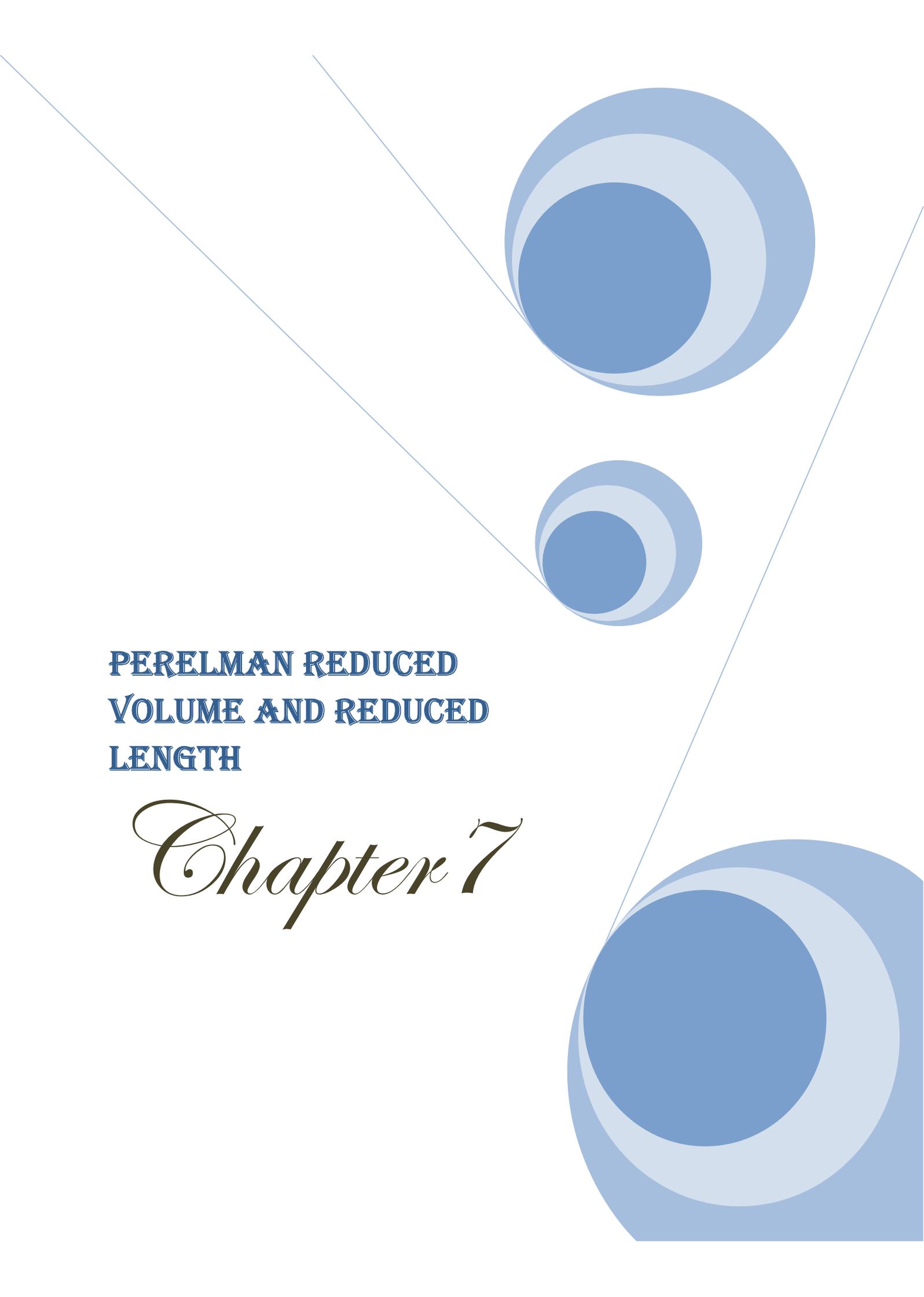

# Chapter 7

# Chapter 7

**Perelman Reduced volume and reduced length**

*In this chapter we discuss Perelman's notions of the $\mathcal{L}-$length in the context of Ricci flows. This is a functional defined on paths in space –time parametrized by backward time, denoted $\tau$. The main purpose of this chapter is to use the Li-Yau-Perelman distance to define the perelman's reduced volume, which was introduced by Perelman, and prove the monotonicity property of the reduced volume under the Ricci flow. The reduced distance, i.e. the $l-$function. The $l-$function is defined in terms of a natural curve energy along the Ricci flow ,which is analogous to the classical curve energy employed in the study of geodesics ,but involves the evolving metric .as well as the scalar curvature as a potential term. There are two applications of this theory of Perelman .We use the theory of $\mathcal{L}-$geodesics and the associated notion of reduced volume to establish non-collapsing results .The second application will be to $\kappa-$non-collapsed solutions of bounded non-negative curvature .A reader who wants to focus on the Poincare conjecture or the Geometrization conjecture could in principle start with Chapter 7.*

*Perelman reduced distance*

*Definition 7.1: suppose that either $M$ is compact or $g_{ij}(\tau)$ are complete and have uniformly bounded curvature. The $\mathcal{L}$ −length of a smooth space curve $\gamma: [\tau_1, \tau_2] \to M$ is defined by*

$$\mathcal{L}(\gamma) = \int_{\tau_1}^{\tau_2} \sqrt{\tau}\big(R(\gamma(\tau)) + |\dot{\gamma}(\tau)|^2\big) d\tau \quad (7.1)$$

*Where $\tau = \tau(t)$ satisfying $\frac{d\tau}{dt} = -1$ and the scalar curvature $R(\gamma(\tau))$ and the norm $|\dot{\gamma}(\tau)|$ are evaluated using the metric at time $t = t_0 - \tau$.*

*We consider an 1-parameter family of curves $\gamma_s : [\tau_1, \tau_2] \to M$ parameterized by $s \in (-\varepsilon, \varepsilon)$. Equivalently, We have a map $\tilde{\gamma}(s, \tau)$ with $s \in (-\varepsilon, \varepsilon)$ and $\tau \in [\tau_1, \tau_2]$. Putting $X = \frac{\partial \tilde{\gamma}}{\partial \tau}$ and $Y = \frac{\partial \tilde{\gamma}}{\partial s}$, We have $[X,Y] = 0$. (Because $\nabla_X Y - \nabla_Y X - [X,Y] = 0$). This implies that $\nabla_X Y = \nabla_Y X$.*

*Writing $\delta_Y$ as shorthand for $\frac{d}{ds}\big|_{s=0}$, and restricting to the curve $\gamma(\tau) = \tilde{\gamma}(0, \tau)$, we have $(\delta_Y Y)(\tau) = Y(\tau)$ and $(\delta_Y X)(\tau) = (\nabla_X Y)(\tau)$. Set $s = \sqrt{\tau}$ one sees immediately with respect to the variable $s$.*

*The $\mathcal{L}$ − functional is*

$$\mathcal{L}(\gamma) = \int_{\sqrt{\tau_1}}^{\sqrt{\tau_2}} \left( \frac{1}{2}\left|\frac{d\gamma}{ds}\right|^2 + 2s^2 R(\gamma(s)) \right) ds \quad (7.2)$$

*Remark 7.2 ($\mathcal{L}$ on Riemannian Products) suppose that we are given a Riemannian product solution $\big(N_1^{n_1} \times N_2^{n_2}, h_1(\tau) + h_2(\tau)\big)$ to the backward Ricci flow and a $C^1$ −Path $\gamma = (\alpha, \beta): [\tau_1, \tau_2] \to N_1 \times N_2$ so $\mathcal{L}_{h_1 + h_2}(\gamma) = \mathcal{L}_{h_1}(\alpha) + \mathcal{L}_{h_2}(\beta)$.*

*Example suppose that Ricci flow is a constant family of Euclidean metrics on $R^n \times [0, T]$. Then we have $R(\gamma(\tau)) = 0$. So according to (7.2) we get*

$$\mathcal{L}(\gamma) = \frac{1}{2} \int_0^{\bar{\tau}} \sqrt{\tau} \left|\frac{d\gamma}{d\tau}\right|^2 d\tau$$

*If we change $= \sqrt{\tau}$, $\frac{d\gamma}{d\tau} = \frac{1}{2\sqrt{\tau}} \frac{d\gamma}{ds}$, $d\tau = 2s\, ds$*

$$\mathcal{L}(\gamma) = \frac{1}{2}\int_0^{\bar{s}^2} s\frac{1}{4s^2}\left|\frac{d\gamma}{ds}\right|^2 2s\,ds = \frac{1}{2}\int_0^{\bar{s}^2}\left|\frac{d\gamma}{ds}\right|^2 ds$$

*So we arrive to (7.2) formula.*

### First variation formulae for $\mathcal{L}$ −geodesics

*At first recall the classical variation formulae that give the first derivative of the metric function $d(x_0, x)$ on a Riemannian manifold $(M, g)$.*

*Let $\gamma: [0,1] \to M$ ranges over all $C^1$ −Curves from $x_0$ to $x$, and the Dirichlet energy $E(\gamma)$ of the curve is given by the formula*

$$E(\gamma) = \frac{1}{2}\int_0^1 |X|_g^2\, dt \quad \text{where } X = \frac{\partial \gamma}{\partial t}$$

*We know that the distance $d(x_0, x)$ on a Riemannian manifold $(M, g)$ can be defined by the energy-minimisation formula*

$$\tfrac{1}{2}d(x_0, x)^2 = \inf_\gamma E(\gamma) \quad (7.3)$$

*Also if $\gamma: [0, \bar{t}] \to M$, $\gamma(0) = p$, , $\gamma(\bar{t}) = q$*

$$\inf\{E(\gamma)|\; \gamma: [0,\bar{t}] \to M, \gamma(0) = p,, \gamma(\bar{t}) = q\} = \frac{d^2(p, q)}{\bar{t}}$$

*Lemma 7.1 The first variation formula for Dirichlet energy $E(\gamma)$ is*

$$\frac{d}{ds}E(\gamma) = \langle X, Y\rangle|_{t=0}^1 - \int_0^1 \langle X, \nabla_Y Y\rangle dt \quad (7.4)$$

*Proof: consider $D_t X$ means $D_t X = \nabla_{\frac{\partial}{\partial t}} X$ and $X = \frac{\partial \gamma}{\partial s}$, $\dot{\gamma} = \frac{\partial \gamma}{\partial t}$.*

$$\delta E(\gamma) = \delta \frac{1}{2}\int_0^1 \langle \dot{\gamma}, \dot{\gamma}\rangle dt$$

$$= \frac{1}{2}\frac{\partial}{\partial s}\Big|_{s=0} \int_0^1 \langle \frac{\partial \gamma}{\partial t}, \frac{\partial \gamma}{\partial t} \rangle \, dt$$

$$= \frac{1}{2}\int_0^1 \langle D_s\left(\frac{\partial \gamma}{\partial t}\right), \frac{\partial \gamma}{\partial t} \rangle + \langle \frac{\partial \gamma}{\partial t}, D_s\left(\frac{\partial \gamma}{\partial t}\right) \rangle \, dt$$

$$= \int_0^1 \langle D_s\left(\frac{\partial \gamma}{\partial t}\right), \frac{\partial \gamma}{\partial t} \rangle \, dt$$

$$= \int_0^1 \langle D_t X, \dot{\gamma} \rangle dt$$

*With integration by parts we get*

$$= \langle X, \dot{\gamma} \rangle|_0^1 - \int_0^1 \langle X, D_t\dot{\gamma} \rangle dt$$

$$= \langle X, Y \rangle|_0^1 - \int_0^1 \langle X, \nabla_Y Y \rangle$$

*Remark7.3 if we fix the end points of γ ,then the first term of the right hand side of (7.4) vanishes .If we consider arbitrary variations X of γ with fixed end points,we thus conclude that in order to be a minimiser for (7.3),that γ must obey the geodesic flow equation*

$$\nabla_Y Y = 0$$

*Now we develop analogous variational formulae for ℒ −length (reduced distance)on Ricci flow*

*Theorem 7.1(Perelman [39]) the first variation formula for ℒ −length is*

$$\delta_Y(\mathcal{L}) = 2\sqrt{\tau}\langle X, Y \rangle\Big|_{\tau_1}^{\tau_2} + \int_{\tau_1}^{\tau_2} \sqrt{\tau} \langle Y, \nabla R - 2\nabla_X X - 4Ric(.,X) - \frac{1}{\tau}X \rangle \, d\tau \quad (7.5)$$

*Where $\langle .,. \rangle$ denotes the inner product with respect to the metric $g_{ij}(\tau)$.*

*Proof let $X(\tau) = \dot{\gamma}(\tau)$,so according to (7.1) we obtain*

$$\delta_Y(\mathcal{L}) = \int_{\tau_1}^{\tau_2} \sqrt{\tau}(\nabla_Y R + 2\langle \nabla_X Y, X \rangle) d\tau$$

$$= \int_{\tau_1}^{\tau_2} \sqrt{\tau}(\langle \nabla R, Y\rangle + 2\langle \nabla_X Y, X\rangle)d\tau$$

*Because $\frac{\partial \tau}{\partial t} = -1$ so $\frac{dg_{ij}}{d\tau} = 2R_{ij}$, using this fact we get*

$$\frac{d\langle Y,X\rangle}{d\tau} = \langle \nabla_X Y, X\rangle + \langle Y, \nabla_X X\rangle + 2Ric(Y,X) \quad (7.6)$$

*Because we can break $\frac{d\langle Y,X\rangle}{d\tau}$ into two parts :first assumes that the metric is constant and the second deals with the variation with $\tau$ of the metric .The first contribution is the usual formula*

$$\frac{d\langle Y,X\rangle}{d\tau} = \langle \nabla_X Y, X\rangle + \langle Y, \nabla_X X\rangle$$

*We show last term in equation (7.6) is that come from differentiating the metric with respect to $\tau$ and Ricci flow equation we get*

$$\frac{d\langle Y,X\rangle}{d\tau} = 2Ric(Y,X)$$

*So we arrive to equality of (7.6)*

*We continue the proof of (7.5)*

$$\int_{\tau_1}^{\tau_2} \sqrt{\tau}(\langle Y, \nabla R\rangle + 2\langle \nabla_X Y, X\rangle)d\tau =$$

$$= \int_{\tau_1}^{\tau_2} \sqrt{\tau}\left(\langle Y, \nabla R\rangle + 2\frac{d\langle Y,X\rangle}{d\tau} - 2\langle Y, \nabla_X X\rangle - 4Ric(Y,X)\right)d\tau$$

$$= 2\sqrt{\tau}\langle X, Y\rangle\Big|_{\tau_1}^{\tau_2} + \int_{\tau_1}^{\tau_2} \sqrt{\tau}\,\langle Y, \nabla R - 2\nabla_X X - 4Ric(X,.) - \frac{1}{\tau}X\rangle\,d\tau$$

*Lemma7.2 (Perelman) we consider a variation $\gamma(\tau,s)$ with fixed endpoints, so that $Y(\tau_1) = Y(\tau_2) = 0$.so from (7.5) the $\mathcal{L}$-shortest geodesic $\gamma(\tau,s)$ for $\tau \in [\tau_1, \tau_2]$ satisfies the following $\mathcal{L}$ − geodesic equation*

$$\nabla_X X - \frac{1}{2}\nabla R + \frac{1}{2\tau}X + 2Ric(X,.) = 0 \quad (7.7)$$

*This equation is called Euler-Lagrange Equation.*

*Therefore we can say for any $\tau_2 > \tau_1 > 0$, there is always an $\mathcal{L}$-geodesic $\gamma(\tau)$ for $\tau \in [\tau_1, \tau_2]$.*

*Written with respect to the variable $s = \sqrt{\tau}$ the $\mathcal{L}$-geodesic equation becomes*

$$\nabla_{\hat{X}}\hat{X} - 2s^2 \nabla R + 4s Ric(\hat{X},.) = 0$$

*Where $\hat{X} = \frac{d\gamma}{ds} = 2sX$*

*Lemma7.3(Perelman) let $\gamma: [\tau_1, \tau_2] \to M$ be an $\mathcal{L}$-geodesic. Then $\lim_{\tau \to 0} \sqrt{\tau}X(\tau)$ exists.*

*Proof: multiplying the $\mathcal{L}$-geodesic equation (7.7) by $\sqrt{\tau}$, we get*

$$\nabla_X(\sqrt{\tau}\,X) = \frac{\sqrt{\tau}}{2}\nabla R - 2\sqrt{\tau}\,Ric(X,.) \; On \; [\tau_1, \tau_2]$$

*or equivalently*

$$\frac{d}{d\tau}(\sqrt{\tau}\,X) = \frac{\sqrt{\tau}}{2}\nabla R - 2\,Ric(\sqrt{\tau}X,.) \; On \quad [\tau_1, \tau_2]$$

*Thus if a continuous curve, defined on $[\tau_1, \tau_2]$, satisfies the $\mathcal{L}$-geodesic equation on every subinterval $0 < \tau_1 \le \tau \le \tau_2$, then $\sqrt{\tau_1}X(\tau_1)$ has a limit as $\tau_1 \to 0^+$.*

*Remark7.4 for a fixed $p \in M$, by taking $\tau_1 = 0$ and $\gamma(0) = p$, the vector $v = \lim_{\tau \to 0} \sqrt{\tau}X(\tau)$ is well-defined in $T_PM$. The $\mathcal{L}$-exponential map $\mathcal{L}exp_\tau: T_PM \to M$ sends $v$ to $\gamma(\tau)$.*

*Note that for any vector $v \in T_PM$, we can find an $\mathcal{L}$-geodesic $\gamma(\tau)$ with $\lim_{\tau \to 0^+} \sqrt{\tau}\dot{\gamma}(\tau) = v$.*

*Now we give an Estimate for speed of $\mathcal{L}$-geodesics*

*Theorem(see[38])7.2: let $(M^n, g(\tau))$, $\tau \in [0, T]$, be a solution to the backward Ricci flow with bounded sectional curvature and $\max\{|Rm|, |Rc|\} \le C_0 < \infty$ on $M \times [0, T]$. There*

*exists a constant $C(n) < \infty$ depending only on $n$ such that given $0 \leq \tau_1 \leq \tau_2 < T$, if $\gamma: [\tau_1, \tau_2] \to M$ is an $\mathcal{L}$ −geodesic with*

$$\lim_{\tau \to \tau_1} \sqrt{\tau} \frac{d\gamma}{d\tau}(\tau) = v \in T_{\gamma(\tau_1)} M$$

*Then for any $\in [\tau_1, \tau_2]$,*

$$\tau \left| \frac{d\gamma}{d\tau}(\tau) \right|^2_{g(\tau)} \leq e^{6C_0 T} |v|^2 + \frac{C(n)T}{\min\{T - \tau_2, C_0^{-1}\}} (e^{6C_0 T} - 1)$$

*Where $|v|^2 = |v|^2_{g(\tau_1)}$.*

*Also we give a lemma about existence of $\mathcal{L}$ −geodesics between any two space-time endpoints.*

*Lemma7.4(see[39]) let $(M^n, g(\tau))$, $\tau \in [0, T]$, be a complete solution to the backward Ricci flow with bounded sectional curvature. Given $p, q \in M$ and $0 < \tau_1 \leq \tau_2 < T$, There exists a smooth path $\gamma(\tau): [\tau_1, \tau_2] \to M$ from $p$ to $q$ such that $\gamma$ has the minimal of $\mathcal{L}$ −length among all such paths. Furthermore, all $\mathcal{L}$ −length minimizing paths are smooth $\mathcal{L}$ −geodesics.*

*Corollary7.1 (see[43]): if we extend the curve $\gamma$ for piecewise smooth curves (where $0 \leq \tau_1 < \tau_2 \leq t_0$) then for first variation formula of breaking points $\tau_0' = \tau_1 < \tau_2' < \cdots < \tau_k' = \tau_2$, we have*

$$\delta(\mathcal{L}(\gamma)) = \int_{\tau_1}^{\tau_2} \sqrt{\tau} \langle \nabla R - 2\nabla_X X - 4\text{Ric}(X) - \frac{1}{\tau} X, Y \rangle d\tau + 2\sqrt{\tau} \langle X, Y \rangle_\tau \Big|_{\tau_1}^{\tau_2} +$$

$$\sum_{i=2}^{k-1} 2\sqrt{\tau_i'} \langle X^-(\tau_i') - X^+(\tau_i'), Y(\tau_i') \rangle_{\tau_i'}$$

*Definition7.2: Fixing a point $p$, we denote by $L(q, \bar{\tau})$ the $\mathcal{L}$ −length of the $\mathcal{L}$ −shortest curve $\gamma(\tau)$, $0 \leq \tau \leq \bar{\tau}$, joining $p$ and $q$. In other words*

$$L(q, \bar{\tau}) = \inf\{\mathcal{L}(\gamma) | \gamma: [0, \bar{\tau}] \to M \text{ with } \gamma(0) = p, \gamma(\bar{\tau}) = q\}$$

*Theorem7.3 (Perelman[39]) suppose that*

$$H(X) = -\frac{\partial R}{\partial \tau} - \frac{R}{\tau} - 2\langle \nabla R, X \rangle + 2Ric(X,X)$$

And

$$K = \int_0^{\bar{\tau}} \tau^{\frac{3}{2}} H(X) d\tau$$

Then

a) $|\nabla L|^2 = -4\bar{\tau}R + \frac{2}{\sqrt{\bar{\tau}}}L - \frac{4}{\sqrt{\bar{\tau}}}K$

b) $\frac{\partial L}{\partial \bar{\tau}} = 2\sqrt{\bar{\tau}}R - \frac{1}{2\bar{\tau}}L + \frac{1}{\bar{\tau}}K$

*Proof: The first variation formula in (7.5) implies that*

$$\nabla_Y L(q,\bar{\tau}) = \langle 2\sqrt{\bar{\tau}}X(\bar{\tau}), Y(\bar{\tau}) \rangle$$

*So* $\nabla L(q,\bar{\tau}) = 2\sqrt{\bar{\tau}}X(\bar{\tau})$

*At first we prove*

$$\frac{\partial L}{\partial \bar{\tau}}(\gamma(\bar{\tau}),\bar{\tau}) = \frac{d}{d\tau}L(\gamma(\tau),\tau)\Big|_{\tau=\bar{\tau}} - \langle \nabla L, X \rangle \quad (7.8)$$

*Since* $\gamma(\bar{\tau}) = -\frac{\partial}{\partial t} + X(\bar{\tau})$, *So the chain rule implies*

$$\frac{d}{d\tau}L(\gamma(\tau),\tau)\Big|_{\tau=\bar{\tau}} = \frac{d}{d\bar{\tau}}L(\gamma(\bar{\tau}),\bar{\tau}) + \langle \nabla L, X \rangle$$

*Therefore we get to (7.8)*

*Also from (7.1) and (7.8) we obtain*

$$\frac{\partial L}{\partial \bar{\tau}}(\gamma(\bar{\tau}),\bar{\tau}) = \sqrt{\bar{\tau}}(R + |X|^2) - \langle \nabla L, X \rangle$$

*And because* $\nabla L(q,\bar{\tau}) = 2\sqrt{\bar{\tau}}X(\bar{\tau})$ *so*

$$\frac{\partial L}{\partial \bar{\tau}}(\gamma(\bar{\tau}),\bar{\tau}) = \sqrt{\bar{\tau}}(R+|X|^2) - \sqrt{\bar{\tau}}|X|^2$$

$$= 2\sqrt{\bar{\tau}}R - \sqrt{\bar{\tau}}(R+|X|^2) \quad (7.9)$$

*Now we compute the $+|X|^2$.*

*According to (7.5) and this fact that $\gamma'(\tau) = \frac{\partial}{\partial \tau} + X(\tau)$, we get*

$$\frac{d}{d\tau}(R(\gamma(\tau),\tau)+|X(\tau)|^2) = \frac{\partial R}{\partial \tau} + \langle \nabla R, X \rangle + |X(\tau)|^2$$

*Also we know $|X(\tau)|^2 = \langle X(\tau), X(\tau) \rangle$*

*So from $\frac{d}{d\tau}|X(\tau)|^2 = \frac{d}{d\tau}\langle X(\tau), X(\tau) \rangle = 2\langle \nabla_X X, X \rangle + 2Ric(X,X)$*

*We get $\frac{d}{d\tau}(R(\gamma(\tau),\tau)+|X(\tau)|^2) = \frac{\partial R}{\partial t} + \langle \nabla R, X \rangle + 2\langle \nabla_X X, X \rangle + 2Ric(X,X)$*

*Also we have $\nabla_X X = \frac{1}{2}\nabla R - \frac{1}{2\tau}X - 2Ric(X,.)$ so*

$$\frac{d}{d\tau}(R(\gamma(\tau),\tau)+|X(\tau)|^2) = \frac{\partial R}{\partial t} + \frac{1}{\tau}R + 2\langle \nabla R, X \rangle - 2Ric(X,X) - \frac{1}{\tau}(R+|X|^2)$$

$$= -H(X) - \frac{1}{\tau}(R+|X|^2)$$

*So $\frac{d}{d\tau}\left(\tau^{\frac{3}{2}}(R+|X|^2)\right)\Big|_{\tau=\bar{\tau}} = \frac{1}{2}\sqrt{\bar{\tau}}(R+|X|^2) - \bar{\tau}^{\frac{3}{2}}H(X)$*

$$= \frac{1}{2}\frac{d}{d\tau}L(\gamma(\tau),\tau)\Big|_{\tau=\bar{\tau}} - \bar{\tau}^{\frac{3}{2}}H(X)$$

*So $\bar{\tau}^{\frac{3}{2}}(R+|X|^2) = \frac{1}{2}L(q,\bar{\tau}) - K$ where*

$$K = \int_0^{\bar{\tau}} \tau^{\frac{3}{2}} H(X) d\tau$$

*Therefore from This fact that $|\nabla L|^2 = 4\bar{\tau}|X|^2 = -4\bar{\tau}R + 4\bar{\tau}(R+|X|^2)$ and () we conclude a and b. So proof is complete.*

*Second variation formulae for ℒ −geodesics*

*Proposition 7.1 if we compute the second variation of energy when γ is a geodesic, we get*

$$\frac{d^2}{ds^2}E(\gamma) = \langle \nabla_Y Y, \gamma \rangle |_0^\tau + \int_0^\tau (\langle \nabla_Y \nabla_Y X, X \rangle + \langle \nabla_Y X, \nabla_Y X \rangle) dt$$

*Now turn to the second spatial variation of the reduced length.*

*Theorem 7.4(Perelman[39]) for any ℒ −geodesic γ , we have*

$\delta_Y^2(\mathcal{L}) =$

$2\sqrt{\tau}\langle \nabla_Y Y, X \rangle |_0^{\bar{\tau}} + \int_0^{\bar{\tau}} \sqrt{\tau}[2|\nabla_X Y|^2 + 2\langle R(Y,X)Y, X \rangle + \nabla_Y \nabla_Y R + 2\nabla_X Ric(Y,Y) - 4\nabla_Y Ric(Y,X)] d\tau$.

*Proof to compute the second variation $\delta_Y^2(\mathcal{L})$. we start the first variation formula .*

*We know*

$$\delta_Y(\mathcal{L}) = \int_0^{\bar{\tau}} \sqrt{\tau}(\nabla_Y R + 2\langle \nabla_X Y, X \rangle) d\tau$$

So $\qquad \delta_Y^2(\mathcal{L}) = Y\left(\int_0^{\bar{\tau}} \sqrt{\tau}(\nabla_Y R + 2\langle \nabla_X Y, X \rangle) d\tau\right)$

*We know $\nabla_Y R = YR$ ,so $\delta_Y^2(\mathcal{L}) = Y\left(\int_0^{\bar{\tau}} \sqrt{\tau}(YR + 2\langle \nabla_X Y, X \rangle) d\tau\right)$*

$$= \left(\int_0^{\bar{\tau}} \sqrt{\tau}\bigl(Y(Y(R)) + 2\langle \nabla_Y \nabla_X Y, X \rangle + 2|\nabla_X Y|^2\bigr) d\tau\right)$$

*Because $\nabla_X Y = \nabla_Y X$ we get*

$$= \left(\int_0^{\bar{\tau}} \sqrt{\tau}\bigl(Y(Y(R)) + 2\langle \nabla_Y \nabla_X Y, X \rangle + 2|\nabla_X Y|^2\bigr) d\tau\right)$$

*Also we know $\nabla_X \nabla_Y Z - \nabla_Y \nabla_X Z - \nabla_{[X,Y]} Z = R(X,Y)Z$ for all vectors $X, Y, Z$.so*

$$2\langle \nabla_Y \nabla_X Y, X \rangle = 2\langle \nabla_X \nabla_Y Y, X \rangle + 2\langle R(Y,X)Y, X \rangle$$

*For continuing proof, we prove the following equality*

$$\frac{\partial}{\partial \tau}\langle \nabla_Y Y, X\rangle = \langle \nabla_X \nabla_Y Y, X\rangle + \langle \nabla_Y Y, \nabla_X X\rangle + 2Ric(\nabla_Y Y, X) + \langle \frac{\partial}{\partial \tau}\nabla_Y Y, X\rangle \quad (7.10)$$

*We can break $\frac{\partial}{\partial \tau}\langle \nabla_Y Y, X\rangle$ into two parts: the first assumes that the metric is constant and the second deals with the variation with $\tau$ of the metric. The first contribution is the usual formula*

$$\frac{\partial}{\partial \tau}\langle \nabla_Y Y, X\rangle = \langle \nabla_X \nabla_Y Y, X\rangle + \langle \nabla_Y Y, \nabla_X X\rangle$$

*This gives the first two terms of the right-hand side of the equation. we show that the last two terms in that equation come from differentiating the metric with respect to $\tau$. To do this recall that in local coordinates, writing the metric as $g_{ij}$, we have*

$$\langle \nabla_Y Y, X\rangle = g_{ij}\left(Y^k \partial_k Y^i + \Gamma^i_{kl} Y^k Y^l\right) X^j$$

*There are two contributions coming from differentiating the metric with respect to $\tau$. The first is when we differentiate $g_{ij}$. This leads to*

$$2Ric_{ij}\left(Y^k \partial_k Y^i + \Gamma^i_{kl} Y^k Y^l\right) X^j = 2Ric(\nabla_Y Y, X)$$

*The other contribution is from differentiating the Christoffel symbols. This yields*

$$g_{ij}\frac{\partial \Gamma^i_{kl}}{\partial \tau} Y^k Y^l X^j$$

*Differentiating the formula $\Gamma^i_{kl} = \frac{1}{2}g^{si}(\partial_k g_{sl} + \partial_l g_{sk} - \partial_s g_{kl})$ leads to*

$$g_{ij}\frac{\partial \Gamma^i_{kl}}{\partial \tau} = -2Ric_{ij}\Gamma^i_{kl} + g_{ij}g^{si}(\partial_k Ric_{sl} + \partial_l Ric_{sk} - \partial_s Ric_{kl})$$

$$= -2Ric_{ij}\Gamma^i_{kl} + \partial_k Ric_{jl} + \partial_l Ric_{jk} - \partial_j Ric_{kl}$$

*Thus, we have*

$$g_{ij}\frac{\partial \Gamma^i_{kl}}{\partial \tau}Y^k Y^l X^j = \left(-2Ric_{ij}\Gamma^i_{kl} + \partial_k Ric_{jl}\right)Y^k Y^l X^j$$

$$= \nabla_X Ric(Y,Y)$$

*So we obtain*

$$2\langle \nabla_Y \nabla_X Y, X\rangle = 2\frac{d}{d\tau}\langle \nabla_Y Y, X\rangle - 4Ric(\nabla_Y Y, X) - 2\langle \nabla_Y Y, \nabla_X X\rangle - 2\langle \frac{\partial}{\partial \tau}\nabla_Y Y, X\rangle$$
$$+ 2\langle \nabla_X \nabla_Y Y, X\rangle + 2\langle R(Y,X)Y, X\rangle$$

*Also we can compute*

$$\langle \frac{\partial}{\partial \tau}\nabla_Y Y, X\rangle = 2(\nabla_Y Ric)(Y,X) - (\nabla_X Ric)(Y,Y) \quad (7.11)$$

*Hence from (7.10) and (7.11) we get*

$$\frac{\partial}{\partial \tau}\langle \nabla_Y Y, X\rangle = \langle \nabla_X \nabla_Y Y, X\rangle + \langle \nabla_Y Y, \nabla_X X\rangle + 2Ric(\nabla_Y Y, X) + 2(\nabla_Y Ric)(Y,X) -$$
$$(\nabla_X Ric)(Y,Y) \hspace{5cm} (7.12)$$

*Suppose $Y(0) = 0$ and the fact that $\sqrt{\tau}X(\tau)$ has a limit as $\tau \to 0$ are used to get the third equality below. So applying (7.12) and integrating by parts, we arrive*

$$\delta_Y^2(\mathcal{L}(\gamma)) = \left(\int_0^{\bar{\tau}}\sqrt{\tau}(Y(Y(R)) + 2\langle R(Y,X)Y,X\rangle + 2|\nabla_Y X|^2)d\tau\right)$$
$$+ 2\int_0^{\bar{\tau}}\sqrt{\tau}\left(\frac{\partial}{\partial \tau}\langle \nabla_Y Y, X\rangle - \langle \nabla_Y Y, \nabla_X X\rangle - 2Ric(\nabla_Y Y, X)\right.$$
$$\left. - 2(\nabla_Y Ric)(Y,X) + (\nabla_X Ric)(Y,Y)\right)d\tau$$

$$= \left(\int_0^{\bar{\tau}}\sqrt{\tau}(Y(Y(R)) + 2\langle R(Y,X)Y,X\rangle + 2|\nabla_Y X|^2)d\tau\right) +$$
$$2\int_0^{\bar{\tau}}\sqrt{\tau}\big(-\langle \nabla_Y Y, \nabla_X X\rangle - 2Ric(\nabla_Y Y, X) -$$
$$2(\nabla_Y Ric)(Y,X) + (\nabla_X Ric)(Y,Y)\big)d\tau + 2\sqrt{\tau}\langle \nabla_Y Y, X\rangle\big|_0^{\bar{\tau}} -$$
$$\int_0^{\bar{\tau}}\frac{1}{\sqrt{\tau}}\langle \nabla_Y Y, X\rangle d\tau = 2\sqrt{\bar{\tau}}\langle \nabla_Y Y, X\rangle + \left(\int_0^{\bar{\tau}}\sqrt{\tau}(Y(Y(R)) - \nabla_Y Y.\nabla R + 2\langle R(Y,X)Y, X\rangle +\right.$$

$2|\nabla_Y X|^2)) \, d\tau +$

$2 \int_0^{\bar\tau} \sqrt{\tau} \left( -\langle \nabla_Y Y, [\nabla_X X + 2Ric(X) - \frac{1}{2}\nabla R + \frac{1}{2\tau}X] \rangle - 2(\nabla_Y Ric)(Y,X) + (\nabla_X Ric)(Y,Y) \right) d\tau$

$$= 2\sqrt{\bar\tau}\langle \nabla_Y Y, X \rangle + \left( \int_0^{\bar\tau} \sqrt{\tau}(\nabla_{Y,Y}^2 R + 2\langle R(Y,X)Y, X\rangle + 2|\nabla_Y X|^2) d\tau \right)$$

$$+ \int_0^{\bar\tau} \sqrt{\tau}\bigl(-4(\nabla_Y Ric)(Y,X) + 2(\nabla_X Ric)(Y,Y)\bigr) d\tau$$

Because we know $\nabla_X X - \frac{1}{2}\nabla R + \frac{1}{2\tau}X + 2Ric(X,.) = 0$ and $Hess(f)(X,Y) = \nabla_X \nabla_Y (f) - \nabla_{\nabla_X Y}(f)$, so proof is complete.

*ℒ and Riemannian distance:*

*Theorem 7.5 (see[43]) let $\gamma: [0,\bar\tau] \to M$, $\bar\tau \in (0,T]$, be a $C^1$-path starting at $p$ and ending at $q$.*

i)     *(bounding Riemannian distance by ℒ) for any $\tau \in [0,\bar\tau]$ we have*

$$d_{g(0)}^2(p,\sigma(\tau)) \leq 2\sqrt{\tau}e^{2C_0\tau}\left(\mathcal{L}(\gamma) + \frac{2nC_0}{3}\tau^{\frac{3}{2}}\right)$$

*Where $(M^n, g(\tau)), \tau \in [0,T]$, denote a complete solution to the backward Ricci flow, and $p \in M$ shall be a base point. also we assume the curvature bound*

$$\max_{(x,t)\in M\times[0,T]}\{|Rm(x,\tau)|, |Rc(x,\tau)|\} \leq C_0 < \infty$$

*Remark for i) when M is noncompact. from i) we conclude for any $\bar\tau \in (0,T]$.*

$$\lim_{q\to\infty} \bar{L}(q,\bar\tau) := \lim_{q\to\infty} 2\sqrt{\bar\tau}L(q,\bar\tau) = \infty$$

ii)     *(bounding speed at some time by ℒ) there exists $t_* \in (0,\bar\tau)$ such that*

$$\tau_* \left|\frac{\partial \gamma}{\partial \tau}(\tau_*)\right|^2_{g(\tau_*)} = \left|\frac{\partial \beta}{\partial \sigma}(\sigma_*)\right|^2_{g(\tau_*)} \leq \frac{1}{2\sqrt{\bar\tau}} \mathcal{L}(\gamma) + \frac{nC_0}{3}\bar\tau$$

Where $(\sigma) := \gamma(\tau)$, $\sigma = 2\sqrt{\tau}$, and $\sigma_* = 2\sqrt{\tau_*}$.

iii)  (bounding L by Riemannian distance) for any $q \in M$ and $\bar\tau > 0$

$$L(q,\bar\tau) \leq e^{2C_0\bar\tau} \frac{d^2_{g(\bar\tau)}(p,q)}{2\sqrt{\bar\tau}} + \frac{2nC_0}{3}\bar\tau^{\frac{3}{2}}$$

For proof of this theorem, we start with a remark

Remark 7.5: for $\tau_1 < \tau_2$ and $x \in M$

$$e^{-2C_0(\tau_2-\tau_1)} g(\tau_2, x) \leq g(\tau_1, x) \leq e^{2C_0(\tau_2-\tau_1)} g(\tau_2, x)$$

Proof i) suppose $\bar\sigma = 2\sqrt{\bar\tau}$ and $\beta(\bar\sigma) := \gamma(\bar\tau)$. At first we compute

$$\int_0^{2\sqrt{\bar\tau}} \left|\frac{\partial \beta}{\partial \bar\sigma}(\bar\sigma)\right|^2 d\bar\sigma = \mathcal{L}(\gamma) - \int_{2\sqrt{\tau}}^{\sqrt{\bar\tau}} \left|\frac{\partial \beta}{\partial \bar\sigma}(\bar\sigma)\right|^2_{g\left(\frac{\bar\sigma^2}{4}\right)} d\bar\sigma - \int_0^{\bar\tau} \sqrt{\bar\tau} R(\gamma(\bar\tau), \bar\tau) d\bar\tau$$

$$\leq \mathcal{L}(\gamma) + \frac{2nC_0}{3}\bar\tau^{\frac{3}{2}}$$

Because $R \geq -nC_0$. Hence, since $g(0) \leq e^{2C_0\bar\tau} g(\bar\tau)$ for $\bar\tau \in [0,\tau]$ we get

$$d^2_{g(0)}(p, \gamma(\tau)) \leq e^{2C_0\bar\tau} \left(\int_0^{2\sqrt{\bar\tau}} \left|\frac{\partial \beta}{\partial \bar\sigma}(\bar\sigma)\right|_{g\left(\frac{\bar\sigma^2}{4}\right)} d\bar\sigma\right)^2$$

$$\leq e^{2C_0\bar\tau} 2\sqrt{\bar\tau} \int_0^{2\sqrt{\bar\tau}} \left|\frac{\partial \beta}{\partial \bar\sigma}(\bar\sigma)\right|^2_{g\left(\frac{\bar\sigma^2}{4}\right)} d\bar\sigma \leq 2\sqrt{\bar\tau} e^{2C_0\bar\tau} \left(\mathcal{L}(\gamma) + \frac{2nC_0}{3}\bar\tau^{\frac{3}{2}}\right)$$

ii) by taking $\tau = \bar\tau$ in proof of i) we get

$$\frac{1}{2\sqrt{\bar\tau}} \int_0^{2\sqrt{\bar\tau}} \left|\frac{\partial \beta}{\partial \bar\sigma}(\bar\sigma)\right|^2_{g\left(\frac{\bar\sigma^2}{4}\right)} d\bar\sigma \leq \frac{1}{2\sqrt{\bar\tau}} \mathcal{L}(\gamma) + \frac{nC_0}{3}\bar\tau$$

Now with using the mean value Theorem for integrals, There exists $t_* \in (0, \bar\tau)$ such that

$$\left|\frac{\partial \beta}{\partial \sigma}(\sigma_*)\right|^2_{g(\tau_*)} \leq \frac{1}{2\sqrt{\bar{\tau}}} \mathcal{L}(\gamma) + \frac{nC_0}{3}\bar{\tau}$$

iii) *Let* $\eta: [0, 2\sqrt{\bar{\tau}}] \to M$ *be a minimal geodesic from* $p$ *to* $q$ *with respect to metric* $g(\bar{\tau})$. *Then because on* $g(\bar{\tau})$ *we have*

$$2\sqrt{\bar{\tau}} L(\gamma(\tau), \tau) = 2\sqrt{\bar{\tau}} \mathcal{L}(\gamma|_{[0,\tau]}) = \frac{\tau}{\bar{\tau}} d^2_{g(\bar{\tau})}(p, q)$$

*So*

$$L(q, \bar{\tau}) \leq \mathcal{L}(\eta) = \int_0^{2\sqrt{\bar{\tau}}} \left(\frac{\sigma^2}{4} R\left(\eta(\sigma), \frac{\sigma^2}{4}\right) + \left|\frac{d\eta}{d\sigma}\right|^2_{g(\tau)}\right) d\sigma$$

$$\leq \int_0^{2\sqrt{\bar{\tau}}} \left(\frac{nC_0 \sigma^2}{4} + e^{2C_0 \bar{\tau}} \left|\frac{d\eta}{d\sigma}\right|^2_{g(\tau)}\right) d\sigma \leq \frac{2nC_0}{3} \bar{\tau}^{\frac{3}{2}} + \frac{e^{2C_0 \bar{\tau}}}{2\sqrt{\bar{\tau}}} d^2_{g(\bar{\tau})}(p, q)$$

*So proof is complete.*

*the* $\mathcal{L}$-*Jacobi equation :consider a family* $\gamma(\tau, u)$ *of* $\mathcal{L}$-*geodesics parameterized by* $s$ *and defined on* $[\tau_1, \tau_2]$ *with* $0 \leq \tau_1 \leq \tau_2$ .*Let* $Y(\tau)$ *be a vector field along* $\gamma$ *defined by*

$$Y(\tau) = \frac{\partial}{\partial s} \gamma(\tau, s)\bigg|_{s=0}$$

*Now from second variation Theorem we obtain* $Y(\tau)$ *satisfies the* $\mathcal{L} -$*Jacobi equation*

$$\nabla_X \nabla_X Y + R(Y, X)X - \frac{1}{2}\nabla_Y(\nabla R) + \frac{1}{2\tau}\nabla_X Y + 2(\nabla_Y Ric)(X, .) + 2Ric(\nabla_X Y, .) = 0 \quad (7.13)$$

*This is a second –order linear equation for Y .supposing that* $\tau_1 > 0$ *,there is a unique vector field* Y *along* $\gamma$ *solving this equation ,vanishing at* $\tau_1$ *with a given first-order derivative along* $\gamma$ *at* $\tau_1$.*similarly ,there is a unique solution* Y *to this equation ,vanishing at* $\tau_2$ *and with a given first order derivative at* $\tau_2$ .

*Definition 7.3 a field* $Y(\tau)$ *along an* $\mathcal{L} -$*geodesic is called an* $\mathcal{L} -$*jacobi field if it satisfies the* $\mathcal{L} -$*jacobi equation ,equation (7.13) ,and if it vanishes at* $\tau_1$.

*Notation: for every vector field Y along γ we denote by $Jac(Y)$ the expression on the left-hand side of equation (7.13).*

*In one of remarks we considered that the vector $\lim_{\tau \to 0} \sqrt{\tau} X(\tau)$ exists. Now we give a similar result even for $\tau_1 = 0$.*

*Lemma 7.5 let γ be an $\mathcal{L}$ − geodesic defined on $[0, \tau_2]$ and let $Y(\tau)$ be an $\mathcal{L}$ −Jacobi field along γ. Then $\lim_{\tau \to 0} \sqrt{\tau} \nabla_X Y$ exists, furthermore, $Y(\tau)$ is completely determined by this limit.*

*Remark 7.6 we can consider that the bilinear pairing*

$$-\int_{\tau_1}^{\tau_2} 2\sqrt{\tau} \langle Jac(Y_1), Y_2 \rangle d\tau$$

*Is a symmetric function of $Y_1$ and $Y_2$. (here we assume that $Y_1(\tau_2) = Y_2(\tau_2) = 0$ and γ is an $\mathcal{L}$ − geodesic and $Y_1, Y_2$ are vector fields along γ vanishing at $\tau_1$).*

**Estimate the Hessian of the $\mathcal{L}$ −distance functions:**

*Here we give an inequality for the hessian of $\mathcal{L}$ involving the integral of the vector field along.*

*Let $\gamma: [0, \bar{\tau}] \to M$ be an $\mathcal{L}$ −shortest curve connecting p and q. we fix a vector Y at $\tau = \bar{\tau}$ with $|Y|_{g_{ij}(\bar{\tau})} = 1$, and extend Y along the $\mathcal{L}$ −shortest geodesic γ on $[0, \bar{\tau}]$ by solving the following ODE*

$$\nabla_X Y = -Ric(Y,.) + \frac{1}{2\tau} Y. \quad (7.14)$$

*Lemma 7.6 (Perelman) suppose $\{Y_1, Y_2, \dots Y_n\}$ is an orthonormal basis at $\tau = \bar{\tau}$ with respect to metric $g_{ij}(\bar{\tau})$ and solve for $Y(\tau)$ in the equation (7.14), then $\{Y_1(\tau), Y_2(\tau), \dots Y_n(\tau)\}$ remains orthogonal on $[0, \bar{\tau}]$ and*

$$\langle Y_i(\tau), Y_j(\tau) \rangle = \frac{\tau}{\bar{\tau}} \delta_{ij}$$

*Proof according to (7.6) and (7.14) we get*

$$\frac{d}{d\tau}\langle Y_i, Y_j\rangle = 2Ric(Y_i, Y_j) + \langle \nabla_X Y_i, Y_j\rangle + \langle Y_i, \nabla_X Y_j\rangle$$

$$\langle Y_i(\tau), Y_j(\tau)\rangle = \frac{\tau}{\bar{\tau}}\delta_{ij} \Rightarrow |Y(\tau)|^2 = \frac{\tau}{\bar{\tau}}$$

So $\{Y_1(\tau), Y_2(\tau), \ldots Y_n(\tau)\}$ remains orthogonal on $[0, \bar{\tau}]$ with $Y_i(0) = 0$, $i = 0, 1, \ldots, n$. So proof is complete.

The main result of this sub section is following Theorem from Perelman.

Theorem7.6 (Perelman[40]) suppose that $|Y|_{g_{ij}(\bar{\tau})} = 1$ at any point $q \in M$, consider an $\mathcal{L}$-shortest geodesic $\gamma$ connecting $p$ to $q$ and extend $Y$ along $\gamma$ by solving (7.14). So the Hessian of the $\mathcal{L}$-distance function L on M with $\tau = \bar{\tau}$ satisfies

$$Hess_L(Y, Y) \leq \frac{1}{\sqrt{\bar{\tau}}} - 2\sqrt{\bar{\tau}}Ric(Y, Y) - \int_0^{\bar{\tau}} \sqrt{\tau}Q(X, Y)d\tau \quad (7.15)$$

Where

$Q(X, Y) = -\nabla_Y \nabla_Y R - 2\langle R(Y, X)Y, X\rangle - 4\nabla_X Ric(Y, Y) + 4\nabla_Y Ric(Y, X) - 2Ric_\tau(Y, Y) + 2|Ric(Y, \cdot)|^2 - \frac{1}{\tau}Ric(Y, Y)$

is the Li-Yau-Hamilton quadratic.

Proof according to the definition $Hess_L(Y, Y)$ we have

$$Hess_L(Y, Y) = Y(Y(L))(\bar{\tau}) - \langle \nabla_Y Y, \nabla L\rangle(\bar{\tau}) \quad (7.16)$$

also $Y(Y(L))(\bar{\tau}) = \delta_Y^2(L) \leq \delta_Y^2(\mathcal{L})$. and $\nabla L(q, \bar{\tau}) = 2\sqrt{\bar{\tau}}X$. So $\langle \nabla_Y Y, \nabla L\rangle = 2\sqrt{\bar{\tau}}\langle \nabla_Y Y, X\rangle$. Therefore from (7.16)

$$Hess_L(Y, Y) \leq \delta_Y^2(\mathcal{L}) - 2\sqrt{\bar{\tau}}\langle \nabla_Y Y, X\rangle(\bar{\tau}) \quad (7.17)$$

So from second variation formula and (7.17) we get

$Hess_L(Y, Y)$
$$\leq \int_0^{\bar{\tau}} \sqrt{\tau}[2|\nabla_X Y|^2 + 2\langle R(Y, X)Y, X\rangle + \nabla_Y \nabla_Y R + 2\nabla_X Ric(Y, Y) - 4\nabla_Y Ric(Y, X)] d\tau$$

$$= \int_0^{\bar{\tau}} \sqrt{\tau} \left[ 2\left| -\text{Ric}(Y,.) + \frac{1}{2\tau}Y \right|^2 + 2\langle R(Y,X)Y,X \rangle + \nabla_Y \nabla_Y R + 2\nabla_X Ric(Y,Y) \right.$$
$$\left. - 4\nabla_Y Ric(Y,X) \right] d\tau$$

*So because $|Y|^2 = 1$ we get*

$$\left| -\text{Ric}(Y,.) + \frac{1}{2\tau}Y \right|^2 = |Ric(Y,.)|^2 - \frac{1}{\tau}Ric(Y,Y) + \frac{1}{4\tau\bar{\tau}}$$

*So*

$Hess_L(Y,Y)$

$$\leq \int_0^{\bar{\tau}} \sqrt{\tau} \left[ 2|Ric(Y,.)|^2 - \frac{2}{\tau}Ric(Y,Y) + \frac{1}{2\tau\bar{\tau}} + 2\langle R(Y,X)Y,X \rangle + \nabla_Y \nabla_Y R \right.$$
$$\left. + 2\nabla_X Ric(Y,Y) - 4\nabla_Y Ric(Y,X) \right] d\tau$$

*Also because*

$$\frac{d}{d\tau} Ric(Y,Y) = \frac{\partial}{\partial \tau} Ric(Y,Y) + \nabla_X Ric(Y,Y) + 2Ric(\nabla_X Y, Y)$$

*And from (7.14)*

$$\frac{d}{d\tau} Ric(Y,Y) = \frac{\partial}{\partial \tau} Ric(Y,Y) + \nabla_X Ric(Y,Y) - 2|Ric(Y,.)|^2 + \frac{1}{\tau}Ric(Y,Y)$$

*Therefore we have*

$Hess_L(Y,Y)$

$$\leq \int_0^{\bar{\tau}} \sqrt{\tau} \left[ 2|Ric(Y,.)|^2 - \frac{2}{\tau}Ric(Y,Y) + \frac{1}{2\tau\bar{\tau}} + 2\langle R(Y,X)Y,X \rangle + \nabla_Y \nabla_Y R \right.$$
$$- 4(\nabla_Y Ric)(Y,Y)$$
$$- \left( 2\frac{d}{d\tau}Ric(Y,Y) - 2\frac{\partial}{\partial\tau}Ric(Y,Y) + 4|Ric(Y,.)|^2 - \frac{2}{\tau}Ric(Y,Y) \right)$$
$$\left. + 4\nabla_X Ric(Y,Y) \right] d\tau$$

$$= -\int_0^{\bar{\tau}} \left[ 2\sqrt{\tau}\frac{d}{d\tau}Ric(Y,Y) + \frac{1}{\sqrt{\tau}}Ric(Y,Y) \right] d\tau + \frac{1}{2\bar{\tau}}\int_0^{\bar{\tau}} \frac{1}{\sqrt{\tau}} d\tau +$$

$$+ \int_0^{\bar{\tau}} \sqrt{\tau} \left[ 2\langle R(Y,X)Y,X \rangle + \nabla_Y \nabla_Y R + \frac{1}{\tau} Ric(Y,Y) + 4\nabla_X Ric(Y,Y) - \nabla_Y Ric(X,Y) \right.$$

$$\left. + \frac{2\partial Ric}{\partial \tau}(Y,Y) - 2|Ric(Y,.)|^2 \right] d\tau$$

$$= \frac{1}{\sqrt{\bar{\tau}}} - 2\sqrt{\bar{\tau}} Ric(Y,Y) - \int_0^{\bar{\tau}} \sqrt{\tau} Q(X,Y) d\tau$$

*So proof is complete.*

*Theorem7.7 (Perelman[43]) suppose that $K = \int_0^{\bar{\tau}} \tau^{\frac{3}{2}} H(X) d\tau$ then we have*

$$\Delta L \leq \frac{n}{\sqrt{\bar{\tau}}} - 2\sqrt{\bar{\tau}} R - \frac{1}{\bar{\tau}} K$$

*Proof  let $\{Y_1, Y_2, ... Y_n\}$ be an orthonormal basis at $\tau = \bar{\tau}$ and with extend them along the shortest $\mathcal{L}-$geodesic $\gamma$ and taking $Y = Y_i$ in (7.15) and summing over $i$, we obtain*

$$\sum Hess_L(Y_i, Y_i) \leq \frac{n}{\sqrt{\bar{\tau}}} - 2\sqrt{\bar{\tau}} R - \sum_{i=1}^n \int_0^{\bar{\tau}} \sqrt{\tau} Q(X,Y_i) \, d\tau$$

*But we know*        $\Delta L = \sum_i Hess(L)(Y_i, Y_i)$

*So*        $\Delta L \leq \frac{n}{\sqrt{\bar{\tau}}} - 2\sqrt{\bar{\tau}} R - \sum_{i=1}^n \int_0^{\bar{\tau}} \sqrt{\tau} Q(X,Y_i) \, d\tau$   *(7.18)*

*Now we prove $\sum_i Q(X,Y_i) = \frac{\tau}{\bar{\tau}} Q(X)$. We know $\langle Y_i(\tau), Y_i(\tau) \rangle = \frac{\tau}{\bar{\tau}} \delta_{ij}$ so $\left\{ \sqrt{\frac{\bar{\tau}}{\tau}} Y_i(\tau) \right\}$ is an orthonormal basis at $\tau$ .also*

$$\sum_i Q(X,Y_i) = -\frac{\tau}{\bar{\tau}} \Delta R + 2\frac{\tau}{\bar{\tau}} Ric(X,X) - 4\frac{\tau}{\bar{\tau}} \langle \nabla R, X \rangle + 4\frac{\tau}{\bar{\tau}} \sum_i \nabla_{Y_i} Ric(Y_i, X)$$

$$- 2\frac{\tau}{\bar{\tau}} \sum_i Ric_\tau(Y_i, Y_i) + 2\frac{\tau}{\bar{\tau}} \sum_i |Ric(Y_i,.)|^2 - \frac{1}{\bar{\tau}} \sum_i Ric(Y_i, Y_i)$$

*Tracing the second Bianchi identity gives*

$$\sum_i \nabla_{Y_i} Ric(Y_i, X) = \frac{1}{2}\langle \nabla R, X \rangle$$

*Also according to following identity*

$$Ric_\tau(Y,Y) = (\nabla R_{ij})Y^i Y^j + 2R_{ikjl}R_{kl}Y^i Y^j - 2R_{ik}R_{jk}Y^i Y^j \text{ we get } \sum_i Ric_\tau(Y_i, Y_i) = \Delta R$$

*Putting this together gives*

$$\sum_i Q(X, Y_i)(\tau) = \frac{\tau}{\bar{\tau}} Q(X) \quad (7.19)$$

*Where* $Q(X) = -\frac{\partial R}{\partial \tau} - \frac{1}{\tau} R - 2\langle \nabla R, X \rangle + 2Ric(X, X)$

*So from (7.18),(7.19) we get*

$$\Delta L \leq \frac{n}{\sqrt{\bar{\tau}}} - 2\sqrt{\bar{\tau}} R - \int_0^{\bar{\tau}} \sqrt{\tau} \left(\frac{\tau}{\bar{\tau}}\right) Q(X) d\tau$$

$$= \frac{n}{\sqrt{\bar{\tau}}} - 2\sqrt{\bar{\tau}} R - \frac{1}{\bar{\tau}} K.$$

*Corollary 7.2 we have*

$$Hess_L(Y,Y) = \frac{1}{\sqrt{\bar{\tau}}} - 2\sqrt{\bar{\tau}} Ric(Y,Y) - \int_0^{\bar{\tau}} \sqrt{\tau}\, Q(X,Y) d\tau$$

*If and only if* $Y(\tau), \tau \in [0, \bar{\tau}]$, *is an* $\mathcal{L}$ −*Jacobian field*

*Proposition 7.2(see[43])we have*

$$\Delta L = \frac{n}{\sqrt{\bar{\tau}}} - 2\sqrt{\bar{\tau}} R - \frac{1}{\bar{\tau}} K$$

*if and only if we are on a gradient shrinking soliton with*

$$R_{ij} + \frac{1}{2\sqrt{\bar{\tau}}} \nabla_i \nabla_j L = \frac{1}{2\bar{\tau}} g_{ij}$$

*Proof when* $Y_i(\tau)\ i = 1, \dots, n$ *are* $\mathcal{L}$ −*Jacobian fields along* $\gamma$ *,we have*

$$\frac{d}{d\tau}\langle Y_i(\tau), Y_j(\tau)\rangle = 2Ric(Y_i, Y_j) + \langle \nabla_X Y_i, Y_j\rangle + \langle Y_i, \nabla_X Y_j\rangle$$

But $\nabla L = 2\sqrt{\tau}X$ so $X = \frac{\nabla L}{2\sqrt{\tau}}$ therefore

$$\frac{d}{d\tau}\langle Y_i(\tau), Y_j(\tau)\rangle = 2Ric(Y_i, Y_j) + \langle \nabla_{Y_i}\left(\frac{\nabla L}{2\sqrt{\tau}}\right), Y_j\rangle + \langle Y_i, \nabla_{Y_j}\left(\frac{\nabla L}{2\sqrt{\tau}}\right)\rangle$$

$$= 2Ric(Y_i, Y_j) + \frac{1}{\sqrt{\tau}}Hess_L(Y_i, Y_j)$$

But we know $\langle Y_i(\tau), Y_j(\tau)\rangle = \frac{\tau}{\bar{\tau}}\delta_{ij}$ so from last equality we get

$$R_{ij} + \frac{1}{2\sqrt{\tau}}\nabla_i\nabla_j L = \frac{1}{2\bar{\tau}}g_{ij}$$

*Also the vice versa of proof is obvious .so proof is complete.*

**Li –Yau –Perelman distance**

*We introduce the Li-Yau-Perelman distance both on the tangent space and on space-time .The reason that the Li-Yau-Perelman distance $l = l(q,\bar{\tau})$ is easier to work with is that it is scale invariant when $\tau_1 = 0$ .*

*Definition7.4 the Li-Yau-Perelman distance $l = l(q,\bar{\tau})$ is defined by*

$$l = l(q,\bar{\tau}) = \frac{L(q,\bar{\tau})}{2\sqrt{\bar{\tau}}}$$

*So in summary if we write the Perelman works on l that we proved in previous theorem, we get the following theorem.*

*Theorem7.8 (Perelman[43]) For the Li-Yau-Perelman distance $l(q,\bar{\tau})$ we have*

a) $\frac{\partial l}{\partial \bar{\tau}} = -\frac{l}{\bar{\tau}} + R + \frac{1}{2\bar{\tau}^{\frac{3}{2}}}K$

b) $|\nabla l|^2 = -R + \frac{l}{\bar{\tau}} - \frac{1}{2\bar{\tau}^{\frac{3}{2}}}K$

c) $\Delta l \leq -R + \frac{n}{2\bar{\tau}} - \frac{1}{2\bar{\tau}^{\frac{3}{2}}}K$

*Now we get the upper bound on the minimum of $l(.,\tau)$ for every $\tau$*

*Lemma 7.7(see[38]) $\min_{x\in M} l(.,\tau) \leq \frac{n}{2}$ for every $\tau$ .*

*Proof we know*

1). $\frac{\partial L}{\partial \tau} = 2\sqrt{\tau}R - \frac{1}{2\tau}L + \frac{1}{\tau}K$

2). $\Delta L \leq \frac{n}{\sqrt{\tau}} - 2\sqrt{\tau}R - \frac{1}{\tau}K$

3). $|\nabla L|^2 = -4\tau R + \frac{2}{\sqrt{\tau}}L - \frac{4}{\sqrt{\tau}}K$

*So 1,2 and 3 gives us*

$$\frac{\partial l}{\partial \tau} - \Delta l + |\nabla l|^2 - R + \frac{n}{2\tau} \leq 0$$

*Also 2 and 3 gives us*

$$2\Delta l - |\nabla l|^2 + R + \frac{l-n}{\tau} \leq 0$$

*Let $\bar{L} = 2\sqrt{\tau}L$ .Therefore 1 and 2 gives us*

$$\frac{\partial \bar{L}}{\partial \tau} + \Delta \bar{L} \leq 2n$$

*So* $\quad \frac{\partial(\bar{L}-2n\tau)}{\partial \tau} + \nabla(\bar{L} - 2n\tau) \leq 0$ .

*Thus, by a standard maximum principle argument, $\min\{\bar{L}(q,\tau) - 2n\tau | q \in M\}$ is non-increasing and therefore $\min\{\bar{L}(q,\tau)|q \in M\} \leq 2n\tau$, so $\min l(.,\tau) \leq \frac{n}{2}$.*

**Estimates on the Li-Yau-Perelman distance**

*Proposition7.3 if the metrics $g_{ij}(\tau)$ have non-negative curvature operator and if the flow exists for $\tau \in [o,\tau_0]$ ,Then*

$$|\nabla l|^2 + R \leq \frac{cl}{\tau}$$

*For some constant c,when ever $\tau$ is bounded away from $\tau_0$ ,say $\tau \leq (1-c)\tau_0$ ,where $c > 0$ .*

*Lemma7.8 (Perelman) if we have a Ricci flow $\frac{\partial g_{ij}}{\partial \tau} = 2R_{ij}$ , then $R(.,\tau) \geq -\frac{n}{2(\bar{\tau}-\tau)}$ whenever the flow exists for $\tau \in [0,\bar{\tau}]$ .*

*Proof we know $\frac{\partial R}{\partial \tau} = -\Delta R - 2|Ric|^2$.Look at the corresponding ODE, $\frac{\partial R}{\partial \tau} = -2|Ric|^2$ .Since $R = tr(Ric)$, $|Ric|^2 \geq \frac{1}{n}R^2$ and therefore ,$-2|Ric|^2 \leq \frac{-2|R|^2}{n}$, i.e .$\frac{\partial R}{\partial \tau} \leq \frac{-2|R|^2}{n}$ .By solving this equation we get that the set $R(.,\tau) \geq -\frac{n}{2(\tau_0-\tau)}$ is preserved by the ODE and therefore it is preserved by the corresponding PDE.*

*Corollary7.3 For every $\bar{\tau} > 0, \exists q \in M$ and $\exists \mathcal{L}$ −geodesic $\gamma: [0,\bar{\tau}] \to M$, $\gamma(0) = p$ and $\gamma(\bar{\tau}) = q$ such that $\mathcal{L}(\gamma) \leq n\sqrt{\bar{\tau}}$ .*

*We proved $\min_{\bar{\tau}>0} l(.,\bar{\tau}) \leq \frac{n}{2}$ .An analogy with this ideas can be found in the original proof of the Harnack inequality given by Li and Yau .They proved that ,under the assumption $Ric(M) \geq -k$,a positive solution of the heat equation $\left(\frac{\partial}{\partial t} + \Delta\right)u = 0$, satisfies the gradient estimate $\frac{|gradu|^2}{u^2} - \frac{u_t}{t} \leq \frac{n}{2t}$ .Along the proof of this fact ,they define the function $F(x,t) = t(|gradf|^2 - f_t)$ .*

*Now we give a short review Rugang Ye works on estimates for reduced lengths ,that is important for the applications to have results on the Lipschitz properties on l or L.*

*Proposition7.4(see47) Let $\bar{\tau} \in (0,T)$.Assume that $Ric \geq -cg$ on $[0,\bar{\tau}]$ for a non-negative constant C.Then $L(.,\tau)$ is locally Lipschitz with respect to the metric $g(\tau)$ for each $\tau \in (0,\bar{\tau}]$.Moreover ,for each compact subset E of M,There are positive constants $A_1$ and $A_2$ such that $\sqrt{\tau}L \leq A_1$ on $E \times (0,\bar{\tau}]$ and*

$$|\dot{\gamma}(s)|^2 \leq \frac{A_2}{s}\left(1 + \frac{1}{\tau}\right)$$

*For $s \in (0,\tau]$,where $\tau \in (0,\bar{\tau}]$ and $\gamma$ denotes an arbitrary $\mathcal{L}_{0,\tau}$ −geodesic from p to $q \in E$ ,where we denote*

$$\mathcal{L}_{a,b}(\gamma) = \int_{\sqrt{a}}^{\sqrt{b}} \left(\frac{1}{2}|\gamma'|^2 + 2Rt^2\right) dt$$

*Proposition7.5(see[47]) assume that the Ricci curvature is bounded from below on $[0, \bar{\tau}]$. Then $L(q,.)$ is locally Lipschitz on $(0, \bar{\tau}]$ for every $q \in M$ .Moreover ,$\tau^{\frac{3}{2}}\left|\frac{\partial L}{\partial \tau}\right|$ is bounded on $E \times (0, \bar{\tau}]$ for each compact subset $E$ of $M$.*

*Proposition7.6 assume that the Ricci curvature is bounded from below on $[0, \bar{\tau}]$ .Then L is locally lipschitz function on $M \times (0, \bar{\tau}]$.*

*Also similar estimates for l hold in the case of bounded sectional curvature.*

*Proposition7.7 Assume that the sectional curvature is bounded on $[0, \bar{\tau}]$ .Then there is a positive constant $c = c(\tau^*)$ for every $\tau^* \in (0, \bar{\tau})$ with the following properties .For each $\tau \in (0, \tau^*]$ we have*

$$|\nabla l|^2 \leq \frac{c}{\tau}(l + \tau + 1)$$

*Almost everywhere in $M$ .For each $q \in M$ we have*

$$\left|\frac{\partial l}{\partial \tau}\right| \leq \frac{c}{\tau}(l + \tau + 1)$$

*Almost everywhere in $(0, \tau^*]$.*

**Perelman's reduced volume**

*The Perelman reduced volume is the fundamental tool which is used to establish non-collapsing results which in turn are essential in proving the existence of geometric limits. Note that reduced volume cannot be defined globally, but only on appropriate open subsets of a time-slice .But as long as one can flow an open set U of a time-slice along minimizing $\mathcal{L}$ –geodesics in the direction of decreasing$\bar{\tau}$, the reduced volumes of resulting family of open sets form a monotone non-increasing function of $\bar{\tau}$. This turns out to be sufficient to extend the non-increasing results to Ricci flow with surgery .Note that Perelman's reduced volume resembles the expression in*

*Huisken's monotonicity formula for the mean curvature flow. In this section we show that the reduced volume is monotonically no increasing in $\tau$ and is finite.*

*Definition 7.5 Given a solution of the backward Ricci flow, the reduced volume function is defined as*

$$\widetilde{V}(\tau) = \int_M \tau^{\frac{-n}{2}} exp(-l(q,\tau)) \, dq$$

*Being $dq$ the volume form of $g(\tau)$.*

*Notice that the integrand is the heat kernel in the Euclidean space. In next propositions we summarize the main properties of this function.*

*Remark 7.1 notice that $(4\pi)^{\frac{n}{2}} = \int_{R^n} \tau^{\frac{-n}{2}} exp(-l(q,\tau)) \, dq$, which is the same integral of the reduced function, but in the Euclidean space.*

*Proof Change of variable $s = \sqrt{\tau}, \frac{d\gamma}{d\tau} = \frac{1}{2\sqrt{\tau}} \frac{d\gamma}{ds}, d\tau = 2s\,ds$ and minimize over $\gamma$ such that $\gamma(0) = p, \gamma(\bar{s}) = q$. Then we get $L(q,\bar{s}) = \frac{1}{2}\frac{d(p,q)^2}{\bar{s}}$, $L(q,\bar{\tau}) = \frac{1}{2}\frac{d(p,q)^2}{\sqrt{\bar{\tau}}}$ and so $l(q,\bar{\tau}) = \frac{1}{4}\frac{d(p,q)^2}{\bar{\tau}}$, therefore we get $\widetilde{V}(\tau) = \int_{R^n} \bar{\tau}^{\frac{-n}{2}} e^{-\frac{|q|^2}{4\bar{\tau}}} dq = (4\pi)^{\frac{n}{2}}$ constant in $\bar{\tau}$.*

*In the case when $M$ is non-compact, it is not clear a priori that the integral defining the reduced volume is finite in general.*

*In fact, as the next Propositions shows, it is always finite and indeed, it is bounded above by the integral for $R^n$.*

*Theorem 7.9 (Perelman [39]) the reduced volume*

$$\widetilde{V}(\tau) = \int_M \tau^{\frac{-n}{2}} exp(-l(q,\tau)) \, dq$$

*is non-increasing along the backward Ricci flow.*

*Proof according to previous Propositions we know*

$$\frac{\partial l}{\partial \tau} = R - \frac{l}{\tau} + \frac{1}{2\tau^{\frac{3}{2}}}K \quad (7.21)$$

$$|\nabla l|^2 = -R + \frac{l}{\tau} - \frac{1}{\tau^{\frac{3}{2}}}K \quad (7.22)$$

*So we get*

$$\frac{\partial l}{\partial \tau} + |\nabla l|^2 = \frac{-K}{\tau^{\frac{3}{2}}} \quad (7.23)$$

*But we know* $\quad \Delta L|_{\bar{\tau}} \leq \frac{n}{\sqrt{\bar{\tau}}} - 2\sqrt{\bar{\tau}}R - \frac{1}{\bar{\tau}}K \quad (7.24)$

*So from (7.23) and (7.24) we arrive*

$$\frac{\partial l}{\partial \tau} - \Delta l + |\nabla l|^2 - R + \frac{n}{2\tau} \geq 0$$

*Let* $\phi = \tau^{\frac{-n}{2}}exp(-l)$ .*Then*

$$\phi_\tau = \left(\frac{-n}{2\tau} - \frac{\partial l}{\partial \tau}\right)\phi$$
$$\leq (-\Delta l + |\nabla l|^2 - R)\phi$$
$$= \Delta \phi - R\phi$$

*This means that we have proved the following inequality*

$$\phi_\tau - \Delta \phi + R\phi \leq 0$$

*Therefore we get*

$$\frac{d}{d\tau}\tilde{V}(\tau) = \int_M (\phi_\tau + R\phi)\,dq \leq \int_M \Delta\phi dq = 0$$

*So $\tilde{V}(\tau)$ is non-increasing.*

*Note that equality hold if and only if g is a Ricci soliton.*

*Definition7.6 the $\mathcal{L}$ −exponential map $\mathcal{L}\exp: T_P M \times R^+ \to M$ is defined as follows:$\forall X \in T_P M$ let:*

$$\mathcal{L}\exp_X(\bar{\tau}) = \gamma(\bar{\tau})$$

*Where $\gamma(\tau)$ is the $\mathcal{L}$ −geodesic, starting at p and having X as the limit of $\sqrt{\tau}\dot{\gamma}$ (τ as $\tau \to 0$).*

Denote by $\mathcal{J}(\tau)$ the jacobian of $\mathcal{L}\exp(\tau): T_P M \to M$. Finally, we have that $\widetilde{V}(\tau) = \int_{U \subset T_p M} \tau^{\frac{-n}{2}} exp(-l(\tau)) \mathcal{J}(\tau) \, dX$. Where $U = U_\tau = \{x \in M | \tau \leq \tau(x)\}$ and $(x) = \sup\{\tau | \gamma(\tau) \text{ is a minimizing geodesic}\}$.

*In next theorem we summarize the main properties of function*

*Theorem 7.10:*

1) If Ric is bounded from below on $[0, \tau]$ for each $\tau$, then $\widetilde{V}(\tau)$ is a non-increasing function.

2) If we assume that at least one of the following conditions hold

a) **Ric** is bounded on $[0, \tau]$ for all $\tau$.

b) The curvature operator is non-negative for each $\tau$

*Then $\widetilde{V}(\tau) \leq (4\pi)^{\frac{n}{2}}$ for all $\tau$*

3) If we assume that either

a) Rm is non-negative

b) The sectional curvature is bounded on $[0, \tau]$ for all $\tau$

*Then we have the following equality*

$$\widetilde{V}(\tau_2) - \widetilde{V}(\tau_1) = -\int_{\tau_1}^{\tau_2} \int_M \left(\frac{\partial l}{\partial \tau} - R + \frac{n}{2\tau}\right) e^{-l} \tau^{\frac{-n}{2}} dq d\tau.$$

# K-NONCOLLAPSING ESTIMATES

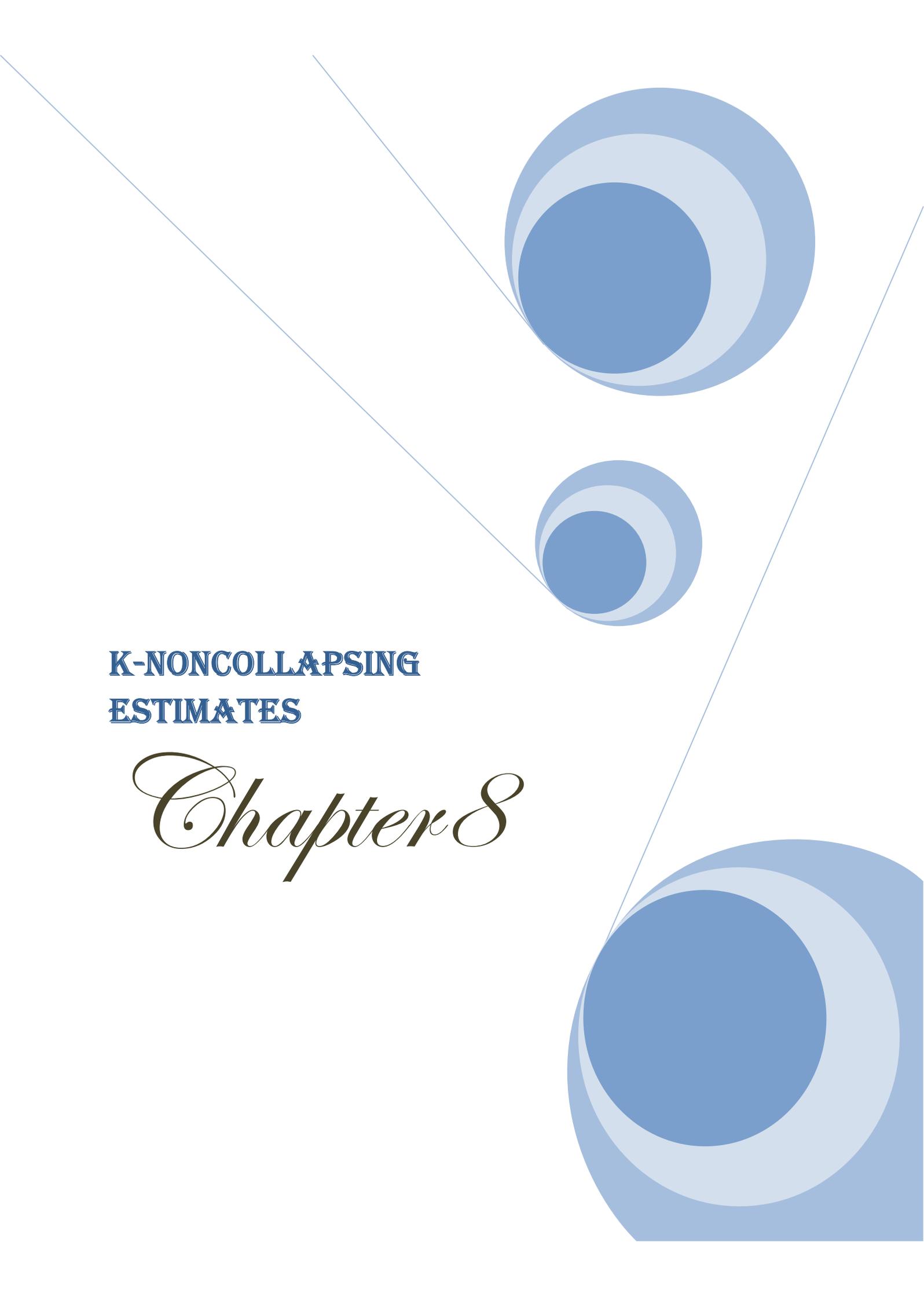

*Chapter 8*

# Chapter8

# $\kappa$ −Noncollapsing Estemates

*A family $(g_\varepsilon)$ of Riemannian metrics on a manifold M called collapses (with bounded curvature )if the injectivity radius of $g_\varepsilon$ goes to zero at each point ,whereas the sectional curvature remains bounded ,say by $1$ .Fix $\rho > 0$ (the scale) and consider all radii $0 < r < \rho$ .One says that a ball $B = B(x,r)$ is admissible if $|Rm| \leq r^{-2}$ on $B(x,r)$ .To understand this condition ,consider the metric $\tilde{g} = r^{-2}g$ .Thus $(B, \tilde{g}) = \widetilde{B}(x, 1)$ has radius 1 and curvature satisfying $|\widetilde{Rm}| = r^2|Rm| < 1$.Now given $\kappa > 0$ ,one says that $g$ is $\kappa$ −collapsed at scale $\rho$ if there exists an admissible ball $B(x,r)$ with $r \leq \rho$ ,such that*

$$r^{-n}Vol(B(x,r)) < \kappa .$$

*Now one says that $g$ is $\kappa$ −Non collapsed at scale $\rho$ if for any admissible ball $B = B(x,r)$, with $0 < r \leq \rho$ ,One has $r^{-n}Vol(B(x,r)) \geq \kappa$.*

*Remark8.1. Note that the constant $\kappa$ is smaller than the Euclidian ratio $r^{-n}Vol_{R^n}(B(x,r))$.*

*Proof. We know for an arbitrary metric ,the limit as $r \to 0$ of this the Euclidean rato is 1.*

## Chapter 8 [κ − Noncollapsing Estemates]

*Note that we will not be able to derive uniform lower bounds on the evolution of the injectivity radius in a Ricci flow .However we will observe if there is an upper bound for the curvature on some parabolic neighborhood around a point , then there is a certain lower bound on the injectivity radius on a local scale. In this section according to no local collapsing theorems of Perelman, we obtain a good relative estimate of the volume element for the Ricci flow.*

*Definition8.1 .if $(M, g(t))_{t \in [a,b)}$ is a Ricci flow, then also $\left(M, \lambda^2 g\left(\frac{1}{\lambda^2} t\right)\right)_{t \in [\lambda^2 a, \lambda^2 b)}$ is a Ricci flow (note that Ric is scale-invariant)*

$$\frac{\partial}{\partial t} \lambda^2 g\left(\frac{1}{\lambda^2} t\right) = -2Ric\left(g\left(\frac{1}{\lambda^2} t\right)\right) = -2Ric\left(\lambda^2 g\left(\frac{1}{\lambda^2} t\right)\right)$$

*Therefore, a rescaling by $\lambda$ in space and $\lambda^2$ in time is called parabolic rescaling with factor $\lambda$ .*

*Definition8.2 (parabolic ball) a parabolic ball $B(x_0, t_0, r, \tau)$ in a Ricci flow $(M, g(.))$ is a space-time product $B(x_0, t_0, r) \times [t_0, t_0 + \tau]$ (respectively $[t_0 + \tau, t_0]$ if $\tau < 0$),where $B(x_0, t_0, r)$ is the $r$−ball around $x_0$ in the $t_0$−time-slice $(M, g(t_0))$.Note that in this section we suppose that $t \in [t_0 - r^2, t_0]$ , $x \in B_{t_0}(x_0, r)$.*

*Definition8.3 .We say that a Ricci flow solution $g(.)$ defined on a time interval $[0, T)$ is $\kappa$ −Non collapsed on the scale $\rho$ if for each $r < \rho$ and all $(x_0, t_0) \in M \times [0, T)$ with $t_0 > r^2$,whenever it is true that $|Rm(x,t)| \leq r^{-2}$ for every $x \in B_{t_0}(x_0, r)$ and $t \in [t_0 - r^2, t_0]$ ,then we also have $Vol\left(B_{t_0}(x_0, r)\right) \geq \kappa r^n$.*

*Here $B_{t_0}(x_0, r)$ is the geodesic ball centered at $x_0 \in M$ and of radius $r$ with respect to the metric $g_{ij}(t_0)$ .*

*Theorem(see[39])8.1 (Perelman:No locall collapsing theorem I) Given any metric$g_{ij}$ on an n-dimensional compact manifold M. Let $g_{ij}(t)$ be the solution to the Ricci flow on $[0, T)$ ,with $T < +\infty$,Starting at $g_{ij}$ .Then there exist positive constants $\kappa$ and $\rho_0$* *such that for any $t_0 \in [0, T)$ and any point $x_0 \in M$, the solution $g_{ij}(t)$ is $\kappa$ −Non collapsed at $(x_0, t_0)$ on all scales less than $\rho_0$.*

*This theorem from Perelman says that if $|Rm| \leq r^{-2}$ on the parabolic ball $\{(x,t) | d_{t_0}(x, x_0) \leq r, t_0 - r^2 \leq t \leq t_0\}$, then the volume of the geodesic ball $B_{t_0}(x_0, r)$ is bounded from below by $\kappa r^n$.*

*Perelman obtained a stronger version of the No local collapsing theorem, where the curvature bound assumption on the parabolic ball is replaced by that on the geodesic ball $B_{t_0}(x_0, r)$ .*

*Theorem(see[39])8.2 (Perelman's strong No local collapsing theorem I) Suppose that M is a compact Riemannian manifold and $g_{ij}(t)$ , $0 \leq t \leq T < \infty$, is a solution to the Ricci flow .Then the solution $g_{ij}(t)$ is $\kappa$ −Non collapsed at $(x_0, t_0) \in M \times [0, T)$ on the scale $\gamma \in (0, \sqrt{T}]$.*

*At first for proof of this theorem we explain some remarks.*

*Let $\Omega$ be a domain (open connected) in M and $C_0^\infty(\Omega)$ be the space of real valued infinitely differentiable functions, compactly supported in $\Omega$, the Sobolev space $H_0^1(\Omega)$ is the closure of $C_0^\infty(\Omega)$ in the Norm*

$$\|f\|^2 = \int_\Omega f^2 + \int_\Omega |\nabla f|^2.$$

*Where the integrations use the volume element arising from the Riemannian structure .for any real valued measurable function f on $\Omega$ ,we say that $f \in L^{P^+}(\Omega)$ if $|f|^q$ is integrable on $\Omega$ for some $q > p$.we use $\|f\|_q$ to denote the $L^q(\Omega)$ norm of f.*

*Definition 8.4 Let H be a non-negative measurable function on $\Omega$ ,for which $\log H \in L^{\frac{n}{2}+}$ .let $\rho$ be a positive real number and define $a_\rho(H)$ as the*

$$\inf \left\{ \int (\rho |\nabla f|^2 - f^2 \log f^2 + f^2 \log H) \ for \ f \in H_0^1, \int f^2 = 1 \right\}$$

## Chapter 8 [κ − Noncollapsing Estemates]

*Definition 8.5* any $f \in H_0^1$ with $\int f^2 = 1$ which attains the minimum of $a_\rho(H)$ will be called a minimizer for $a_\rho(H)$.

*Show can consider the following results*

*Result 1* for $f \in H_0^1$ with $\int f^2 = 1$, the functional $\int (\rho|\nabla f|^2 - f^2 \log f^2)$ is bounded below

*Result 2.* $a_\rho(H)$ is finite.

*Result 3.* $a_\rho(H)$ is an attained minimum.

*Result 4.* $\mu(g_{ij}, \tau)$ is achieved by a smooth minimizer $f$ from result 3.

*Now we come back to previous to the proof of Perelman's theorem for No locall collapsing theorem I.*

*Proof.* Our aim is to prove

$$Vol_{t_0}\left(B_{t_0}(x_0, a)\right) \geq \kappa a^n, \text{ for all } 0 < a \leq r \quad (8.1)$$

*We know* $\mu(g_{ij}, \tau - t)$ *is non decreasing in $t$ where*

$$\mu(g_{ij}, t) = \inf\left\{ \mathcal{W}(g_{ij}, f, \tau) \Big| \int_M (4\pi\tau)^{\frac{-n}{2}} e^{-f} d\vartheta = 1 \right\}$$

*If we assume that $g_{ij}(t) = g_{ij}(0)$ for all $t \in R$ then $\mu(g(0), 2T) \leq \mu(g(0), \tau)$ for all $0 < \tau \leq 2T$.*

*Let $f$ be the minimizer of $\mu(g(0), 2T)$, we said that $f$ is smooth. Since M is compact, we get $|\mu(g(0), 2T)| \leq a$ for some $a \in R$. Let*

$$\mu_0 = \inf_{0 \leq \tau \leq 2T} \mu(g_{ij}(0), \tau) \geq \mu(g_{ij}(0), 2T) > -\infty$$

*according to property of $\mu$ we have*

$$\mu(g_{ij}(t_0), b) \geq \mu(g_{ij}(0), t_0 + b) \geq \mu_0 \quad (8.2)$$

*For $0 < b < r^2$. Let $0 < \xi \leq 1$ be a positive smooth function on $\mathbb{R}$ where $\xi(s) = 1$ for $|s| \leq \frac{1}{2}$, $\frac{|\dot\xi|^2}{\xi} \leq 20$, everywhere, and $\xi(s)$ is very close to zero for $|s| \geq 1$. Define a function $f$ on $M$ by*

$$(4\pi r^2)^{\frac{-n}{2}} e^{-f(x)} = e^{-c}(4\pi r^2)^{\frac{-n}{2}} \xi\left(\frac{d_{t_0}(x,x_0)}{r}\right),$$

*Where the constant $c$ is chosen so that $\int (4\pi r^2)^{\frac{-n}{2}} e^{-f} d\vartheta_{t_0} = 1$.*

*Then from (8.2) we get*

$$\mathcal{W}(g_{ij}(t_0), f, r^2) = \int_M [r^2(|\nabla f|^2 + R) + f - n](4\pi r^2)^{\frac{-n}{2}} e^{-f} d\vartheta_{t_0} \geq \mu_0. \quad (8.3)$$

*By (8.3), we obtain*

$$(c - n) + \int_M [r^2(|\nabla f|^2 + R) - \log\xi](4\pi r^2)^{\frac{-n}{2}} e^{-f} d\vartheta_{t_0} \geq \mu_0 \quad (8.4)$$

*so we see that*

$$1 = \int_M (4\pi r^2)^{\frac{-n}{2}} e^{-c} \xi\left(\frac{d_{t_0}(x,x_0)}{r}\right) dV_{t_0}$$

$$\geq \int_{B_{t_0}(x_0, \frac{r}{2})} (4\pi r^2)^{\frac{-n}{2}} e^{-c} \xi\left(\frac{d_{t_0}(x,x_0)}{r}\right) dV_{t_0}$$

$$= (4\pi r^2)^{\frac{-n}{2}} e^{-c} Vol_{t_0}\left(B_{t_0}\left(x_0, \frac{r}{2}\right)\right) \text{ (Because } \frac{d_{t_0}(x,x_0)}{r} \leq \frac{1}{2}\text{ )}$$

$$|\nabla f|^2 = |\nabla(-\log\xi)|^2 = \frac{(\dot\xi)^2}{\xi^2} \frac{1}{r^2} \quad (8.5)$$

*So from (8.4) and (8.5) we get*

$$c \geq -\int_M \left(\frac{(\dot\xi)^2}{\xi} - \log\xi.\xi\right) e^{-c}(4\pi r^2)^{\frac{-n}{2}} dV_{t_0} + (n-1) + \mu_0$$

$$\geq -2(20 + e^{-1})e^{-c}(4\pi r^2)^{\frac{-n}{2}} Vol_{t_0}\left(B_{t_0}(x_0, r)\right) + (n-1) + \mu_0$$

## Chapter 8 [κ − Noncollapsing Estemates]

$$\geq -2(20 + e^{-1}) \frac{(Vol_{t_0}(B_{t_0}(x_0, r)))}{\left(Vol_{t_0}\left(B_{t_0}\left(x_0, \frac{r}{2}\right)\right)\right)} + (n-1)_+ \mu_0$$

*Where we used the fact that $\xi(s)$ is very close to zero for $|s| \geq 1$.*

*Note also that*

$$1 = \int_M (4\pi r^2)^{\frac{-n}{2}} e^{-f} d\vartheta_{t_0}$$

$$= \int_M (4\pi r^2)^{\frac{-n}{2}} e^{-c} \xi\left(\frac{d_{t_0}(x, x_0)}{r}\right) d\vartheta_{t_0} \leq 2 \int_{B_{t_0}(x_0, r)} e^{-c} (4\pi r^2)^{\frac{-n}{2}} d\vartheta_{t_0}.$$

*Let us set*

$$\kappa = \min\left\{\frac{1}{2} \exp(-2(20 + e^{-1})3^{-n} + (n-1)_+ \mu_0), \frac{1}{2} \alpha_n\right\}$$

*Where $\alpha_n$ is the volume of the unit ball in $\mathbb{R}^n$. Then we obtain*

$$Vol_{t_0}\left(B_{t_0}(x_0, r)\right) \geq \frac{1}{2} e^c (4\pi r^2)^{\frac{n}{2}}$$

$$\geq \frac{1}{2} (4\pi)^{\frac{n}{2}} \exp(-2(20 + e^{-1})3^{-n} + (n-1)_+ \mu_0). r^n \geq \kappa r^n$$

*Provided $Vol_{t_0}(B_{t_0}\left(x_0, \frac{r}{2}\right)) \geq 3^{-n} Vol_{t_0}\left((B_{t_0}(x_0, r))\right)$*

*Note that the above argument also works for any smaller radius $a \leq r$. Thus we have proved the following assertion*

$$Vol_{t_0}\left(B_{t_0}(x_0, a)\right) \geq \kappa a^n \quad (8.6)$$

*Where $a \in (0, r]$ and $Vol_{t_0}(B_{t_0}\left(x_0, \frac{a}{2}\right)) \geq 3^{-n} Vol_{t_0}\left((B_{t_0}(x_0, a))\right)$. Now we argue by contradiction to prove the assertion (8.1) for any $a \in (0, r]$. Then by (8.6) we have*

$$Vol_{t_0}(B_{t_0}\left(x_0, \frac{a}{2}\right)) < 3^{-n} Vol_{t_0}\left((B_{t_0}(x_0, a))\right)$$

$$< 3^{-n} \kappa a^n < \kappa \left(\frac{a}{2}\right)^n$$

*This say that (8.1) for $\frac{a}{2}$ would also fail. By induction, we deduce that,*

$$Vol_{t_0}\left(B_{t_0}\left(x_0, \frac{a}{2^k}\right)\right) < \kappa \left(\frac{a}{2^k}\right)^n \quad \text{for all} \quad k \geq 1$$

*This is a contradiction since*

$$\lim_{k \to \infty} \frac{Vol_{t_0}\left(B_{t_0}\left(x_0, \frac{a}{2^k}\right)\right)}{\left(\frac{a}{2^k}\right)^n} = \alpha_n$$

## No local collapsing theorem II

*In this section we will extend the No local collapsing theorem to any complete solution with bounded curvature .In some sense, the second No local collapsing theorem gives a good relative estimate of the volume element for the Ricci flow.*

*Theorem(see[49])8.3 (Perelman) Let $g_{ij}(x,t)$ be a solution to the Ricci flow $\left(g_{ij}\right)_t = -2R_{ij}$ and $Ric(x, t_0) \leq (n-1)k$ for $dist_{t_0}(x, x_0) < r_0$ .Then the distance function $d(x,t) = dist_t(x, x_0)$ satisfies at $t = t_0$ , outside $B(x_0, r_0)$ ,the differential inequality*

$$\partial_t - \Delta d \geq -(n-1)\left(\frac{2}{3}\kappa r_0 + r_0^{-1}\right)$$

*Where $d_t = \frac{d}{dt} dist_t$.*

*Proof. We will assume that $x$ and $x_0$ are not conjugate in metrics $g(t_0)$,because otherwise the inequality that we want to prove can be understood in a barrier sence .let $\gamma(s) = exp_x(sX)$ ,for $s \in [0, L]$ ,where $X = \dot{\gamma}(0)$, be a minimal geodesic between $x$ and $x_0$ ,such that $\gamma(0) = x_0$ and $\gamma(L) = x$ .Let $(X, e_1, \dots, e_{n-1})$ be the orthonormal basis of $T_{x_0}M$ .Let*

*$E_i$ (with $1 \leq i \leq n-1$) be the parallel vector fields along $\gamma(s)$ such that $E_i(0) = e_i$ .Let $X_i(s)$ be the Jacobi fields along $\gamma(s)$ ,such that , $X_i(L) = E_i(L)$ and $X_i(0) = 0$ (they exists ,since we assumed $x$ and $x_0$ are not conjugate points ).The formula for Laplacian of the distance function*

$$\Delta dist(x, x_0) = \sum_{i=1}^{n-1} S''_{X_i}(\gamma)$$

## Chapter 8 [κ − Noncollapsing Estemates]

Where $S''_{X_i}(\gamma)$ is the second variation along $X_i$ of the length of $\gamma$, where $X_i(t)$ are Jacobi field constructed above.

$$S''_{X_k}(\gamma) = \int_0^L \left(|X'_k(s)|^2 - R(\dot\gamma, X_k, \dot\gamma, X_k)\right) ds$$

The RHS of the equality above is usually denoted by $I(X_k, X_k)$.

Define the vector fields $Y_k$ as follows:

$$Y_k = \begin{cases} \dfrac{s}{r_0} E_k(s), if\ s \in [0, r_0] \\ \\ E_k(s), if\ s \in [r_0, L] \end{cases}$$

We can notice that the vector fields $Y_k$ have the same values at the ends of $\gamma$ as a Jacobi field $X_k$ and it is the known fact that $I(X_k, X_k) \leq I(Y_k, Y_k)$. Now we can compute

$$\Delta dist(x, x_0) \leq \sum_{i=1}^{n-1} I(Y_i, Y_i)$$

$$\leq \sum_{i=1}^{n-1} \left( \int_0^{r_0} \frac{1}{r_0^2} - \frac{s^2}{r_0^2} R(X, E_k, X, E_k) + \int_{r_0}^{dist(x,x_0)} (-R(X, E_k, X, E_k)) \right)$$

$$= \int_\gamma (-Ric(X,X)) + \int_0^{r_0} Ric(X,X)\left(1 - \frac{s^2}{r_0^2}\right) + \frac{n-1}{r_0}$$

$$\leq d_t + (n-1)\left(\frac{2}{3}\kappa r_0 + r_0^{-1}\right)$$

Where we have used the fact that

$$d_t(x, x_0) = \frac{d}{dt} dist(x, x_0)\Big|_{t=t_0} = \frac{d}{dt}\int_0^L |\dot\gamma_t(s)|_{g(t)}\Big|_{t=t_0} = \int_0^L Ric(X,X) ds$$

*Where X is the unit tangent vector to the part $\gamma_t$ in metric $g(t)$. because the Ricci flow equation yields*

$$\frac{d}{dt}\left(\langle\vartheta,\vartheta\rangle_{g(t)}\right) = -2Ric(\vartheta,\vartheta).$$

*Theorem(see[48-49])8.4 (Perelman) Let $g_{ij}(x,t)$ be a solution to the Ricci flow on a $n$−dimensional manifold M and denote by $d_t(x,x_0)$, the distance between $x$ and $x_0$ with respect to the metric $g_{ij}(t)$. If $Ric(.,t_0) \leq (n-1)\kappa$ on $B_{t_0}(x_0,r_0) \cup B_{t_0}(x_1,r_0)$ for some $x_0, x_1 \in M$ and some positive constants $\kappa$ and $r_0$. then, at $t = t_0$*

$$\frac{d}{dt}d_t(x_0,x_1) \geq -2(n-1)(\frac{2}{3}\kappa r_0 + r_0^{-1})$$

*Proof. The proof of this theorem is similar to previous theorem*

*Definition8.6 :*

*Let $f:[a,b] \to \mathbb{R}$ be a continuous function on an interval. we say that the forward difference quotient of f at a point $t \in [a,b)$, denoted $\frac{df}{dt}(t)$, is less than c provided that*

$$\overline{\lim_{\Delta t \to 0^+}} \frac{f(t+\Delta t) - f(t)}{\Delta t} \leq c$$

*We say that it is greater than or equal to $c'$ if*

$$c' \leq \underline{\lim_{\Delta t \to 0^+}} \frac{f(t+\Delta t) - f(t)}{\Delta t}$$

*Lemma8.1 suppose that $f:[a,b] \to \mathbb{R}$ is a continuous function. suppose that $\psi$ is a $C^1$ −function on $[a,b] \times \mathbb{R}$ and suppose that $\frac{df}{dt}(t) \leq \psi(t,f(t))$ for every $t \in [a,b)$ in the sense of forward difference quotients. suppose also that there is a function $G(t)$ defined on $[a,b]$ that satisfies the differential equation $G'(t) = \psi(t,G(t))$ and has $f(a) \leq G(a)$. Then $f(t) \leq G(t)$ for all $\in [a,b]$.*

## Chapter 8 [κ − Noncollapsing Estemates]

*Claim8.1. suppose that for every minimal $g(t_0)$-geodesic $\gamma$ from $x_0$ to $x_1$ the function $l_t(\gamma)$ which is the $g(t)$ −length of $\gamma$ satisfies*

$$\left.\frac{d(l_t(\gamma))}{dt}\right|_{t=t_0} \geq c$$

*Then*

$$\left.\frac{d(d_t(x_0,x_1))}{dt}\right|_{t=t_0} \geq c$$

*Where, if the distance function is not differentiable at $t_0$, then the inequality in the conclusion is interpreted by replacing the derivative on the left –hand side with the Liminf of the forward difference quotients of $d_t(x_0, x_1)$ at $t_0$.*

*Theorem(see[38])8.5 (result of previous Perelman's theorem) Let $t_0 \in \mathbb{R}$ and $(M, g(t))$ be a Ricci flow defined for $t$ in an interval containing $t_0$ and with $(M, g(t))$ complete for every $t$ in this interval. Fix a constant $\kappa < \infty$. suppose that $Ric(x, t_0) \leq (n-1)\kappa$ for all $x \in M$. then for every points $x_0, x_1 \in M$ we have*

$$\left.\frac{d(d_t(x_0,x_1))}{dt}\right|_{t=t_0} \geq -4(n-1)\sqrt{\frac{2\kappa}{3}}$$

*In the sense of forward difference quotients.*

*Proof. We divide the proof for two cases*

I)     *Let $d_t(x_0, x_1) \geq \sqrt{\frac{6}{\kappa}}$*

*We take $r_0 = \sqrt{\frac{3}{2\kappa}}$ in theorem8.4, and we conclude that the liminf at $t_0$ of the difference quotients for $d_t(x_0, x_1)$ is at most $-4(n-1)\sqrt{\frac{3}{2\kappa}}$.*

II)  Let $d_t(x_0, x_1) < \sqrt{\frac{6}{\kappa}}$

*Let $\gamma(u)$ be any minimal $g(t_0)$-geodesic from $x_0$ to $x_1$ parametrized by arc length. Since*

$$\left.\frac{d(l_t(\gamma))}{dt}\right|_{t=t_0} = -\int_\gamma Ric_{g(t_0)}(\gamma'(u), \gamma'(u))du$$

*We see that*

$$\left.\frac{d(l_t(\gamma))}{dt}\right|_{t=t_0} \geq -(n-1)\kappa\sqrt{\frac{6}{\kappa}} = -(n-1)\sqrt{6\kappa}$$

*But according to claim8.1 we get the liminf of the forward difference quotient of $d_t(x_0, x_1)$ at $t = t_0$ is at least $-(n-1)\sqrt{6k} \geq -4(n-1)\sqrt{\frac{2\kappa}{3}}$.*

*Result. Let $(M, g(t))$, $a \leq t \leq b$ be a Ricci flow with $(M, g(t))$ complete for every $t \in [0, T)$. Fix a positive function $\kappa(t)$, and suppose that $Ric_{g(t)}(x, t) \leq (n-1)\kappa(t)$ for all $x \in M$ and all $t \in [a, b]$. Let $x_0$ and $x_1$ be two points of $M$. Then*

$$d_a(x_0, x_1) \leq d_b(x_0, x_1) + 4(n-1)\int_a^b \sqrt{\frac{2\kappa(t)}{3}}dt$$

*Proof. By applying previous theorem we have*

$$\left.\frac{d(d_t(x_0, x_1))}{dt}\right|_{t=t'} \geq -4(n-1)\sqrt{\frac{2\kappa(t')}{3}}$$

*In the sense of forward difference quotients. thus this result is an immediate consequence of lemma8.1*

# Chapter 8 [κ − Noncollapsing Estemates]

*Theorem(see [39])8.6 (No local collapsing theorem II) for any $A > 0$ there exists $\kappa = \kappa(A) > 0$ With the following property .If $g_{ij}(t)$ is a complete solution to the Ricci flow on $0 \leq t \leq r_0^2$ with bounded curvature and satisfying*

$$|Rm|(x,t) \leq r_0^2 \quad On \; B_0(x_0, r_0) \times [0, r_0^2]$$

*and*

$$Vol_0(B_0(x_0, r_0)) \geq A^{-1} r_0^n,$$

*The $g_{ij}(t)$ is $\kappa$−noncollapsed on all scales less than $r_0$ at every point $(x, r_0^2)$ with $d_{r_0^2}(x, x_0) \leq A r_0$.*

# Biography